\documentclass[a4paper,10pt]{article}
\usepackage{a4,amsxtra,amsmath,amsfonts,amscd, amssymb, mathrsfs}
\usepackage[all]{xy}
\newtheorem{thm}{Theorem}[section]
\newtheorem{prop}[thm]{Proposition}
\newtheorem{cor}[thm]{Corollary}
\newtheorem{lem}[thm]{Lemma}
\newtheorem{ex}[thm]{Example}

\newtheorem{rem}[thm]{Remark}
\newtheorem{art}[thm]{}

\newcommand{\codim}{{\rm codim}}

\newcommand{\Div}{{\rm div}}

\newcommand{\Pic}{{\rm Pic}}

\newcommand{\Spec}{{\rm Spec}}
\newcommand{\Spf}{{\rm Spf}}

\newcommand{\id}{{\rm id}}

\newcommand{\Gal}{{\rm Gal}}
\newcommand{\supp}{{\rm supp}}

\newcommand{\Acal}{{\mathscr A}}
\newcommand{\Amcal}{{\mathcal A}}
\newcommand{\Bcal}{{\mathscr B}}
\newcommand{\Bmcal}{{\mathcal B}}
\newcommand{\Ccal}{{\mathscr C}}
\newcommand{\Cmcal}{{\mathcal C}}
\newcommand{\Dcal}{{\mathscr D}}
\newcommand{\Dmcal}{{\mathcal D}}
\newcommand{\Ecal}{{\mathscr E}}
\newcommand{\Emcal}{{\mathcal E}}

\newcommand{\Gcal}{{\mathscr G}}
\newcommand{\Hcal}{{\mathscr H}}

\newcommand{\Lcal}{{\mathscr L}}
\newcommand{\Lmcal}{{\mathcal L}}
\newcommand{\Mcal}{{\mathscr M}}

\newcommand{\Ocal}{{\mathscr O}}

\newcommand{\Scal}{{\mathscr S}}

\newcommand{\Ucal}{{\mathscr U}}
\newcommand{\Vcal}{{\mathscr V}}
\newcommand{\Xcal}{{\mathscr X}}
\newcommand{\Xmcal}{{\mathcal X}}
\newcommand{\Ycal}{{\mathscr Y}}
\newcommand{\Zcal}{{\mathscr Z}}
\newcommand{\Zmcal}{{\mathcal Z}}

\newcommand{\qdop}{{\mathbb Q}}
\newcommand{\ndop}{{\mathbb N}}
\newcommand{\rdop}{{\mathbb R}}

\newcommand{\kdop}{{\mathbb K}}
\newcommand{\ldop}{{\mathbb L}}

\newcommand{\adop}{{\mathbb A}}
\newcommand{\zdop}{{\mathbb Z}}

\newcommand{\tdop}{{\mathbb T}}

\newcommand{\metr}{{\|\hspace{1ex}\|}}
\newcommand{\canmetr}{{\|\hspace{1ex}\|_{\rm can}}}

\newcommand{\X}{{\frak X}}
\newcommand{\Y}{{\frak Y}}

\newcommand{\Tfrak}{{\frak T}}

\newcommand{\ghp}{{\hat{\frak g}_{X}^+}}

\newcommand{\proof}{\noindent {\bf Proof: \/}}
\newcommand{\qed}{{ \hfill $\square$}}

\newcommand{\val}{{\rm val}}
\newcommand{\Val}{{\rm Val}}
\newcommand{\vol}{{\rm  vol}}
\newcommand{\ub}{{\mathbf u}}
\newcommand{\ubb}{{\overline{\mathbf u}}}
\newcommand{\wb}{{\mathbf w}}

\newcommand{\ab}{{\mathbf a}}
\newcommand{\bb}{{\mathbf b}}
\newcommand{\cb}{{\mathbf c}}
\newcommand{\mb}{{\mathbf m}}
\newcommand{\nb}{{\mathbf n}}

\newcommand{\xb}{{\mathbf x}}
\newcommand{\yb}{{\mathbf y}}
\newcommand{\zb}{{\mathbf z}}
\newcommand{\Tor}{{\mathbb G}_m^n}
\newcommand{\rtor}{{\rdop^n/\Lambda}}

\newcommand{\Deltabar}{{\overline{\Delta}}}
\newcommand{\Ccalbar}{{\overline{\Ccal}}}

\newcommand{\Sigmabar}{{\overline{\Sigma}}}
\newcommand{\sigmabar}{{\overline{\sigma}}}
\newcommand{\valbar}{{\overline{\val}}}
\newcommand{\Xan}{{X^{\rm an}}}

\newcommand{\relint}{{\rm relint}}

\newcommand{\str}{{\rm str}}
\newcommand{\kcirc}{{\rm \kdop^\circ}}
\newcommand{\ktilde}{{\rm \tilde{\kdop}}}
\newcommand{\xskel}{{S(\Xcal')}}
\newcommand{\xstr}{{\str(\tilde{\Xcal}')}}
\newcommand{\Max}{{\rm Max}}
\newcommand{\Xfrak}{{\frak X}}
\title{Non-archimedean canonical measures\\ on abelian varieties}
\author{Walter Gubler}
\date{\today}
\begin{document}

\maketitle

\begin{abstract}
For a closed $d$-dimensional subvariety $X$ of an abelian variety $A$ and  a canonically metrized line bundle $L$ on $A$, Chambert-Loir has introduced measures $c_1(L|_X)^{\wedge d}$ on the Berkovich analytic space associated to $A$ with respect to the discrete valuation of the ground field. In this paper, we give an explicit description of these canonical measures in terms of convex geometry. We use a generalization of the tropicalization related to the Raynaud extension of $A$ and Mumford's construction. The results have applications to the equidistribution of small points. 
\end{abstract}

\section{Introduction}

Let $K$ be a field with a discrete valuation $v$, valuation ring $K^\circ$ and residue field $\tilde{K}$.  We denote the completion of the algebraic closure of the completion of $K$ by $\kdop$. This algebraically closed complete field is used for analytic considerations on algebraic varieties defined over $K$. For the analytic facts, we refer to \S 2.

In non-archimedean analysis, there is no analogue known for the first Chern form of a metrized line bundle. However, Chambert-Loir \cite{Ch} has introduced measures $c_1(\overline{L_1}) \wedge \dots \wedge c_1(\overline{L_d})$ on the Berkovich analytic space $\Xan$ associated to a $d$-dimensional projective variety $X$. The analogy to the corresponding  forms in differential geometry comes from Arakelov geometry. These measures are best understood in case of metrics induced by line bundles $\Lcal_1, \dots, \Lcal_d$ on a projective $K^\circ$-model $\Xcal$ of $X$, with generic fibres $L_1, \dots ,L_d$. In this standard situation from Arakelov geometry, $c_1(\overline{L_1}) \wedge \dots \wedge c_1(\overline{L_d})$ is a discrete measure on $\Xan$ with support  and multiplicities determined by the irreducible components of $\tilde{\Xcal}$ and their degrees with respect to $\Lcal_1, \dots, \Lcal_d$. However, the canonical metric on an ample line bundle of an abelian 
  variety $A$ over $K$ is  given by such models only if $A$ has potential good reduction. In general, a variation of Tate's limit argument shows that the canonical metric is a uniform limit of roots of model metrics and hence the corresponding  canonical measure is given as a limit of discrete measures. We recall the theory of Chambert-Loir's measures in  \S 3. 

We consider an irreducible $d$-dimensional closed subvariety $X$  of the abelian variety $A$. Using the Raynaud extension of $A$, there is a complete lattice $\Lambda$ in $\rdop^n$ and a map $\valbar:A^{\rm an} \rightarrow \rtor$, where $n$ is the torus rank of $A$. We call  $\valbar(\Xan)$ the tropical variety associated to $X$. This analytic analogue of tropical algebraic geometry is described in \S4. Let $b$ be the dimension of the abelian part  of good reduction in the Raynaud extension of $A$ and hence $\dim(A)=b+n$. For a simplex $\Delta$ in $\rdop^n$, we denote by $\delta_\Deltabar$ the Dirac measure in $\Deltabar$, i.e. the pushforward of the Lebesgue measure on $\Delta$ to $\rtor$. The main result of this paper is the following explicit description of  canonical measures in terms of convex geometry:

\begin{thm} \label{main theorem}
There are rational simplices $\overline{\Delta}_1, \dots \overline{\Delta}_{N}$ in $\rtor$ with the following properties:
\begin{itemize}
\item[(a)] For $j=1, \dots,N$, we have $\dim(\overline{\Delta}_j)\in \{d-b, \dots, d\}$.
\item[(b)] $\valbar(\Xan)= \bigcup_{j=1}^N \Deltabar_j$.
\item[(c)] For canonically metrized line bundles $\overline{L_1}, \dots, \overline{L_d}$ on $A$, there are $r_j \in \rdop$ with
$$\valbar_*\left(c_1(\overline{L_1}|_X) \wedge \dots \wedge c_1(\overline{L_d}|_X) \right)
= \sum_{j=1}^{N} r_j \cdot \delta_{\overline{\Delta}_j}.$$
\item[(d)] If all line bundles in (c) are ample, then $r_j>0$ for $j \in \{1, \dots , N\}$.
\end{itemize}
\end{thm}
{\it Erratum: }  In the preprint version \cite{Gu7} of this paper, it was claimed that the tropical variety $\valbar(\Xan)$ is of pure dimension. However the referee has found a gap in the argument (see \ref{wrong argument}) and so this question remains open. As a consequence, in  Theorem 1.2 of \cite{Gu6}, one should omit to claim that  the tropical variety is of pure dimension. All other claims remain valid.
\vspace{3mm}

Theorem \ref{main theorem} was proved in \cite{Gu4} for abelian varieties which are totally degenerate at the place $v$. This special case  is equivalent to $b=0$ which makes the arguments easier. In particular, the tropical variety $\valbar(\Xan)$ is of pure dimension $d$. In the general case, we can still show in Theorem \ref{dimension of tropical variety} that the tropical variety $\valbar(\Xan)$ is a polytopal set with the above properties (a) and (b). The most serious problem is that the tropical dimension may be strictly smaller than $d$. This leads to the unpleasant fact that the canonical measure  in Theorem \ref{main theorem}(c) may have singular parts in lower dimensions which is in sharp contrast to the totally degenerate case. 

Using a semistable alteration, we will give in \S 6 an explicit description of  the canonical measure $c_1(\overline{L_1}|_X) \wedge \dots \wedge c_1(\overline{L_d}|_X) $ on $\Xan$ in terms of convex geometry. It relies on our study of  Mumford models of $A$ in \S 4 and on the properties of the skeleton of the strictly semistable model from the alteration given in \S 5. A Mumford model is associated to a rational  $\Lambda$-periodic polytopal decomposition of $\rdop^n$ such that the reduction modulo $v$ brings  toric varieties and convex geometry into play.  In Theorem \ref{piecewise linear}, we show that the support of this canonical measure is a canonical subset of $\Xan$ which does not depend on the choice of the ample line bundles $L_j$ and which has a canonical piecewise linear structure. Finally, the proof of Theorem \ref{main theorem} will be finished in \S 7 and we will show in 2 examples how these canonical measures can look like. In the appendix, we study building blocks of strongly non-degenerate strictly pluristable formal models. This is the background for the generalization of our results in \S 5 to such models which is required only in the proof of Theorem \ref{piecewise linear}.

Theorem \ref{main theorem} has the following application to diophantine geometry. Let $K$ be either a number field or the function field of an irreducible projective variety $B$ of positive dimension over a field $k$. In the latter case, we assume that $B$ is regular in codimension $1$ and we count the prime divisors $v$ of $B$ with multiplicity $\deg_{\mathbf c}(v)$ for a fixed ample class $\mathbf c$ on $B$. In any case, $K$ satisfies the product formula and hence we get absolute heights on projective varieties over $K$ (see \cite{BG}). In particular, we have the N\'eron--Tate height $\hat{h}$ on the abelian variety $A$ with respect to a fixed ample symmetric line bundle $L$. Note that $\hat{h}$ is a positive semi-definite quadratic form on $A(\overline{K})$ and hence defines a semi-distance. By Arakelov geometry, there is an extension of the N\'eron--Tate height to all closed subvarieties of $A$ defined over $\overline{K}$ (see \cite{Gu3}).

Let $X$ be an irreducible $d$-dimensional closed subvariety of the abelian variety $A$ over $K$. We choose a small generic net $(P_m)_{m \in I}$ in $X(\overline{K})$. Here, small means 
$$\lim_{m} \hat{h}(P_m)=\frac{\hat{h}(X)}{(d+1)\deg_L(X)}$$ 
and generic means that for every proper closed subset $Y$ of $X$, there is $m_0 \in I$ such that $P_m \not \in Y$ for all $m \geq m_0$. The absolute Galois group $G:=\Gal(\overline{K}/K)$ acts on $X(\overline{K})$ and $O(P_m)$ denotes the orbit of $P_m$.

We fix a discrete valuation $v$ on $K$ and we form $X_v^{\rm an}$ and the tropical variety $\valbar(X_v^{\rm an})$ with respect to $v$. We fix an embedding  $\overline{K} \hookrightarrow \kdop_v$ to identify $A(\overline{K})$ with a subset of $A_v^{\rm an}$. On $\valbar({X_v^{\rm an}})$,  we consider the  discrete probability measures 
$$\nu_m:= \frac{1}{|O(P_m)|} \cdot \sum_{P_m^\sigma \in O(P_m)} \delta_{\valbar(P_m^\sigma)}.$$ 

\vspace{3mm}
\noindent {\bf Tropical equidistribution theorem}{\it \: There is a regular probability measure $\nu$ on $\rdop^n$ with support equal to the tropical variety $\valbar({X_v^{\rm an}})$ such that $\nu_m\stackrel{w}{\to} \nu$ as a weak limit of Borel measures. More precisely, if we endow $L$ with a canonical metric $\metr_v$, then we have $\nu=\deg_L(X)^{-1}\valbar_*(c_1((L|_X, \metr_v))^{\wedge d})$.}

\vspace{3mm}
Note that this statement is only useful if we have the positivity of $\nu$ from the explicit description in Theorem \ref{main theorem}. The tropical equidistribution theorem follows from the equidistribution theorem 
\begin{equation} \label{equidistribution theorem}
\frac{1}{|O(P_m)|} \cdot \sum_{P_m^\sigma \in O(P_m)} \delta_{P_m^\sigma} \stackrel{w}{\rightarrow} 
\frac{1}{\deg_L(X)}\cdot c_1(L|_X,\metr_{v})^{\wedge d}
\end{equation}
on ${X_v^{\rm an}}$. For an archimedean place $v$ of a number field $K$ and for a metrized ample line bundle with positive curvature on a smooth projective variety, equidistribution \eqref{equidistribution theorem} was proved by Szpiro--Ullmo--Zhang \cite{SUZ}. This was generalized by Yuan (\cite{Yu}, Theorem 5.1) to semipositively metrized ample line bundles on  projective varieties over a number field and also to non-archimedean places. In \cite{Gu6}, Theorem 1.1, Yuan's generalization was proved in the function field case.

The potential applications of the tropical equidistribution theorem are related to the Bogomolov conjecture. The latter claims that the N\'eron--Tate height has a positive lower bound on $X(\overline{K})$ outside an explicit exceptional set. In the number field case, the Bogomolov conjecture was proved by Ullmo \cite{Ul} for curves and by Zhang \cite{Zh} in general. The main tool was the archimedean version of \eqref{equidistribution theorem}. For function fields, the Bogomolov conjecture is still open. In \cite{Gu5}, it was proved for abelian varieties which are totally degenerate with respect to a place $v$.  The proof relied on the tropical equidistribution theorem for totally degenerate abelian varieties (\cite{Gu5}, Theorem 5.5). For an arbitrary abelian variety, it is clear that the tropical equidistribution theorem can not imply the Bogomolov conjecture since the dimension of the tropical variety may decrease. However, it is plausible that it can be used once the case 
  of abelian varieties with everywhere good reduction is understood. 

\vspace{3mm}

\centerline{\it Terminology}

In $A \subset B$, $A$ may be  equal to $B$. The complement of $A$ in $B$ is denoted by $B \setminus A$ \label{setminus}. The zero is included in $\ndop$ and in $\rdop_+$.

All occuring rings and algebras are commutative with $1$. If $A$ is such a ring, then the group of multiplicative units is denoted by $A^\times$. A variety over a field is a separated reduced scheme of finite type over that field. However,  a formal analytic variety is not necessarily reduced (see \S 2). For the degree of a map $f:X \rightarrow Y$ of irreducible varieties, we use either $\deg(f)$ or $[X:Y]$. 

Let $Y$ be a variety over a field. Following \cite{Ber4}, \S 2, we use the following canonical stratification of $Y$. We start with $Y^{(0)}:=Y$. For $r \in \ndop$, let $Y^{(r+1)} \subset Y^{(r)}$ be the complement of the set of normal points in $Y^{(r)}$. Since the set of normal points is open and dense, we get a chain of closed subsets:
$$Y=Y^{(0)} \supsetneq Y^{(1)} \supsetneq Y^{(2)} \supsetneq \dots \supsetneq Y^{(s)} \supsetneq Y^{(s+1)} = \emptyset.$$
The irreducible components of $Y^{(r)} \setminus Y^{(r+1)}$ are called the {\it strata} of $Y$. The set of strata is denoted by $\str(Y)$. It is partially ordered by $S \leq T$ if and only if $\overline{S} \subset \overline{T}$. A {\it strata subset} is a union of stratas. A {\it strata cycle} is a cycle whose components are strata subsets.

For $\mb \in \zdop^n$, let $\xb^\mb:=x_1^{m_1} \cdots x_n^{m_n}$ and $|\mb|:=m_1+\dots+m_n$. The standard scalar product of $\ub,\ub' \in \rdop^n$ is denoted by  $\ub \cdot \ub':=u_1 u_1' + \dots + u_n u_n'$. For the notation used from convex geometry, we refer to \ref{convex geometry}.

\vspace{3mm}
\small
The author thanks Vladimir Berkovich for answering a question related to the proof of Theorem \ref{piecewise linear} and the referee for his precious comments and suggestions.
\normalsize

\section{Analytic and formal geometry}

In this section, we recall results  from Berkovich analytic spaces and formal geometry. 
Our base field $\kdop$ is  algebraically closed  with a non-trivial non-archimedean complete absolute value $|\phantom{a}|$,  valuation ring $\kdop^\circ$ and residue field  $\tilde{\kdop}$. 

\begin{art} \rm \label{affinoid algebras}
The {\it Tate algebra} $\kdop \langle x_1, \dots, x_n \rangle$ is the completion of  
$\kdop[x_1, \dots, x_n]$ with respect to the Gauss norm. Its elements are the power series in the variables $x_1, \dots,x_n$ with coefficients $a_\mb \in K$ such that $|a_\mb| \to 0$ for $m_1+ \dots+m_n \to \infty$. A {\it $\kdop$-affinoid algebra} $\Acal$ is isomorphic to $ \kdop \langle x_1, \dots, x_n \rangle/I$ for some ideal $I$ in $\kdop \langle x_1, \dots, x_n \rangle$. The maximal spectrum $\Max(\Acal)$ corresponds to the zero set  of  $I$ in the closed unit ball ${\mathbb B}^n:=\{\mathbf x \in \kdop^n \mid \max_j |x_j| \leq 1\}$. The supremum semi-norm of $\Acal$ on $\Max(\Acal)$ is denoted by $|\phantom{a}|_{\rm sup}$. We define
$$\Acal^\circ:= \{a \in \Acal \mid |a|_{\rm sup} \leq 1 \}, 
\quad \Acal^{\circ \circ}:=\{a \in \Acal \mid |a|_{\rm sup} < 1 \}$$
and the finitely generated reduced $\tilde{\kdop}$-algebra $\tilde{\Acal}:= \Acal^\circ / \Acal^{\circ \circ}$ (see \cite{BGR}). 
\end{art}

\begin{art} \rm \label{Berkovich spectrum} 
For an affinoid algebra $\Acal$, the {\it Berkovich spectrum} $X=\Mcal(\Acal)$ is  the set of semi-norms $p$ on $\Acal$ with
$p(ab)=p(a)p(b)$, $p(1)=1$ and $p(a) \leq |a|_{\rm sup}$ 
for all $a,b \in \Acal$. We use the coarsest topology on $X$ such that the maps $p \mapsto p(a)$ are continuous for all $a \in \Acal$. 
The affine $\tilde{\kdop}$-variety $\tilde{X}=\Spec(\tilde{\Acal})$ is called the {\it reduction} of $X$.  The reduction map $X \rightarrow \tilde{X}$, $p \mapsto \tilde{p}:= \{p<1\}/\Acal^{\circ \circ}$,  is surjective. If $Y$ is an irreducible component of $\tilde{X}$, then there is a unique $\xi_Y \in X$ with $\tilde{\xi}_Y$ equal to the generic point of $Y$. For details, we refer to \cite{Ber}. Note that our definition of an affinoid algebra is the same as in \cite{BGR}, but Berkovich calls them strictly affinoid algebras.
\end{art}

\begin{art} \rm \label{subdomains}
A {\it rational subdomain} of $X:=\Mcal(\Acal)$ is defined by 
$$X \left(\frac{\mathbf f}{g}\right) := \{x \in X \mid |f_j(x)| \leq |g(x)|, \, j=1,\dots,r\},$$
where $g,f_1, \dots, f_r \in \Acal$ are without common zero. It is the Berkovich spectrum of the affinoid algebra
$$\Acal \left\langle \frac{\mathbf f}{g}\right\rangle := \kdop \langle \mathbf x, y_1, \dots, y_r \rangle / \langle I, g(\mathbf x) y_j - f_j \mid j=1, \dots ,r \rangle,$$
where we use the description $\Acal=\kdop \langle \mathbf x \rangle / I$ from \ref{affinoid algebras}. 

More generally, one defines an {\it affinoid subdomain} of $X$ as the Berkovich spectrum of an affinoid algebra characterized by a certain universal property (see \cite{BGR}, 7.2.2). By a theorem of Gerritzen and Grauert, an affinoid subdomain is a 
finite union of rational domains. For more details, we refer to \cite{BGR}, Chapter 7, and \cite{Ber}, \S 2.2.
\end{art}

\begin{art} \rm \label{analytic spaces}
An {\it analytic space} $X$ over $\kdop$ is given by an atlas of affinoid subdomains $U=\Mcal(\Acal)$. For the precise definition, we refer to \cite{Ber2}, \S 1 (where they are called strictly analytic spaces).
\end{art}

\begin{art} \rm \label{formal analytic varieties}
The {\it formal topology} on $X=\Mcal(\Acal)$ is given by the preimages of the open subsets of $\tilde{X}$ with respect to the reduction map. This quasi-compact topology is generated by affinoid subdomains and hence we get a canonical ringed space called  a {\it formal affinoid variety} over $\kdop$ which we denote by $\Spf(\Acal)$. By definition, a morphism of affinoid varieties over $\kdop$ is induced by a reverse homomorphism of the corresponding  $\kdop$-affinoid algebras (see \cite{Bo} for details). 

A {\it formal analytic variety} over $\kdop$ is a $\kdop$-ringed space $\Xfrak$ which has a locally finite open atlas of formal affinoid varieties $\Spf(\Acal_i)$ over $\kdop$ called a {\it formal affinoid atlas}.  The {\it generic fibre} $\Xfrak^{\rm an}$ (resp. the {\it reduction} $\tilde{\Xfrak}$) is locally given by $\Mcal(\Acal_i)$ (resp. $\Spec(\tilde{\Acal}_i)$). By \ref{Berkovich spectrum}, we get a surjective reduction map $\Xfrak^{\rm an} \rightarrow \tilde{\Xfrak}$,  $x \mapsto \tilde{x}$. For every irreducible component $Y$, there is a unique $\xi_Y \in \Xfrak^{\rm an}$ such that $\tilde{\xi}_Y$ is the generic point of $Y$. 
\end{art}

\begin{art} \rm \label{admissible formal schemes}
An  {\it admissible $\kcirc$-algebra} is a $\kcirc$-algebra $A$ without $\kcirc$-torsion which is isomorphic to $\kdop^\circ\langle x_1, \dots , x_n \rangle / I$ for an ideal $I$. An {\it admissible formal scheme} $\Xcal$ over $\kdop^\circ$ is a formal scheme over $\kcirc$ which has a locally finite atlas of open subsets $\Spf(A_i)$  for  admissible $\kdop^\circ$-algebras $A_i$  (see \cite{BL3}, \cite{BL4} for details). 

The {\it special fibre} $\tilde{\Xcal}$ of $\Xcal$  is a scheme  over $\tilde{\kdop}$ locally given by $\Spec(\tilde{A}_i)$. It is locally of finite type over $\ktilde$ and  not necessarily reduced. The latter is in sharp contrast to the reduction of  formal analytic varieties. These categories are related by the following functors:

The formal analytic variety $\Xcal^{\rm f-an}$ associated to $\Xcal$ is locally given by $\Spf(\Acal_i)$ for the affinoid algebra $\Acal_i:=A_i \otimes_\kcirc \kdop$.  The canonical morphism $(\Xcal^{\rm f-an})\sptilde  \rightarrow \tilde{\Xcal}$ is finite and surjective (see \cite{BL1}, \S 1). 

The  generic fibre of $\Xcal^{\rm f-an}$  is also called the {\it generic fibre} of $\Xcal$ and is denoted by $\Xcal^{\rm an}$. Note that $\Xcal^{\rm f-an}$ and $\Xcal^{\rm an}$ are equal as a set but  $\Xcal^{\rm an}$ has a finer topology. Using composition of the reduction map for $\Xcal^{\rm f-an}$ (see \ref{formal analytic varieties}) with  the canonical morphism above, we get a surjective {\it reduction map} $\pi:\Xcal^{\rm an} \rightarrow \tilde{\Xcal}$. 

If $\Xfrak$ is a formal analytic variety over $\kdop$ given by the formal affinoid atlas $\Spf(\Acal_i)$, then the {\it associated formal scheme} $\X^{\rm f-sch}$ is locally given by  $\Spf(\Acal_i^\circ)$. 

It is often useful to flip between formal analytic varieties over $\kdop$ and admissible formal schemes over $\kcirc$. This is possible because the functors $\Xcal \rightarrow \Xcal^{\rm f-an}$ and $\Xfrak \rightarrow \Xfrak^{\rm f-sch}$ give an equivalence between the category of admissible formal schemes over $\kdop^\circ$ with reduced special fibre and the category of reduced formal analytic varieties over $\kdop$. Moreover, the canonical map $(\Xcal^{\rm f-an})\sptilde  \rightarrow \tilde{\Xcal}$ is an isomorphism. For details, see \cite{BL1}, \S 1, and \cite{Gu1}, \S1. 
\end{art}

\begin{art} \rm \label{schemes}
For a scheme $X$   of finite type over a subfield $K$ of $\kdop$, there is an   analytic space $X^{\rm an}$ over $\kdop$ associated to $X$. The construction is similar as for complex varieties. Moreover, most algebraic properties hold also analytically and conversely, there is  a GAGA principle. For details, we refer to \cite{Ber}, 3.4.

If  $\Xmcal$ is a a flat scheme of finite type over the valuation ring $K^\circ$ with generic fibre $X$, then the associated formal scheme   $\hat{\Xmcal}$ over  $\kdop^\circ$ is obtained by the $\pi$-adic completion of $X$ for any $\pi \in K$ with $|\pi|<1$. The special fibre of $\hat{\Xmcal}$ is the base change of the special fibre of $\Xmcal$ to $\tilde{\kdop}$. The generic fibre  $\hat{\Xmcal}^{\rm an}$ is an analytic subdomain of $X^{\rm an}$ such that $\hat{\Xmcal}^{an}(\kdop)=\Xmcal(\kcirc)$. If $\Xmcal$ is proper over $K^\circ$, then $\hat{\Xmcal}^{\rm an}=X^{\rm an}$. For details, we refer to \cite{Gu4}, 2.7. 
\vspace{3mm}

For convenience of the reader, we summarize here the {\it notational policy} of the paper: $\Xmcal$ denotes a flat algebraic scheme over $\kcirc$, $\Xcal$ is used for an admissible formal scheme over $\kcirc$ and $\X$ denotes a formal analytic variety over $\kdop$. The generic fibre in either case is usually denoted by $X$.
\end{art}

\section{Chambert-Loir's measures}

In this section, $\kdop$ denotes an algebraically closed field which is complete with respect to the non-trivial non-archimedean absolute value $|\phantom{a}|$.  Let $K$ be a subfield of $\kdop$ such that the valuation $v:= - \log |\phantom{a}|$ restricts to a discrete valuation on $K$. We will assume, as in our applications later on, that varieties are defined over $K$ and we will perform analytic considerations over $\kdop$ using \ref{schemes}.

First, we will characterize admissible metrics on a line bundle by their fundamental properties. As in Arakelov geometry, metrics associated to $\kcirc$-models are admissible and we want to include also canonical metrics on an abelian variety. Then we will recall the basic properties of Chambert-Loir's measures with respect to line bundles endowed with admissible metrics. These analogues of top-dimensional wedge products of first Chern forms were introduced in \cite{Ch} and later generalized in \cite{Gu4}. 

\begin{art} \label{metrics} \rm 
We recall some facts about metrics on line bundles. Let $X$ be a proper scheme over ${K}$ and let $L$ be a line bundle on $X$. We consider metrics $\metr$, $\metr'$ on $L$ which are continuous with respect to the analytic topology on $L^{\rm an}$. Then we have the {\it distance of uniform convergence}
$$d(\metr, \metr'):= \sup_{x \in \Xan} \left|\log\left(\|s(x)\|/\|s(x)\|'\right)\right|.$$
The definition is independent of the choice of a non-zero $s(x)\in L_x$.
\end{art} 

\begin{art} \label{models} \rm
A {\it formal $\kdop^\circ$-model} of $X$ is an admissible formal scheme over $\kcirc$   together with an isomorphism $\Xcal^{\rm an} \cong \Xan$. 
A {\it formal $\kdop^\circ$-model of $L$ on $\Xcal$} is a line bundle $\Lcal$ on  $\Xcal$ together with an isomorphism $\Lcal^{\rm an} \cong L^{\rm an}$.

For notational simplicity, we usally ignore the isomorphism between the generic fibre $\Xcal^{\rm an}$ and $\Xan$. We simply identify  $\Xcal^{\rm an}$ with $\Xan$.

An {\it algebraic model} $\Xmcal$ of $X$ over the discrete valuation ring $K^\circ$ of $K$ is a  scheme $\Xmcal$ which is flat and proper over $K^\circ$ and which has generic fibre (isomorphic to) $X$. We will also use formal $\kcirc$-models for analytic spaces and line bundles in the same sense as above.
\end{art}

\begin{ex} \rm \label{formal models}
If $\Lcal$ is a formal $\kcirc$-model of $L$ on $\Xcal$, then the {\it associated  formal metric}  $\metr_\Lcal$ on $L$ is defined in the following way: If $\Ucal$ is a formal trivialization of $\Lcal$ and if $s \in \Gamma(\Ucal,\Lcal)$ is corresponding to $\gamma \in \Ocal_\Xcal(\Ucal)$, then 
$\|s(x)\|=|\gamma(x)|$ on $\Ucal^{\rm an}$. Obviously, $\metr_{\Lcal}$ is continuous and independent of the choice of the trivialization. By \ref{schemes}, every algebraic model of $L$ over $K^\circ$ induces a formal $\kcirc$-model and hence an associated formal metric.
\end{ex}


\begin{prop} \label{ghp-candidates}
For every line bundle $L$ on a proper scheme $X$ over $K$, there is a set $\ghp(L)$ of continuous metrics on $L^{\rm an}$ with the following properties:
\begin{itemize}
\item[(a)] If $\metr_i$ is a $\ghp(L_i)$-metric for $i=1,2$, then $\metr_1 \otimes \metr_2$ is a $\ghp(L_1 \otimes L_2)$-metric.
\item[(b)] For any $n \in \ndop \setminus \{0\}$, a metric $\metr$ on $L$ is a $\ghp(L)$-metric if and only if $\metr^{\otimes n}$ is a $\ghp(L)^{\otimes n}$-metric.
\item[(c)] If $\varphi:Y \rightarrow X$ is a morphism of proper schemes over $K$ and $\metr$ is a $\ghp(L)$-metric, then $\varphi^*\metr$ is a ${\hat{\frak g}_{Y}^+}(\varphi^*L)$-metric.
\item[(d)] If $\varphi$ is also surjective and $\metr$ is any metric on $L$ such that $\varphi^*\metr$ is a ${\hat{\frak g}_{Y}^+}(\varphi^*L)$-metric, then $\metr$ is a $\ghp(L)$-metric.
\item[(e)] The set $\ghp(L)$ is closed with respect to uniform convergence.
\item[(f)] If $\Lcal$ is a formal $\kcirc$-model of $L$ with numerically effective reduction $\tilde{\Lcal}$, then the associated formal metric $\metr_\Lcal$ is in $\ghp(L)$.
\end{itemize}
\end{prop}

\proof See \cite{Gu3}, Remark 10.3 and Proposition 10.4. \qed

\begin{art} \rm \label{admissible metrics}
Taking the intersection over all possible $\ghp(L)$ in Proposition \ref{ghp-candidates}, we get a smallest set of continuous metrics on $L^{\rm an}$ satisfying the properties of Proposition \ref{ghp-candidates}. Such a metric is called a {\it semipositive admissible metric}. A (continuous) metric $\metr$ on $L^{\rm an}$ is called an {\it admissible metric} if and only if there is a surjective morphism $\varphi:X'\rightarrow X$ of proper schemes over ${K}$, line bundles $M$, $N$ on $X'$ with $\varphi^*(L) \cong M \otimes N^{-1}$ and semipositive admissible metrics $\metr_1$, $\metr_2$ on $M$, $N$ such that $\varphi^* \metr = \metr_1/\metr_2$. 
\end{art}

\begin{prop} \label{properties of admissible metrics}
Admissible metrics of line bundles on a proper scheme $X$ over ${K}$ have the following properties:
\begin{itemize}
\item[(a)] The tensor product of admissible metrics is again admissible.
\item[(b)] The dual metric of an admissible metric is admissible.
\item[(c)] The pull-back of an admissible metric with respect to a morphism $\varphi:Y\rightarrow X$ of proper schemes over ${K}$ is an admissible metric.
\item[(d)] Every formal metric is admissible.
\end{itemize}
\end{prop}

\proof The base change of a proper surjective morphism is again proper and surjective which proves easily (a) and (c). Property (b) is trivial and (d) follows from \cite{Gu3}, Proposition 10.4. \qed

\begin{ex} \label{canonical metrics on abelian var} \rm 
Let $L$ be a line bundle on an abelian variety $A$ over ${K}$. We will define a canonical metric  on $L$ and then we will show that it is admissible.

We choose a rigidification $\rho$ of $L$, i.e. $\rho \in L_0({K}) \setminus \{0\}$. We assume first that $L$ is even. Then the theorem of the cube yields a canonical identification $[m]^*(L)=L^{\otimes m^2}$ of rigidified line bundles for every $m \in \zdop$. There is a unique bounded metric $\metr_\rho$ on $L$ such that for all $m\in \zdop$, we have $[m]^* \metr_\rho = \metr_\rho^{\otimes m^2}$. In fact, a variation of Tate's limit argument shows that
\begin{equation} \label{Tate's limit} 
\metr_\rho = \lim_{m \to \infty} [m]^*\metr^{\otimes \frac{1}{m^2}} 
\end{equation}
for every continuous metric $\metr$ on $L^{\rm an}$ (see \cite{BG}, Theorem 10.5.7). If $L$ is odd, then the same applies with $m^2$ replaced by $m$. Since any line bundle on $A$ is the tensor product of an even one with an odd one, unique up to $2$-torsion, we get a canonical metric $\metr_\rho$ on every rigidified line bundle $(L,\rho)$ of $A$.

A metric $\metr$ on $L$ is said to be {\it canonical} if there is a rigidification $\rho$ of $L$ such that $\metr$  is equal to $\metr_\rho$. A canonical metric is unique up to positive rational multiples (\cite{BG}, Remark 9.5.9) and we usually denote it by $\canmetr$. We claim that $\canmetr$ is admissible.

To see this, we assume first that $L$ is ample and even. By  Proposition \ref{ghp-candidates}(f), there is a semipositive admissible metric $\metr$ on $L$. Then Proposition \ref{ghp-candidates}(b) and \eqref{Tate's limit} yield that $\canmetr$ is a semipositive admissible metric.
If $L$ is just even, then there are ample even line bundles $M$, $N$ on $A$ with $L \cong M \otimes N^{-1}$ and we deduce that $\canmetr$ is admissible from the special case above and from Proposition \ref{properties of admissible metrics}. 


If $L$ is odd, then $L$ is algebraically equivalent to $0$. By definition, the latter means that we have $K$-rational points $P$ and $P_0$ on  an irreducible smooth projective curve $C$ over $K$ and a correspondence $E$ in $C \times A$ such that the line bundle associated to the divisor  $E([P]-[P_0])$ is isomorphic to $L$. Here, we use $E([P]-[P_0]):=(p_2)_*(E.p_1^*([P]-[P_0]))$, where $p_i$ are the projections of $C \times A$. There are semistable $\kcirc$-models $\Ccal$ and $\Acal$ of $C$ and $A$ (for curves, this is well-known and for abelian varieties, see Examples \ref{canonical measure on abelian variety} and \ref{spectrum example}). They are defined over the valuation ring $F^\circ$ of a finite extension $F$ over the completion $K_v$. More precisely, there are semistable algebraic $F^\circ$-models $\Cmcal$ and $\Amcal$ of $C_F$ and $A_F$  such that the associated formal schemes over $\kcirc$ are $\Ccal$ and $\Acal$, respectively (see \ref{schemes}). 

There is a divisor $\Dmcal$ on $\Cmcal$ with horizontal part $[P]-[P_0]$ and whose vertical part has rational coefficients such that the intersection numbers $\Dmcal\cdot Y$ are $0$ for all irreducible components $Y$ of $\tilde{\Cmcal}$. There is an $F^\circ$-model $\Zmcal$ of $C \times A$ with  $\kcirc$-morphisms $p_1:\Zmcal\rightarrow \Cmcal$ and  $p_2:\Zmcal \rightarrow \Amcal$ (extending the corresponding projections) such that the correspondence $E$ extends to a correspondence $\Emcal$ of  $\Zmcal$. We define $\Emcal(\Dmcal):=(p_2)_*(\Emcal.p_1^*(\Dmcal))$ as a $\qdop$-divisor on $\Amcal$. It is well-known that $\Emcal(Z)$  induces the canonical metric $\metr_{\rm can}$ of $L$. More precisely, if $N$ is a common denominator for the coefficients of the vertical part of $\Dmcal$, then the line bundle associated to the divisor $N\Dmcal$ induces a formal $\kcirc$-model of $L^{\otimes N}$ and $\metr_{\rm can}^{\otimes N}$ is the associated formal metric. 
Moreover, we deduce that $\metr_{\rm can}$ is a semipositive admissible metric. For more details about these constructions, we refer to  \cite{Gu3}, Theorem 8.9 and Example 10.11.

If $L$ is any line bundle on $A$, then we deduce that $\metr_{\rm can}$ is admissible by linearity and by the special cases above.
\end{ex}

In non-archimedean analysis, no good definition of Chern forms of metrized line bundles is known. However, Chambert-Loir has introduced a measure which is the analogue of top-dimensional wedge products of such Chern forms.

\begin{prop} \label{Chambert-Loir's measures}
There is an unique way to associate to any $d$-dimensional geometrically integral proper variety $X$ over $K$ and to any family of  admissibly metrized line bundles $\overline{L_1}, \dots, \overline{L_d}$ on $X$  a  regular Borel measure $c_1(\overline{L_1}) \wedge \cdots \wedge c_1(\overline{L_d})$ on $X^{\rm an}$ such that the following properties hold:
\begin{itemize}
\item[(a)] $c_1(\overline{L_1}) \wedge \cdots \wedge c_1(\overline{L_d})$ is multilinear and symmetric in $\overline{L_1}, \dots, \overline{L_d}$.
\item[(b)] If $\varphi:Y \rightarrow X$ is a morphism of $d$-dimensional geometrically integral proper varieties over ${K}$, then 
$$\varphi_* \left( c_1(\varphi^* \overline{L_1}) \wedge \dots \wedge c_1(\varphi^* \overline{L_d}) \right) = \deg(\varphi) c_1(\overline{L_1}) \wedge \dots \wedge c_1(\overline{L_d}).$$
\item[(c)] If the metrics of $\overline{L_1},\dots, \overline{L_d}$ are semipositive and $g \in C(\Xan)$, then
$$\left| \int_{X^{\rm an}} g \,c_1(\overline{L_1}) \wedge \dots \wedge c_1(\overline{L_d}) \right| \leq |g|_{\rm sup} \deg_{L_1, \dots , L_d}(X).$$
\item[(d)] If $\Xcal$ is a formal $\kcirc$-model of $X$ with reduced special fibre and if the metric of $\overline{L_j}$ is induced by a formal $\kcirc$-model $\Lcal_j$ of $L$ on $\Xcal$ for every $j=1, \dots, d$, then
$$c_1(\overline{L_1}) \wedge \dots \wedge c_1(\overline{L_d})= 
\sum_Y\deg_{\tilde{\Lcal}_1, \dots, \tilde{\Lcal}_d}(Y) \delta_{\xi_Y},$$
where $Y$ ranges over the irreducible components of $\tilde{\Xcal}$ and $\delta_{\xi_Y}$ is the Dirac measure in the unique point $\xi_Y$ of $\Xan$ which reduces to the generic point of $Y$.
\item[(e)] If $L_1, \dots, L_d$ are endowed with semipositive admissible metrics $\metr_j$, then $\mu=c_1(\overline{L_1}) \wedge \cdots \wedge c_1(\overline{L_d})$ is a positive Borel measure  with $\mu(\Xan)=\deg_{L_1,\dots,L_d}(X)$. 
\item[(f)] If we endow the set of positive regular Borel measures on $\Xan$ with the weak topology and if we fix the line bundles $L_1,\dots,L_d$ on $X$, then $c_1(\overline{L_1}) \wedge \dots \wedge c_1(\overline{L_d})$ is continuous with respect to the vector $(\metr_1, \dots, \metr_d)$ of semipositive admissible metrics on $L_1,\dots ,L_d$.
\end{itemize}
\end{prop}

\proof For existence, we refer to \cite{Gu4}, \S 3. Uniqueness for formal metrics is clear  by (d). The general case will be skipped. It follows from a repeated application of the minimality of semipositive admissible metrics in  \ref{ghp-candidates}. \qed

\begin{art} \rm \label{canonical measures}
We consider a  $d$-dimensional geometrically integral closed subvariety $X$ of the abelian variety $A$ and canonically metrized line bundles  $\overline{L_1}, \dots, \overline{L_d}$  on $A$. Then $\mu:=c_1(\overline{L}|_X) \wedge \cdots c_1(\overline{L}|_X)$ is called a {\it canonical measure} on $X$. It does not depend on the choice of the canonical metrics. Moreover, if one line bundle is odd, then $\mu=0$ (see \cite{Gu4}, 3.15).
\end{art}

\begin{rem} \rm 
By finite base change of $K$ and then using linearity in the irreducible components, we may extend Chambert-Loir's measures to all proper schemes $X$ over $K$ of pure dimension $d$.
\end{rem}

\section{Raynaud extensions and Mumford models}

In this section, $\kdop$ denotes an algebraically closed field with a non-trivial non-archimedean complete absolute value $|\phantom{a}|$,  valuation $v:=-\log |\phantom{a}|$ and value group  $\Gamma$. We fix an abelian variety $A$ over $\kdop$. 

First, we recall some results of Bosch and L\"utkebohmert about the Raynaud extension of $A$. To simplify the exposition, we will replace cubical line bundles  by the use of metrics.  Then we explain Mumford's construction, which gives an admissible formal $\kdop^\circ$-model $\Acal$ associated to certain polytopal decompositions of $\rdop^n$. Moreover, we will relate ample line bundles on $A$ and their models on $\Acal$  to affine convex functions. At the end, we will define the tropical variety of a closed subvariety of $A$ which is a periodic polytopal set in $\rdop^n$.

\begin{art} \rm  \label{Raynaud extension}
To define the Raynaud extension of $A$, we will follow the rigid analytic presentation of Bosch and L\"utkebohmert (\cite{BL2}, \S 1) 
and adapt it to Berkovich analytic spaces as in \cite{Ber}, \S 6.5. 
There is a formal group scheme $\Acal_1$ over $\kdop$ with generic fibre $A_1$ and a homomorphism $A_1 \rightarrow A^{\rm an}$ of analytic groups over $\kdop$ inducing an isomorphism onto an analytic subdomain of $A^{\rm an}$ such that $\Acal_1$ has semiabelian reduction. Moreover, $\Acal_1$ and the homomorphism $A_1 \rightarrow A^{\rm an}$ are unique up to isomorphism and hence we may identify $A_1$ with a compact subgroup of $A^{\rm an}$. 

It is convenient here to work in the category of formal analytic varieties as we may identify the objects with its generic fibres using a coarser topology (see \ref{formal analytic varieties}). Since $\Acal_1$ has semistable reduction, the special fibre is reduced and  $A_1$ has the structure of a formal analytic group. 
Let $T_1$ be the maximal formal affinoid subtorus in $A_1$. Then semistable reduction means that there is a unique formal abelian scheme $\Bcal$ over $\kdop^\circ$ with generic fibre $B$ such that we have an exact sequence
\begin{equation} \label{Raynaud extension 1}
1 \rightarrow T_1 \rightarrow A_1 \stackrel{q_1}{\rightarrow} B \rightarrow 0
\end{equation}
of formal analytic groups. Note that we may identify $T_1$ with the compact analytic subgroup $\{|x_1|= \dots = |x_n|=1\}$ of $T= (\Tor)^{\rm an}$. 
The {\it uniformization} $E$ of $A$ is given as an analytic group by $E:=(A_1 \times T)/T_1$, where $T_1$ acts on $A_1 \times T$ by $t_1 \cdot (a,t):=(t_1+a,t_1^{-1} \cdot t)$. Using the canonical maps, we get an exact sequence
\begin{equation} \label{Raynaud extension 2}
1 \rightarrow T \rightarrow E \stackrel{q}{\rightarrow} B\rightarrow 0
\end{equation}
of analytic groups. The closed immersion $T_1 \rightarrow A_1$ extends uniquely to a homomorphism $T \rightarrow A^{\rm an}$ of analytic groups and hence we get a unique extension  of $A_1 \hookrightarrow A^{\rm an}$ to a homomorphism $p:E \rightarrow A^{\rm an}$ of analytic groups. The kernel $M$ of $p$ is a discrete subgroup of $E(\kdop)$ and the homomorphism $E/M \rightarrow A^{\rm an}$ induced by $p$ is an isomorphism. The exact sequences \eqref{Raynaud extension 1} and \eqref{Raynaud extension 2} are called the {\it Raynaud extensions} of $A$. We will write the group structure on the uniformization $E$ multiplicatively. We call $n$ the {\it torus rank} of $A$.

By \cite{BGR}, Theorem 6.13, the formal abelian scheme $\Bcal$ is algebraizable and the GAGA principle shows that the same is true for the Raynaud extension \eqref{Raynaud extension 2}. 

There are two extreme cases of abelian varieties over $K$. First, we have abelian varieties of {\it  good reduction} at $v$ which means that the torus part $T$ of the Raynaud extension is trivial. On the other hand, we have the abelian varieties with {\it totally degenerate reduction} at v which means that the abelian part $B$ of the Raynaud extension is trivial.
\end{art}

\begin{art} \rm \label{Raynaud trivializations}
The Raynaud extension \eqref{Raynaud extension 1} is locally trivial, i.e. there is an open atlas $\Tfrak$ of $\Bcal^{\rm f-an}$ by formal affinoid varieties $V$ such that $q_1^{-1}(V) \cong V \times T_1$. This follows easily from the corresponding fact for semiabelian varieties applied to the reduction of \eqref{Raynaud extension 1} (see \cite{BL2}, p. 655). For every $V$, we fix such a trivialization given by a section $s_V:V \rightarrow A_1$. The transition functions $g_{VW}:=s_V -s_W$ are maps from $V \cap W$ to $T_1$. As usual, we fix coordinates $x_1, \dots , x_n$ on  $T=(\Tor)^{\rm an}$. The functions $x_1, \dots , x_n$ are defined on the trivialization $V \times T$ of $E$ by pull-back, but they do not extend to $E$. However, the functions $|x_1|, \dots, |x_n|$ are well defined on $E$ independently of the choice of the formal affinoid atlas $\Tfrak$. Using  $p(x_j)=|x_j|(p)$, we get a well-defined continuous surjective map
\begin{equation*}
\val: E \longrightarrow \rdop, \quad p \mapsto (-\log p(x_1), \dots , -\log p(x_n)). 
\end{equation*}
We will see at the end of this section that this map has similar properties as in tropical algebraic geometry, where one considers the special case $T=E$ and where no abelian variety is behind the construction. (In tropical algebraic geometry, this map is called the tropicalization map and it is also denoted by $\val$ to emphasize that it is obtained on rational points by applying the valuation to the coordinates.) Note that $\val$ maps the discrete subgroup $M$ of $A^{\rm an}$ isomorphically onto a complete lattice $\Lambda$ in $\rdop^n$ (\cite{BL2}, Theorem 1.2) and hence $\val$ induces a continuous surjective map 
$$\valbar: A^{\rm an} \rightarrow \rdop^n/\Lambda.$$ 
We will construct in Example \ref{canonical measure on abelian variety} a natural homeomorphism $\iota$ of $\rdop^n/\Lambda$ onto a compact subset $S(A)$ (called the skeleton) of $A^{\rm an}$. By \S 6.5 of \cite{Ber},  $\valbar \circ \iota$ gives a proper deformation retraction of $A$ onto $S(A)$.

If $\chi$ is an element of the character group $\check{T}$ of $T$, then the units $\chi^{-1}(g_{VW})$ are transition functions of a formal line bundle $\Ocal_\chi$ on $\Bcal$. Obviously, $s_V$ induces a trivialization $s_{V+b}(x)=s_V(x-b)+a$ of $A_1$ for all $a \in A_1$ with $q_1(a)=b$ and $x \in V+b$. Hence $O_\chi$ is a translation invariant line bundle proving that $O_\chi \in \Pic^\circ(B)$ and the same argument shows that the special fibre $\tilde{\Ocal}_\chi \in \Pic^\circ(\tilde{\Bcal})$. The translation invariance of $\Ocal_\chi$ can be also seen from the fact that  $\Ocal_\chi$ is given by the formal group schemes extension of $\Bcal$ by the formal multiplicative group obtained from the push-forward of the Raynaud extension by the character $\chi$.
We have the description
$$E= \Spec \left( \bigoplus_{\chi \in \check{T}} O_\chi \right)$$
of the Raynaud extension which is easily obtained by using the Laurent series development on the trivializations $V \times T$. Note that $q^* O_\chi$ is trivial on $E$ with canonical nowhere vanishing section $e_\chi$ given  by the functions $\chi$ on the trivialization $V \times T$ of $E$. Additional information for this and the next paragraph can be also found in the book of Fresnel and van der Put \cite{FvdP}, Chapter 6.
\end{art}

\begin{art} \label{line bundles} \rm
Next, we describe a line bundle $L$ on $A$ using the uniformization $E$. By $A^{\rm an}= E/M$, we see that $p^*(L^{\rm an})$ is equipped with an $M$-action $\alpha$ such  $L^{\rm an}$ may be recovered from $p^*(L^{\rm an})$ by passing to the quotient with respect to $\alpha$. There is a formal line bundle $\Hcal$ on $\Bcal$ with generic fibre denoted by $H$ such that $q^*(H)$ is isomorphic to $p^*(L^{\rm an})$ (see \cite{BL2}, Proposition 4.4). We fix such an isomorphism to get the identification $q^*(H)=p^*(L^{\rm an})$. Then  $q^*(\Hcal)$ is a formal $\kcirc$-model of $p^*(L^{\rm an})$ and as in Example \ref{formal models} , we get  a formal metric $q^*\metr_\Hcal$ on $p^*(L^{\rm an})$. There is a cocycle $Z$ of $H^1(M,(\rdop^\times)^E)$ such that
$$\left(q^*\|\alpha_\gamma(w)\|_\Hcal\right)_{\gamma \cdot x} = Z_\gamma(x)^{-1} \cdot \left(q^*\|w \|_\Hcal\right)_x$$
for all $\gamma \in M$, $x \in E$ and $w \in (p^*L^{\rm an})_x$. By  the description of the action $\alpha$ given in \cite{BL2}, Proposition 4.9, it is easy to deduce that $Z_\gamma(x)$ depends only on $\val(x)$. For $\lambda \in \Lambda$, we get a unique function $z_\lambda: \rdop^n \rightarrow \rdop$ with
$$z_\lambda(\val(x))=- \log(Z_\gamma(x)) \quad (\gamma \in M, x \in E, \lambda=\val(\gamma)).$$
Moreover, the same consideration shows that
$$z_\lambda(\ub) = z_\lambda(\mathbf 0) + b(\ub,\lambda) \quad (\ub \in \rdop^n, \lambda \in \Lambda)$$
for a symmetric bilinear form $b:\Lambda \times \Lambda \rightarrow \zdop$. By \cite{BL2}, Theorem 6.13, $L$ is ample if and only if $H$ is ample on $B$ and $b$ is positive definite on $\Lambda$. We note also that the bilinear form $b$ is trivial if $L \in \Pic^\circ(A)$ (use \cite{BL2}, Corollary 4.11).
\end{art}

\begin{art} \label{convex geometry} \rm 
We fix now the notation used from convex geometry (see  \cite{Gu4}, \S 6.1 and App. A, for more details). A polytope $\Delta$ of $\rdop^n$ is called {\it $\Gamma$-rational} if it may be given as an intersection of half spaces of the form $\{\ub \in \rdop^n \mid \mb \cdot \ub \geq c\}$ for suitable $\mb \in \zdop^n$ and $c \in \Gamma$. If $\Gamma = \qdop$, then such a polytope is called {\it rational}. The relative interior of $\Delta$ is denote by $\relint(\Delta)$. A {\it closed face} of $\Delta$ is either the polytope $\Delta$ itself or is equal to $H \cap \Delta$ where $H$ is the boundary of a half-space of $\rdop^n$ containing $\Delta$. An {\it open face} of $\Delta$ is the relative interior of a closed face.

A {\it polytopal decomposition} of $\Omega \subset \rdop^n$ is a locally finite family $\Ccal$ of polytopes contained in $\Omega$ which includes all faces, which is face to face and which covers $\Omega$. A {\it subdivision} $\Dcal$ of $\Ccal$ is a polytopal decomposition of $\Omega$ such that every $\Delta \in \Ccal$  has a polytopal decomposition in $\Dcal$.

We use the quotient map $\rdop^n \rightarrow \rtor, \, \ub \mapsto \ubb,$ to translate the above notions. A {\it polytope} $\Deltabar$ in $\rtor$ is given by a polytope $\Delta$ in $\rdop^n$ which maps bijectively onto $\Deltabar$. A {\it polytopal decomposition} $\Ccalbar$ of $\rtor$ is a a finite family of polytopes in $\rtor$ induced by a $\Lambda$-periodic polytopal decomposition $\Ccal$ of $\rdop^n$. 

We define convex functions as in analysis (and not as in the theory of toric varieties). A  convex function $f: \rdop^n \rightarrow \rdop$ is called  {\it strongly polyhedral} with respect to the polytopal decomposition $\Ccal$ of $\rdop^n$ if the $n$-dimensional polytopes in $\Ccal$ are the maximal subsets of $\rdop^n$ where $f$ is affine. 
\end{art}

\begin{art} \rm \label{polytopal domains}
A $\Gamma$-rational polytope $\Delta$ induces a {\it polytopal domain} $U_\Delta:=\val^{-1}(\Delta)$ of the torus $T$ with affinoid algebra
$$\kdop \langle U_\Delta \rangle := \left \{ \sum_{\mb \in \zdop^n} a_\mb x_1^{m_1} \dots x_n^{m_n} \mid \lim_{|\mb| \to \infty} v(a_\mb) + \mb \cdot \ub =\infty \, \; \forall \ub \in \Delta \right\}$$
(see \cite{Gu4}, Proposition 4.1). We need the following generalization:
\end{art}

\begin{lem} \label{relative polytopal domain}
Let $V$ be an affinoid variety with affinoid algebra $\Ocal(V)$. Then every $h \in \Ocal(V \times U_\Delta)$ has a Laurent series development
\begin{equation} \label{Laurent series development}
h= \sum_{\mb \in \zdop^n} a_\mb x_1^{m_1} \cdots x_n^{m_n} 
\end{equation}
for uniquely determined $a_\mb \in \Ocal(V)$ and the supremum semi-norm is given by
\begin{equation} \label{Laurent norm} 
|h|_{\sup} =\sup_{\ub \in \Delta,\mb \in \zdop^n} |a_\mb|_{\sup}e^{-\mb \cdot \ub}.
 \end{equation}
The supremum in \eqref{Laurent norm} is a maximum achieved in a vertex $\ub$ of $\Delta$. If $V$ is connected, then $h$ is a unit in $\Ocal(V \times U_\Delta)$ if and only if there is $\mb_0 \in \zdop^n$ such that $|a_{\mb_0}(y)\xb^{\mb_0}|>|a_\mb(y)\xb^\mb|$ for all $\xb \in U_\Delta$, $y \in V$ and $\mb \neq \mb_0$.

Conversely, a Laurent series as in \eqref{Laurent series development} is in $\Ocal(V \times U_\Delta)$ if and only if $-\log \|a_\mb\|+ \mb \cdot \ub$ tends to $\infty$ for $|\mb|\to \infty$, where $\metr$ is any Banach norm on the affinoid algebra $\Ocal(V)$. 
\end{lem}

\proof The description of $\Ocal(V\times U_\Delta)$ as the set of Laurent series \eqref{Laurent series development} is a direct generalization of \cite{Gu4}, Proposition 4.1. The proof follows the same arguments and will be omitted. It remains to prove the characterization of the units.

If $h \in \Ocal(V)$ has such a dominant term $a_{\mb_0}(y)\xb^{\mb_0}$, then $a_{\mb_0}$ has no zeros on $V$ and hence Hilbert's Nullstellensatz for affinoid algebras (\cite{BGR}, Proposition 7.1.3/1) shows that $a_{\mb_0} \in \Ocal(V)^\times$. We may assume $\mb_0=\mathbf 0$ and $a_{\mathbf 0}=1$. Then we have $h=1-h_1$ for $h_1 \in \Ocal(V \times U_\Delta)$ with $|h_1|_{\rm sup}<1$ and hence 
$$h^{-1}= \sum_{m=0}^\infty h_1^m \in \Ocal(V\times U_\Delta).$$
If $h \in \Ocal(V \times U_\Delta)$ has no such dominant term, then there is $\xb \in U_\Delta$, $y \in V$ and $\mb_0 \neq \mb_1 \in \zdop^n$ with
$$|a_{\mb_0}(y) \xb^{\mb_0}| = |a_{\mb_1}(y) \xb^{\mb_1}|=|h|_{\rm sup}.$$
Let $\ub := \val(\xb)$ and let $W$ be the affinoid subdomain of $V \times U_\Delta$ given by $\val^{-1}(\ub)$. It is isomorphic to $V \times T_1$ for the affinoid torus $T_1=\{|x_1|= \dots = |x_n|=1\}$ in $T$. Since the restriction of $h$ to $W$ has no dominant term as well, we get $h|_W \not \in \Ocal(W)^\times$ (\cite{BGR}, Lemma 9.7.1/1) and hence $h \not \in \Ocal(V \times U_\Delta)^\times$. \qed

\begin{art} \rm \label{Mumford's construction}
Next, we 
define a formal $\kdop^\circ$-model $\Acal$ of $A$ associated to a $\Gamma$-rational polytopal decomposition $\Ccalbar$ of $\rtor$. In the algebraic framework, this is a construction of Mumford \cite{Mu} which is useful for compactifying moduli spaces of abelian varieties (see \cite{FC}).
We denote by $\Ccal$ the $\Gamma$-rational $\Lambda$-periodic polytopal decomposition of $\rdop^n$ which induces $\Ccalbar$. 

We choose a formal affinoid atlas $\Tfrak$ as in \ref{Raynaud trivializations}. For $V \in \Tfrak$ with trivialization $q_1^{-1}(V) \cong V \times T_1$ and $\Delta \in \Ccal$, we define the affinoid subdomain
\begin{equation} \label{relative polytopal}
U_{V,\Delta}:= q^{-1}(V) \cap \val^{-1}(\Delta) \cong V \times U_\Delta
\end{equation}
of $E$, where the term on the right is in the trivialization $q^{-1}(V) \cong V \times T$. The sets  $U_{V,\Delta}$ form a formal analytic atlas on $E$ inducing a formal analytic variety $\Ecal^{\rm f-an}$ with corresponding formal $\kdop^\circ$-model $\Ecal$ of $E$. We note that $\Ecal$ has a formal affine open covering by the sets $\Ucal_{V,\Delta}:=U_{V,\Delta}^{\rm f-sch}$.

We may assume that $\Tfrak$ is closed under translation with elements of $q(M)$. We may form the quotient of $\Ecal^{\rm f-an}$  by $M$ leading to a formal analytic structure on $A^{\rm an}$. The associated formal $\kdop^\circ$-model $\Acal$ of $A$ (see \ref{admissible formal schemes}) is called the {\it Mumford model} associated to $\Ccalbar$. It has a covering by formal affine open subsets $\Ucal_{[V,\Delta]}$ obtained by gluing $\Ucal_{V+q(\gamma),\Delta+\val(\gamma)}$ for all $\gamma \in M$. Obviously, $\Acal$ is independent of the choice of $\Tfrak$. Note that we have canonical morphisms $q: \Ecal \rightarrow \Bcal$ and $p:\Ecal \rightarrow \Acal$ extending the corresponding maps on generic fibres.
\end{art}

Recall that the strata of a variety were introduced at the end of \S 1. The next result describes the strata of the special fibre of a Mumford model.

\begin{prop} \label{strata on Mumford models}
Let $\Acal$ be the Mumford model of $A$ associated to the $\Gamma$-rational polytopal decomposition $\Ccalbar$ of $\rtor$. Let  $\Ecal$ be the formal $\kdop^\circ$-model  of $E$ associated to the polytopal decomposition $\Ccal$ of $\rdop^n$ which was used in \ref{Mumford's construction} to construct $\Acal$.  
\begin{itemize}
\item[(a)] The formal torus $\tdop_1= \Spf(\kdop^{\circ} \langle x_1^{\pm1}, \dots,x_n^{\pm 1} \rangle)$ acts canonically on $\Ecal$ inducing a $(\Tor)_{\tilde{\kdop}}$-action on the special fibre $\tilde{\Ecal}$. 
\item[(b)] There is a bijective order reversing correspondence between strata $Z$ of $\tilde{\Ecal}$ and open faces $\tau$ of $\Ccal$. It is given by 
$$\tau= \val(\pi^{-1}(Z)), \quad Z = \pi(\val^{-1}(\tau)),$$
where $\pi: E \rightarrow \tilde{\Ecal}$ is the reduction map. We have $\dim(Z) + \dim(\tau) = \dim(A)$.
\item[(c)] There is a bijective order reversing correspondence between strata $W$ of $\tilde{\Acal}$ and open faces $\overline{\tau}$ of $\Ccalbar$. It is given by 
$$\overline{\tau}= \valbar(\pi^{-1}(W)), \quad W = \pi(\valbar^{-1}(\overline{\tau})),$$
where $\pi: A \rightarrow \tilde{\Acal}$ is the reduction map. We have $\dim(W) + \dim(\overline{\tau})= \dim(A)$.
\item[(d)] Every irreducible component $Y'$ of $\tilde{\Ecal}$ is mapped isomorphically onto an irreducible component $Y$ of $\tilde{\Acal}$. By (c), we get a bijective correspondence between irreducible components of $\tilde{\Acal}$ and vertices of $\Ccalbar$. Moreover,  $\tilde{q}:Y'\rightarrow \tilde{\Bcal}$ is a fibre bundle whose fibre is a $(\Tor)_{\tilde{\kdop}}$-toric variety.
\end{itemize}
\end{prop}

\proof By construction, $\tdop_1^{\rm f-an}$ acts on $\Ecal^{\rm f-an}$ and (a) follows. 
To prove (b), we note that strata are compatible with localization and hence we may consider a formal affinoid chart $U_{V,\Delta} \cong V \times U_\Delta$ as in \eqref{relative polytopal}. 
By \cite{Gu4}, Proposition 4.4, it follows that 
the strata of $\tilde{U}_\Delta$ are the same as the $(\Tor)_{\ktilde}$-orbits and they correspond to the open faces of $\Delta$. 
The strata of $\tilde{U}_{V,\Delta}$ are the preimages of the strata of $\tilde{U}_\Delta$ leading to the desired correspondence. The other claims in (b) follow also from the corresponding statements for $U_\Delta$ given in \cite{Gu4}, Proposition 4.4.


By (b) and the construction of $\Acal$, $\tilde{p}$ maps a stratum of $\tilde{\Ecal}$  isomorphically onto a stratum of $\tilde{\Acal}$ and hence (c) follows easily from (b).
To prove (d), let $\ub$ be the vertex of $\Ccal$ corresponding to the irreducible component $Y'$ by (b). Let $\Delta \in \Ccal$ with vertex $\ub$. In the trivialization \eqref{relative polytopal}, $Y'$ is given by $\tilde{V} \times Y_{\Delta,\ub}$, where $Y_{\Delta,\ub}$ is the affine toric variety given by the  local cone of $\Delta$ in $\ub$ (see \cite{Gu2}, Proposition 4.4(d)). If $\Delta$ ranges over $\Ccal$, we see that $Y'$ has over $V$ the form $V \times Y_\ub$, where $Y_\ub$ is the $(\Tor)_{\tilde{\kdop}}$-toric variety given by the fan of local cones of the polytopes $\Delta \in \Ccal$ in the vertex $\ub$. 
This can be done for every $V \in \Tfrak$ to cover $Y'$. We note that $Y'$ is the union of the strata corresponding to the open faces $\tau$ of $\Ccal$ with vertex $\ub$. Since $Y'$ is locally isomorphic to $Y$ and no gluing arises with respect to the $M$-action, we easily deduce (d). \qed

\begin{rem} \rm \label{strata as torsors}
Let $\Delta \in \Ccal$ with relative interior $\tau$. We denote by $\ldop_\Delta$ the linear subspace of $\rdop^n$ generated by $\Delta - \ub$, $\ub \in \Delta$. Then $N_\Delta:= \ldop_\Delta \cap \zdop^n$ is a complete lattice in $\ldop_\Delta$ and we get a subtorus $H_\Delta$ of $\tilde{T}=(\Tor)_{\tilde{\kdop}}$ with $H_\Delta(\tilde{\kdop})=N_\Delta \otimes_\zdop \tilde{\kdop}^\times$. It follows from the above proof and \cite{Gu4}, Remark 4.8, that the stratum of $\tilde{\Ecal}$ associated to $\tau$ is a $\tilde{T}/H_\Delta$-torsor over $\tilde{\Bcal}$ with respect to $\tilde{q}$.
\end{rem}

\begin{art} \rm \label{Mumford line bundles}
Next, we describe  $\kdop^\circ$-models of the line bundle $L$ on $A$. They should be defined on a given Mumford model $\Acal$ of $A$ associated to the $\Gamma$-rational polytopal decomposition $\Ccalbar$ of $\rtor$. As in \ref{line bundles}, we choose a formal line bundle $\Hcal$ on $\Bcal$ with $q^*(H)=p^*(L^{\rm an})$ on the uniformization $E$ of $A$ such that $L^{\rm an}= q^*(H)/M$ on $A^{\rm an}=E/M$ and which leads to a cocycle $z_\lambda(\ub)$ with respect to $\lambda \in \Lambda$. We fix a formal affine atlas  of $\Bcal$ which trivializes the line bundle $\Hcal$ and which induces a formal affinoid trivialization $\Tfrak$ for $q_1:A_1 \rightarrow B$.  
\end{art}

\begin{prop} \label{Mumford line models}
There is a bijective correspondence between isomorphism classes of formal $\kdop^\circ$-models $\Lcal$ of $L$ on $\Acal$ which are trivial over the formal open subsets $\Ucal_{[V,\Delta]}$, where $\Delta \in \Ccal, V \in \Tfrak$, and continuous real functions $f$ on $\rdop^n$ satisfying the following two conditions: 
\begin{itemize}
\item[(a)] For $\Delta \in \Ccal$, there are $\mb_\Delta \in \zdop^n$ and $c_\Delta \in \Gamma$ with $f(\ub)= \mb_\Delta \cdot \ub + c_\Delta$ on $\Delta$.
\item[(b)] $f(\ub + \lambda) = f(\ub) + z_\lambda(\ub) \quad (\lambda \in \Lambda, \ub \in \rdop^n)$.
\end{itemize}
Let $\metr_\Lcal$ be the formal metric of $L$ associated to $\Lcal$ (see Example \ref{formal models}). Then the corresponding $f_\Lcal:\rdop^n \rightarrow \rdop$ is uniquely determined by 
\begin{equation} \label{f-char}
f_\Lcal \circ \val = - \log \circ \left( p^* \metr_\Lcal / q^* \metr_\Hcal \right),
\end{equation}
where the quotient of the metrics on $q^*(H)=p^*(L^{\rm an})$ is evaluated at any non-zero local section. 
\end{prop}

\proof Let $\Lcal$ be a formal $\kdop^\circ$-model of $L$ on $A$ which is trivial on every $\Ucal_{[V,\Delta]}$. Using the identification $q^*(H)=p^*(L^{\rm an})$, we may view $p^* \metr_\Lcal / q^* \metr_\Hcal$ as a metric on $O_E$. The corresponding real function is obtained by evaluating this metric at the constant section $1$. Since formal metrics are continuous, the right hand side of \eqref{f-char} is a continuous function on $E$. 

Our first goal is to show that this function descends to $\rdop^n$, i.e. there is $f_\Lcal:\rdop^n \rightarrow \rdop$  with \eqref{f-char}. We choose a connected $V \in \Tfrak$ and $\Delta \in \Ccal$. By assumption, the formal affine open subset $\Vcal$ with generic fibre $V$ trivializes the formal line bundle $\Hcal$, i.e. there is a nowhere vanishing section $s_V \in \Gamma(\Vcal,\Hcal)$. We may consider $s_V$ as a section of $H|_V$ with $\|s_V\|_\Hcal=1$ on $V$. By assumption, $\Lcal$ is trivial on $\Ucal_{[V,\Delta]}$ and hence there is a nowhere vanishing section $t_V \in \Gamma(V,L)$ with $\|t_V\|_\Lcal=1$ on $U_{[V,\Delta]}$. We apply Lemma \ref{relative polytopal domain} to the unit $h:=q^*(s_V)/p^*(t_V)$ on $U_{V,\Delta} \cong V \times U_\Delta$, hence there is $\mb_\Delta \in \zdop^n$ and $a_{V,\Delta} \in \Ocal(V)^\times$ with
\begin{equation} \label{quotient and h}
p^* \metr_\Lcal / q^* \metr_\Hcal = |h| = |a_{V,\Delta} \xb^{\mb_\Delta}|
\end{equation}
on $U_{V,\Delta}$. A priori, $\mb_\Delta$ depends also on $V$, but since the functions $|x_i|$ are well-defined on $E$, it follows easily from \eqref{quotient and h} that we may select $\mb_\Delta$ independently from $V \in \Tfrak$. We conclude that $|a_{V,\Delta} |= |a_{W,\Delta} |$ on $V \cap W$ for every $V,W \in \Tfrak$. If we vary $V \in \Tfrak$ keeping $\Delta$ fixed, we get a formal $\kdop^\circ$-model $\Gcal$ of $O_B$ on $\Bcal$, given by trivializations $a_{V,\Delta} \in \Ocal(V)^\times$. Since the special fibre $\tilde{\Bcal}$ is smooth, the formal metric $\metr_\Gcal$ on $O_B$ induces a constant function $\|1\|_\Gcal$ (see \cite{Gu2}, Proposition 7.6). This means that $|a_{V,\Delta} |$ is constant on $V$ and hence there is $a_\Delta \in \kdop^\times$ with $|a_{V,\Delta}|=|a_\Delta|$ for all $V \in \Tfrak$. For $x \in U_{V,\Delta}$, we conclude that $|h(x)|$ in \eqref{quotient and h} depends only on $\val(x)$ and hence there is a unique function $f_\Lcal$ with \eqref{f-char}.  
 Moreover, we have proved that (a) holds with $c_\Delta:=v(a_\Delta)$. 
Since $\Ccal$ is a polytopal complex, it is clear that continuity of $f_\Lcal$ follows from (a).

Finally, we prove (b) for $f_\Lcal$. Let $x \in E$ with $\val(x)=\ub$ and let $\gamma \in M$ with $\val(\gamma)=\lambda$. Then (b) follows from
\begin{align*}
f_\Lcal \circ \val(\gamma \cdot x) &= -\log \left( \left(p^*\metr_\Lcal \right)_{\gamma\cdot x} / \left( q^* \metr_\Hcal \right)_{\gamma \cdot x} \right)\\
&\stackrel{\ref{line bundles}}{=} - \log \left( \left(p^*\metr_\Lcal \right)_{x} / \left( e^{z_\lambda(\ub)} q^*\metr_{\Hcal} \right)_x \right)\\
&= f_\Lcal(\ub) + z_\lambda(\ub).
\end{align*}
Conversely, let $f:\rdop^n \rightarrow \rdop$ be a continuous function satisfying (a) and (b). We define a metric $\metr'$ on $p^*(L^{\rm an})=q^*(H)$ by 
$$\metr'/q^*\metr_\Hcal = e^{-f \circ \val}.$$
As a consequence of (b), $\metr'$ descends to a metric $\metr_f$ on $L^{\rm an}=p^*(L^{\rm an})/M$. It is uniquely characterized by the property
\begin{equation} \label{metr_f}
f \circ \val = - \log \left( p^*\metr_f / q^* \metr_\Hcal \right).
\end{equation}
We choose $\mb_\Delta \in \zdop^n$ and $c_\Delta \in \Gamma$ from (a). There is $a_\Delta \in \kdop^\times$ with $c_\Delta = v(a_\Delta)$. For $V \in \Tfrak$ and $\Delta \in \Ccal$, the metric $p^*\metr_f$ is given on $U_{V,\Delta}$ by
\begin{equation} \label{quotient and a}
p^*\metr_f / q^*\metr_\Hcal = |a_\Delta| \cdot |\xb^{\mb_\Delta}|
\end{equation}
as a consequence of \eqref{metr_f}. Using $s_V \in \Gamma(V, H)$ from above, we deduce that the nowhere vanishing section $t_{V,\Delta}:= \left(a_\Delta \xb^{\mb_\Delta} \right)^{-1} \cdot q^*(s_V) \in \Gamma(U_{V,\Delta},q^*(H))$ satisfies
$$p^*\|t_{V,\Delta}\|_f =1$$
on $U_{V,\Delta}$. We may view  $(t_{V,\Delta})_{V \in \Tfrak, \Delta \in \Ccal}$ as
a family of frames of $L^{\rm an}$ of constant $\metr_f$-metric $1$ and hence $\metr_f$ is  the metric on $L$ associated to a unique formal $\kdop^\circ$-model $\Lcal_f$ of $L$ (\cite{Gu2}, Proposition 7.5) as desired. 

It remains to show that $f \mapsto \Lcal_f$ is inverse to $\Lcal \mapsto f_\Lcal$. Here, the same argument as for Proposition 6.6 in \cite{Gu4} applies. \qed

\begin{prop} \label{relatively ampel}
Let $\Ccalbar$ be a $\Gamma$-rational polytopal decomposition of $\rtor$ with associated Mumford model  $\Acal$ of $A$ over $\kdop^\circ$ obtained from the formal $\kdop^\circ$-model $\Ecal$ of the Raynaud extension $E$ as in \ref{Mumford's construction}. We assume that there is a $\kdop^\circ$-model $\Lcal$ of $L$ on $\Acal$ as in Proposition \ref{Mumford line models} corresponding to the affine function $f_\Lcal$. Let  $p:\Ecal \rightarrow \Acal$ be the quotient map. Then $(p^*\Lcal)\sptilde$ is relatively ample with respect to the canonical reduction $\tilde{q}:\tilde{\Ecal} \rightarrow \tilde{\Bcal}$ if and only if $f_\Lcal$ is a strongly polyhedral convex function with respect to $\Ccal$ (see \ref{convex geometry}). 
\end{prop}

\proof Let $\ub$ be a vertex of $\Ccal$. By Proposition \ref{strata on Mumford models}, we get a corresponding irreducible component $Y_\ub'$ of $\tilde{\Ecal}'$. Moreover, we have seen that $Y_\ub'$ is  a fibre bundle over $\tilde{\Bcal}$ which is trivial over $\tilde{\Vcal}=\tilde{V}$ for every $V \in \Tfrak$ with associated formal scheme $\Vcal$ over $\kcirc$. The fibre $Z_\ub$ of the bundle is the $(\Tor)_{\tilde{\kdop}}$-toric variety associated to the local cones of  the polytopes $\Delta \in \Ccal$ with vertex $\ub$. We claim that the restriction of $(p^*\Lcal)\sptilde$ to the trivialization $\tilde{V} \times Z_\ub$ is the pull-back of a line bundle on $Z_\ub$. Indeed, the toric variety $Z_\ub$ is given by pasting the family
$$(U_\Delta)\sptilde := \Spec\left(\tilde{\kdop}[\tilde{\xb}^{S_\Delta}]\right),$$
where $\Delta$ is ranging over the polytopes of $\Ccal$ with vertex $\ub$ and 
$$\tilde{\xb}^{S_\Delta}:= \{ \tilde{\xb}^\mb \mid \mb \in \zdop^n, \, \ub' \cdot \mb \geq 0 \; \forall \ub' \in \Delta - \ub\}.$$
For such a $\Delta$, we use the presentation  $f_\Lcal(\ub')=\mb_\Delta \cdot \ub' + c_\Delta$ from Proposition \ref{Mumford line models}(a). If we change the identification $p^*(L^{\rm an})=q^*(H)$, then $\metr_\Lcal$ is replaced by a positive multiple and hence we may assume $f_\Lcal(\ub)=0$.
There is $a_\Delta \in \kdop^\times$ with $c_\Delta= v(a_\Delta)$. The functions $\tilde{a}_\Delta\tilde{\xb}^{\mb_\Delta}|_{(U_\Delta)\sptilde}$ define a $(\Tor)_{\tilde{\kdop}}$-equivariant Cartier divisor $D$ on $Z_\ub$. Since $\Hcal$ is trivial over $\Vcal$ (see \ref{Mumford line bundles}), we deduce easily from Proposition \ref{Mumford line models} that $\tilde{p}^*\tilde{\Lcal}|_{\tilde{V} \times Z_\ub}$ is isomorphic to the pull-back of $O_{Z_\ub}(D)$ with respect to the second projection.

The claim follows from the fact that $D$ is ample if and only if $f_\Lcal$ is a strongly polyhedral convex function in the vertex $\ub$ (\cite{Fu2}, \S 3.4). \qed

\begin{art} \label{tropical variety} \rm 
Let $X$ be a closed subscheme of $A$. Then the subset $\valbar(X^{\rm an})$ of $\rtor$ is called the {\it tropical variety} associated to $X$. 

We note the analogue to tropical algebraic geometry, where one studies the tropical variety associated to an algebraic subvariety of $\Tor$. However, the lift $p^{-1}(X^{\rm an})$ of $X$ to the Raynaud extension $E$ is only an analytic subvariety and hence our tropical varieties are best studied in the framework of Berkovich analytic spaces (see \cite{Gu4} for details about tropical analytic geometry).
\end{art}

\begin{prop} \label{tropical polytopal set}
The tropical variety $\valbar(X^{\rm an})$ is a finite union of $\Gamma$-rational polytopes in $\rtor$ of dimension at most $\dim(X)$. If $X$ is connected, then the tropical variety is also connected.
\end{prop} 

\proof Let $E$ be the Raynaud extension of $A$ and let $\Tfrak$ be an atlas of trivalizations of $E$ over $B$ as in \ref{Raynaud trivializations}. 
We choose any $\Gamma$-rational polytope $\Delta$ of $\rdop^n$ inducing a polytope $\Deltabar$ of $\rtor$ and $V \in \Tfrak$. The  trivialization leads to $U_{V, \Delta} \cong V \times U_\Delta$. By \cite{Ber5}, Corollary 6.2.2, $\val(U_{V, \Delta} \cap p^{-1}(X^{\rm an}))$ is a finite union of $\Gamma$-rational polytopes in $\rdop^n$ of dimension at most $\dim(X)$. Since $A^{\rm an}$ is covered by finitely many $U_{[V,\Delta]}$, we conclude easily that $\valbar(X^{\rm an})$ is a finite union of $\Gamma$-rational polytopes in $\rtor$ of dimension at most $\dim(X)$. If $X$ is connected, then $X^{\rm an}$ is also connected (\cite{Ber}, Theorem 3.4.8). By continuity of $\valbar$, we see that $\valbar(X^{\rm an})$ is also connected. \qed

\begin{thm} \label{dimension of tropical variety}
Let $X$ be purely $d$-dimensional closed subscheme of $A$ and let $b$ be the dimension of the abelian variety $B$ of good reduction in the Raynaud extension \eqref{Raynaud extension 1} of $A$. Then  the tropical variety $\valbar(X^{\rm an})$ is a finite union of $\Gamma$-rational polytopes in $\rtor$ of dimension at least $d-b$ and at most $d$.
\end{thm}

\proof By Proposition \ref{tropical polytopal set}, there are $\Gamma$-rational polytopes $\overline{\Delta_1}, \dots \overline{\Delta_k} \in \rtor$ of dimension at most $d$ with $\valbar(X^{\rm an})=\overline{\Delta_1} \cup \dots \cup \overline{\Delta_k}$. We may assume that  no polytope $\overline{\Delta_j}$ may be omitted in this decomposition. We have to prove $\dim(\overline{\Delta_j}) \geq d-b$. For $\ubb \in \valbar(X^{\rm an})$, it is enough to show that the dimension of the polytopal set $\valbar(X^{\rm an})$ in a neighbourhood of $\ubb$ is at least $d-b$. 

We choose a lift $\ub$ of $\ubb$ to $\rdop^n$ and an $n$-dimensional $\Gamma$-rational polytope $\Delta$ with $\ub \in \relint(\Delta)$. Using the notation of the proof of Proposition \ref{tropical polytopal set}, we choose $V \in \Tfrak$ such that $\ub \in \val(X_{V,\Delta})$ for $X_{V,\Delta}:=U_{V, \Delta} \cap p^{-1}(\Xan)$. The claim follows now from the following more general result:

{\it Let $Y$ be any closed analytic subvariety of $U_{V,\Delta}$ of pure dimension $d$ such that $\ub \in \val(Y)$ and
let $N$ be the dimension of $\val(Y)$ in a neighbourhood of $\ub$. Then we have  $N \geq d-b$.}

The proof is by induction on $N$ and follows  \cite{Gu4}, Proposition 5.4.
If $N=0$, then we may assume that $\val(Y)=\{\ub\}$ by passing to a smaller $\Delta$. By our choice of $V$, we have the  trivialization $U_{V,\Delta}\cong V \times U_\Delta$, where $U_\Delta$ is the polytopal domain in $(\Tor)^{\rm an}$ associated to $\Delta$. Passing to the associated admissible formal affine $\kcirc$-schemes, we get $\Ucal_{V,\Delta}\cong \Vcal \times \Ucal_\Delta$.  By abuse of notation, we will use the  projection $p_2$  also on $\Ucal_{V,\Delta}$. Let $\Ycal$ be the closure of $Y$ in $\Ucal_{V,\Delta}$. 
By \cite{Gu4}, Proposition 4.4, the open face $\tau:=\relint(\Delta)$ induces the stratum $Z_\tau:=\pi(\val^{-1}(\tau))$ of $\tilde{U}_\Delta$ of dimension $n-\dim(\tau)=0$, where $\pi:U_\Delta \rightarrow \tilde{U}_\Delta$ is the reduction map. Now we use that $\pi \circ p_2^{\rm an} = \tilde{p}_2 \circ \pi$ on $Y$ and that the reduction  $\pi$ on the right hand side is a surjective map from $Y$ onto the special fibre $\tilde{\Ycal}$ (see \ref{admissible formal schemes}). We conclude from $\val(Y)=\{\ub\}$ that $\tilde{p}_2$ maps $\tilde{\Ycal}$  to the closed point $Z_\tau$. Since $Y$ is  of pure dimension $d$, the same is true for the special fibre $\tilde{\Ycal}$ and hence we get $$d \leq \dim(\tilde{p}_2^{-1}(Z_\tau)) \leq \dim(\tilde{\Vcal})=b.$$
This proves the claim for $N=0$.

Now we prove the case $N>0$. By \cite{Gu4}, Proposition 5.2, $\val(Y)$ is a finite union of $\Gamma$-rational polytopes. We conclude that $\ub$ is contained in an $N$-dimensional $\Gamma$-rational polytope $\sigma \subset \val(Y)$. Note that any point $\ub' \in \sigma$ contained in a sufficiently small neighbourhood of $\ub$ has also local dimension $N$. 
By density of $Y(\kdop)  \in Y^{\rm an}$ (\cite{Ber}, Proposition 2.1.15), we find such an $\ub'$ with  $\ub' =\val(y)$ for some $y \in Y(\kdop)$. Moreover, we may assume that $\ub'$ has an $n$-dimensional $\Gamma$-rational polytope $\Delta'$ as a neighbourhood such that $\Delta' \cap \val(Y) = \Delta' \cap \sigma$.
There are $\alpha \in \kdop$ and $\mb \in \zdop^n$ such that the hyperplane $H=\{\xb^\mb=\alpha\}$  passes through $y$ and such that  $\val(H^{\rm an})=\{\omega \cdot \mb =v(\alpha)\}$ intersects $\val(Y)\cap \Delta'$ transversally. By Krull's Hauptidealsatz, the closed analytic subvariety $Y':= Y \cap H^{\rm an} \cap U_{V,\Delta'}$ of $U_{V,\Delta'}$ has pure dimension $d-1$. We deduce from
$$\val(Y') \subset \val(Y) \cap \{\omega \cdot \mb = v(\alpha)\} \cap \Delta',$$
 that $\val(Y')$ has dimension $N' \leq N-1$ in a neighbourhood of $\ub'$. By induction applied to $Y'$, we get $N-1 \geq N' \geq d-1-b$ proving the claim. \qed

\begin{rem} \label{wrong argument} \rm
We assume now that $X$ is an irreducible $d$-dimensional closed subvariety of $A$. In the preprint  \cite{Gu7} of this paper, it was claimed in Theorem \ref{dimension of tropical variety} that $\valbar(\Xan)$ is of pure dimension. As pointed out by the referee, the argument was based on a wrong application of Chevalley's theorem which does not hold in the category of analytic spaces and so this question remains open. 

{\it However, if $A$ is isogeneous to $B_1 \times B_2$, where $B_1$ (resp. $B_2$) is an abelian variety with good  (resp. with totally degenerate) reduction at $v$, then $\valbar(\Xan)$ has indeed pure dimension $d-e$ for some $e \in \{0,\min(b,d)\}$.}

To prove this, let $\varphi:A \rightarrow B_1 \times B_2$ be an isogeny. By \cite{BL2}, Theorem 1.2, $\varphi$ lifts to an isogeny $\phi:E \rightarrow B_1^{\rm an} \times (\Tor)^{\rm an}$ between the associated uniformizations of the Raynaud extension. Obviously, $(\Tor)^{\rm an}$ is also the torus part in the Raynaud extension of $A$ and $\phi$ restricts to an isogeny $(\Tor)^{\rm an} \rightarrow (\Tor)^{\rm an}$. On the other hand, an (analytic) endomorphism of $\Tor$ is given by $\phi^*(x_j)=\xb^{\mb_j}$ for some $\mb_j \in \zdop^n$, $j=1,\dots,n$. We conclude that the linear isomorphism $\phi_{\rm aff}$, given by the matrix $(\mb_1, \dots, \mb_n)^t$, maps $\valbar(\Xan)$ onto $\valbar(\phi(X)^{\rm an})$. Hence we may assume that $A=B_1 \times B_2$. Let $p_2:A \rightarrow B_2$ be the second projection. Then $Y=p_2(X)$ is an irreducible closed subvariety of $B_2$ of dimension $d-e$ for some $e \in \{0,\min(b,d)\}$ with $b=\dim(B_1)$. By construction, we have $\valbar(\Xan)=\valbar(Y^{\rm an})$ and hence the claim follows from the fact that the tropical variety of an irreducible $d'$-dimensional closed subvariety of a totally degenerate abelian variety has pure dimension $d'$ (see Theorem \ref{dimension of tropical variety} or \cite{Gu4}, Theorem 6.9). \qed
\end{rem}

\section{Subdivisions of the skeleton}

In this section, $\kdop$ denotes an algebraically closed field endowed with a non-trivial non-archimedean complete absolute value $|\phantom{a}|$. Let $v:=-\log |\phantom{a}|$ be the valuation with value group $\Gamma := v(\kdop^\times)$, valuation ring $\kdop^\circ$ and  residue field $\tilde{\kdop}$. 

A smooth variety $X'$ over $\kdop$ has not always a smooth formal $\kcirc$-model and hence we study a strictly semistable $\kcirc$-model $\Xcal'$. Its special fibre $\tilde{\Xcal}'$ may be viewed as a divisor with normal crossings on $\Xcal'$. The skeleton of $\Xcal'$ is a metrized polytopal set of $(X')^{\rm an}$ closely related to the stratification of $\tilde{\Xcal}'$. 
We will see that the skeleton has similar properties as a tropical variety.

We will study the effect of subdivisions on the models. In particular, this is interesting if $X'$ maps to an abelian variety. The most important result of this somehow technical section is at the end, where we will compute the degree of an irreducible component of $\tilde{\Xcal'}$ in this setting under a certain transversality assumption. In the next section, this result is used to compute canonical measures on abelian varieties.

\begin{art} \rm  \label{strictly semistable}
Let $\Xcal'$ be a {\it strictly semistable} admissible formal scheme over $\kcirc$, i.e. $\Xcal'$ is covered by formal open subsets $\Ucal'$ with an \'etale morphism 
$$\psi: {\Ucal'} \longrightarrow {\Scal} :=\Spf \left( \kdop^\circ \langle x_0', \dots, x_d' \rangle / \langle x_0' \dots x_r' - \pi \rangle \right)$$
for  $r \leq n$ and $\pi \in \kdop^\times$ with $|\pi|<1$. The generic fibre $U'$ of $\Ucal'$ is smooth and hence the generic fibre $X'$ of $\Xcal'$ is a smooth analytic space. For simplicity, we assume that $X'$ is connected. Then $X'$ is $d$-dimensional, but $r$ and $\pi$ may depend on the choice of $\Ucal'$.

Note 
$\Scal = {\Spf \left( \kcirc \langle x_0', \dots , x_r' \rangle / \langle x_0' \cdots x_r' - \pi \rangle \right)} \times {\Spf \left( \kcirc \langle x_{r+1}', \dots , x_d' \rangle \right)}$.
For $i=1,2$, we denote the $i$-th factor by $\Scal_i$ and the corresponding projection by $p_i: \Scal \rightarrow \Scal_i$. The second factor $\Scal_2$ is  the affine formal scheme associated to the closed unit ball of dimension $d-r$. The first factor $\Scal_1$ is isomorphic to the affine formal scheme over $\kcirc$ associated to the polytopal domain $U_\Delta$ in ${\mathbb G}_m^r$, where $\Delta$ is the  simplex $\{u_1'+ \dots + u_r' \leq v(\pi)\}$ in $\rdop^r_+$. 

We will use the strata of the special fibre $\tilde{\Xcal}$ which were introduced at the end of \S 1. We will always normalize the formal open covering as in the following proposition. The reason will become clear in the construction of the skeleton.
\end{art}

\begin{prop} \label{normalization of covering}
Any formal open covering of $\Xcal'$ admits a refinement $\{\Ucal'\}$ by formal open subsets $\Ucal'$ as in \ref{strictly semistable} and which has the following  properties:
\begin{itemize}
\item[(a)] Every $\Ucal'$ is a formal affine open subscheme of $\Xcal'$.
\item[(b)] There is a distinguished stratum $S$ of $\tilde{\Xcal}'$ associated to $\Ucal'$ such that for any stratum $T$ of $\tilde{\Xcal}'$, we have $S \subset \overline{T}$ if and only if $\tilde{\Ucal}' \cap \overline{T} \neq \emptyset$.
\item[(c)]  $\tilde{\psi}^{-1}( \{\tilde{\mathbf 0}\} \times \tilde{\Scal}_2)$ is the stratum of $\tilde{\Ucal}'$ which is equal to $\tilde{\Ucal}' \cap S$ for the distinguished stratum $S$ from (b).
\item[(d)] Every stratum of $\tilde{\Xcal}'$ is the distinguished stratum of a suitable $\Ucal'$.
\end{itemize}
\end{prop}

\proof We start with the formal open covering $\{\Ucal'\}$ from \ref{strictly semistable}. We will refine it successively to get the claim. First, we may assume that  the covering is a refinement of the given formal open covering of $\Xcal'$. Let $\tilde{P}$ be any point of $\tilde{\Xcal}'$ and let $S$ be the stratum of $\tilde{\Xcal}'$ which contains $\tilde{P}$. There is an $\Ucal'$ with $\tilde{P} \in \tilde{\Ucal}'$. We remove from $\Ucal'$ the closure of all strata $T$ of $\tilde{\Xcal}'$ with $S \cap \overline{T} = \emptyset$. Note that the closure of a stratum in $\tilde{\Xcal}$ is a strata subset (see \cite{Ber4}, Proposition 2.1) and that the closures of two strata are either disjoint or one closure is contained in the other. Hence we get from $\Ucal'$ a formal open subset which contains $\tilde{P}$ and  which has property (b) for  our $S$. By passing to a formal affine open subset containing $\tilde{P}$, we get also (a). If we do this for every point $\tilde{P}$, we get a formal open subcovering with  properties (a),(b) and (d). So we may assume that the covering $\{
 \Ucal'\}$ satisfies (a),(b) and (d). We will show that this implies (c).

By \cite{Ber4}, Lemma 2.2, the restriction of $\tilde{\psi}$ to a stratum of $\tilde{\Ucal}'$ induces an \'etale morphism to a stratum of $\tilde{\Scal}$ and hence the preimage of a stratum of $\tilde{\Scal}$ is a stratum of $\tilde{\Ucal}'$. We conclude that $\tilde{\psi}^{-1}(\{\tilde{\mathbf 0}\} \times \tilde{\Scal}_2)$ is the union of $d-r$-dimensional strata $S_i$ of $\tilde{\Ucal}'$. Let $S$ be the distinguished stratum of $\tilde{\Xcal}'$ associated to $\Ucal'$. By (b), $S$ is contained in the closure of every $S_i$.  Since $\{\tilde{\mathbf 0}\} \times \tilde{\Scal}_2$ is a closed stratum of $\tilde{\Scal}$,  $\tilde{\psi}(S \cap \tilde{\Ucal}')$ is contained in $\{\tilde{\mathbf 0}\} \times \tilde{\Scal}_2$. This proves $S=S_i$ for some $i$. By dimensionality reasons, we get $S=S_i$ for every $i$  proving (c).  \qed



\begin{art} \rm \label{skeleton}
Next, we describe the 
{\it skeleton} of a strictly semistable formal scheme $\Xcal'$ over $\kcirc$. For details, we refer to  \cite{Ber4}, \S 4, \cite{Ber5}, \S 4,  and \cite{Gu4}, 9.1. 

We start with the model example $\Scal$ from \ref{strictly semistable}.  Replacing $x_0'$ by $\pi/(x_1' \dots x_r')$, every analytic function $f$ on $\Scal^{\rm an}$ has a unique representation as a convergent Laurent series of the form
$$f= \sum_{m_1,\dots,m_r \in \zdop} \sum_{m_{r+1},\dots,m_d \in \ndop} a_\mb (x_1')^{m_1} \dots (x_d')^{m_d}.$$
For every $\ub$ in the  simplex $\Delta:= \{\mathbf \ub \in \rdop^r_+ \mid u_1'+ \dots + u_r' \leq v(\pi)\}$, we get an element $\xi_{\ub} \in \Scal^{\rm an}$ using the bounded multiplicative seminorm
$$|f(\xi_{\ub})|:= \max_{\mb} |a_\mb| e^{-m_1 u_1- \dots - m_r u_r}.$$ 
We define the skeleton of $\Scal$ as  $\{\xi_\ub \mid \ub \in \Delta \}$. It is a closed subset of $\Scal^{\rm an}$ homeomorphic to $\Delta$. To omit the preference of the coordinate $x_0'$, it is better to identify  the skeleton of $\Scal'$ with the  simplex $\{u_0'+ \dots + u_r' = v(\pi)\}$ in $\rdop^{r+1}_+$. 

Next, we consider a formal open subset $\Ucal'$ of $\Xcal'$ as in Proposition \ref{normalization of covering}. Then the skeleton $S(\Ucal')$  of $\Ucal'$ is defined as the preimage of the skeleton of $\Scal$ with respect to the  morphism $\psi^{\rm an}$. It is a closed subset of the generic fibre $U'$ of $\Ucal'$. Using (b) of Proposition \ref{normalization of covering}, one can show that $\psi^{\rm an}$ induces a homeomorphism of $S(\Ucal')$ onto the skeleton of $\Scal$. Using the above, we may identify $S(\Ucal')$ again with the metrized simplex $\{u_0'+ \dots + u_r' = v(\pi)\}$ in $\rdop^{r+1}_+$. This is independent  of $\psi$ up to reordering  the coordinates $u_0', \dots, u_r'$. 

Finally, the skeleton $S(\Xcal')$ of $\Xcal'$ is the union of all skeletons $S(\Ucal')$. Berkovich has shown that $S(\Xcal')$ is a closed subset of the generic fibre $X'$ which depends only on the formal model $\Xcal'$, but neither on the choice of the formal covering $\{\Ucal'\}$ nor on the choice of the \'etale morphisms $\psi$.

The skeleton $S(\Xcal')$ has a canonical  structure as an abstract metrized simplicial set which reflects the incidence relations between the strata of $\tilde{\Xcal}'$: For every stratum $S$ of codimension $r$, there is a {\it canonical simplex} $\Delta_S$ in $\xskel$ defined in the following way. We choose a formal affine open subset $\Ucal'$ as in Proposition \ref{normalization of covering} such that $S$ is the  distinguished stratum associated to $\Ucal'$. Then we define $\Delta_S$ as the skeleton of $\Ucal'$. It is easy to see that $\Delta_S$ does not depend on the choice of $\Ucal'$ and hence we may identify $\Delta_S$ with the  simplex $\{u_0'+ \dots + u_r' = v(\pi)\}$ in $\rdop^{r+1}_+$.  The canonical simplices have the properties:
\begin{itemize}
\item[(a)] The canonical simplices $(\Delta_S)_{S \in \str(\tilde{\Xcal}')}$ cover $\xskel$. 
\item[(b)]  For $S \in \xstr$, the map $T \mapsto \Delta_T$ gives a bijective order reversing correspondence between  $T \in \xstr$ with $S \subset \overline{T}$ and  closed faces of $\Delta_S$.
\item[(c)] For $R, S \in \xstr$, $\Delta_{R} \cap \Delta_S$ is the union of all $\Delta_T$ with $T \in \xstr$ and $\overline{T} \supset R \cup S$.
\end{itemize}
There is a continuous map $\Val: X' \rightarrow \xskel$ which restricts to the identity on $\xskel$. It is enough to define it for $p \in U'$, where $U'$ is the generic fibre of a formal affine open subset $\Ucal'$ as above. Using the identification $\Delta_S= \{u_0'+ \dots + u_r' = v(\pi)\}$, we set  $\Val(p):=(-\log \circ p(\psi^*(x_0')), \dots , - \log \circ p(\psi^*(x_r'))) \in \Delta_S$. By \cite{Ber4}, Theorem 5.2, the map $\Val$ gives a proper strong deformation retraction of $X'$ to the skeleton $S(\Xcal')$.
\end{art}

\begin{art} \rm \label{subdivision}
It would be tempting to call the family of canonical simplices a polytopal decomposition of $\xskel$. However, we note that the family is not necessarily face to face, only the weaker property \ref{skeleton}(c) holds instead. 

In the following, we consider now a $\Gamma$-rational {\it polytopal subdivision} $\Dcal$ of   $\xskel$. This means that $\Dcal$ is a family of $\Gamma$-rational polytopes, each contained in a canonical simplex, such that 
$\Dcal \cap \Delta_S:= \{\Delta \in \Dcal \mid \Delta \subset \Delta_S\}$
is a polytopal decomposition of $\Delta_S$ for every $S \in \xstr$.
\end{art}

\begin{prop} \label{refinement and formal analytic structure}
There is a coarsest formal analytic structure $\X''$ on $X'$ which refines $(\Xcal')^{\rm f-an}$ in such a way that $\Val^{-1}(\Delta)$ is a formal open subset for every $\Delta \in \Dcal$. 
\end{prop}

\proof Let $S \in \xstr$ and let $U'$ be the generic fibre of a set $\Ucal'$ as in Proposition \ref{normalization of covering}. We note that such sets $U'$, for varying $S$, form a formal affinoid atlas of $(\Xcal')^{\rm f-an}$. To prove the claim, it is enough to show that the sets
\begin{equation} \label{formal affinoid atlas of U'}
\left(U' \cap \Val^{-1}(\Delta) \right)_{\Delta \in \Dcal\cap \Delta_S}
\end{equation}
define a formal affinoid atlas on $U'$. The polytope $\Delta \in \Dcal\cap \Delta_S$ is given by finitely many inequalities of the form $\mb \cdot \ub'+v(\lambda) \geq 0$ for some $\mb \in \zdop^{r+1}$ and $\lambda \in \kdop^\times$. In terms of the semistable coordinates $x_0', \dots ,x_r'$ of $U'$, the subset $U' \cap \Val^{-1}(\Delta)$ is given by finitely many inequalities of the form $|\lambda \psi^*(\xb')^\mb| \leq 1$ and hence it is an affinoid subdomain of $U'$. This description yields easily that \eqref{formal affinoid atlas of U'} is a formal affinoid atlas of $U'$ proving the claim. \qed

\begin{rem} \label{local description by base change} \rm 
Let $\Ucal'$ be  a formal open subset of $\Xcal'$ as in Proposition \ref{normalization of covering} with \'etale morphism $\psi: \Ucal' \rightarrow \Scal$. The generic fibre $U'$ is a formal open subset of the formal analytic variety $\X''$ from Proposition \ref{refinement and formal analytic structure} and we write suggestively $U' \cap \X''$ for the formal analytic structure on $U'$ induced by $\X''$. 

Let $S$ be the distinguished stratum of $\tilde{\Xcal}'$ associated to ${\Ucal}'$.  
The first factor $\Scal_1$ from \ref{strictly semistable} is the formal scheme over $\kcirc$ associated to the polytopal domain $\val^{-1}(\Delta_S)$ in $\{\xb' \in {\mathbb G}_m^{r+1} \mid x_0' \cdots x_r'=v(\pi)\}$. 

The polytopal decomposition $\Dcal \cap \Delta_S$ of $\Delta_S$ leads to a formal analytic refinement of $\Scal_1^{\rm f-an}$ inducing an admissible formal scheme $\Scal_1'$ over $\kcirc$ and a canonical morphism $\iota_1:\Scal_1' \rightarrow \Scal_1$ extending the identity from the generic fibre. By base change, we get a morphism $\iota: \Scal' \rightarrow \Scal$ with the same property. Note that $\Scal'=\Scal_1' \times \Scal_2$ has reduced special fibre (see \ref{admissible formal schemes}). 

Since the base change $\psi': \Ucal'' \rightarrow \Scal'$ of $\psi$ with respect to $\iota$ is \'etale,  $\Ucal''$ has also reduced special fibre (see \cite{EGA IV}, 17.5.7). By \ref{admissible formal schemes}, $\Ucal''$ is an admissible formal scheme over $\kcirc$ associated to a formal analytic variety. By the proof of Proposition \ref{refinement and formal analytic structure}, the latter is $U' \cap \X''$ and hence $\Ucal'' = (U' \cap \X'')^{\rm f-sch}$. 
\end{rem}

In the following, $\Xcal''$ denotes the admissible formal $\kcirc$-scheme associated to $\X''$ and hence we may identify the special fibre $\tilde{\X}''$ with the reduction $\tilde{\Xcal}''$ (see \ref{admissible formal schemes}). Since $\Dcal$ is a polytopal subdivision of the skeleton, the identity is a formal analytic morphism $\X'' \rightarrow (\Xcal')^{\rm f-an}$ and hence we get a unique morphism $\iota':\Xcal'' \rightarrow \Xcal'$ extending the identity on the generic fibre.

Recall that the order on the strata is given by inclusion of closures. Similarly, we define an order on the open faces of $\Dcal$.

\begin{prop} \label{orbits and refinement}
Let $\X''$ be the formal analytic variety associated to $\Dcal$ as described in Proposition \ref{refinement and formal analytic structure}. Then there is a bijective  correspondence between open faces $\tau$ of $\Dcal$ and strata $R$ of $\tilde{\X}''$, given by
\begin{equation} \label{orbit correspondence}
R= \pi\left(\Val^{-1}(\tau)\right), \quad \tau= \Val\left(\pi^{-1}(Y)\right),
\end{equation}
where $\pi: X' \rightarrow \tilde{\X}''$ is the reduction map and $Y$ is any non-empty subset of $R$.
\end{prop}

\proof Let $\tau$ be an open face of $\Dcal$. We have to prove that $R:=\pi\left(\Val^{-1}(\tau)\right)$ is a stratum of $\tilde{\X}''$. There is a unique $S \in \xstr$ such that $\tau$ is contained in $\relint(\Delta_S)$. 
We choose a formal affine open subset $\Ucal'$ as in Proposition \ref{normalization of covering} such that $S$ is the distinguished stratum associated to $\Ucal'$. Note that strata are compatible with localization and hence we may assume $\Xcal' = \Ucal'$. By Remark \ref{local description by base change}, we have a Cartesian diagram  of admissible formal schemes over $\kcirc$
\begin{equation*}
\begin{CD} 
\Xcal'' @>{\psi'}>> \Scal'  @>{p_1'}>>  \Scal_1' \\
@VV{\iota'}V  @VV{\iota}V  @VV{\iota_1}V \\
\Xcal'  @>{\psi}>>  \Scal @>{p_1}>>  \Scal_1
\end{CD}
\end{equation*}
with $\psi$ and $\psi'$ \'etale. Let $\psi_1:=p_1 \circ \psi$ and $\psi_1':=p_1' \circ \psi'$. 

The idea of the proof is to use $\psi_1$ to reduce the claim to the corresponding statement for the polytopal domain $\Scal_1$ in $\mathbb G_m^r$. We describe this result here in terms of the torus $\mathbb G_m^{r+1}$ and  in terms of the  valuation map 
$$\val: (\mathbb G_m^{r+1})^{\rm an} \rightarrow \rdop^{r+1},\quad p \mapsto \left(-\log \circ p(x_0'), \dots, - \log \circ p(x_r')\right)$$ 
to omit the preference of the first coordinate. 
By \cite{Gu4}, Propositions 4.4 and 4.7, there is a bijective correspondence between open faces $\sigma$ of $\Dcal$ (which is a polytopal decomposition of $\Delta_S=\{{\mathbf u}\in \rdop^{r+1}\mid u_0'+ \dots + u_r' = v(\pi)\}$ by the assumption $\Ucal'=\Xcal'$) and strata $T_1'$ of $\tilde{\Scal}_1'$, given by
\begin{equation} \label{correspondence for factor}
T_1'= \pi \left( \val^{-1}(\sigma) \cap \Scal_1^{\rm an}\right), \quad
\sigma = \val\left(\pi^{-1}(T_1')\right),
\end{equation}
where $\pi: \Scal_1^{\rm an} =({\Scal_1'})^{\rm an} \rightarrow \tilde{\Scal}_1'$ denotes the reduction map. In fact, one can replace $T_1'$ in the second formula of \eqref{correspondence for factor} by any non-empty subset of $T_1'$. To see this, note that the formal affinoid subtorus $D=\{|x_0|= \dots =|x_r|=1,\, x_0 \cdots x_r=1\}$ of ${\mathbb G}_m^{r+1}$ acts on $\Scal_1^{\rm an}$ and this extends to an action of the formal torus on $\Scal_1'$. The strata of $\Scal_1'$ are the same as the torus orbits. We conclude that $D$ acts transitively on the set $\{\pi^{-1}(\tilde{P})\mid \tilde{P} \in T_1'(\ktilde)\}$. Since the map $\val$ is invariant under the $D$-action, we conclude that $\val(\pi^{-1}(T_1'))=\val(\pi^{-1}(\tilde{P}))$ for any $\ktilde$-rational point $\tilde{P}$ of $T_1'$. Note that we may use base extension to make any non-closed point rational; therefore since the map $\val$ is invariant under base extension we conclude that $\val(\pi^{-1}(T_1'))=\val(\pi^{-1}(\tilde{P}))$ holds for {\it any} (i.e. not necessarily closed) point $\tilde{P}$ of $T_1'$. This proves the second formula in \eqref{correspondence for factor} with $T_1'$ replaced by any
  non-empty subset.

Now let $T_1'$ be the stratum of $\tilde{\Scal}_1'$ corresponding to the given open face $\tau$. Obviously, $T':=(\tilde{p}_1')^{-1}(T_1')=T_1' \times \tilde{\Scal}_2$ is a stratum of $\tilde{\Scal}'=\tilde{\Scal}_1' \times \tilde{\Scal}_2$. 
We would like to prove that the preimage of $T'$ with respect to $\tilde{\psi}'$ is equal to $R$. Using \eqref{correspondence for factor}, we first note that
\begin{equation} \label{psi1-formel}
(\tilde{\psi}')^{-1}(T') = (\tilde{\psi}_1')^{-1} \left(\pi \left( \val^{-1}(\tau) \cap \Scal_1^{\rm an}\right)\right). 
\end{equation}
Next, we prove the following formula
\begin{equation} \label{vertauschen von psi1 mit pi}
(\tilde{\psi}_1')^{-1} \left(\pi \left( \val^{-1}(\tau) \cap \Scal_1^{\rm an}\right)\right)
=\pi \left( (\psi_1^{\rm an})^{-1} (\val^{-1}(\tau) ) \right) .
\end{equation}
The inclusion "$\supset$"  follows immediately from $\pi \circ \psi_1^{\rm an}=\tilde{\psi}_1' \circ \pi$. Here, we have used that $\psi_1^{\rm an}=(\psi_1')^{\rm an}$. To prove the reverse inclusion, let us choose a point $\tilde{x}' \in (\tilde{\psi}_1')^{-1} (\pi \left( \val^{-1}(\tau) \cap \Scal_1^{\rm an}\right))$. The reduction map is surjective, hence there is $x'\in X'$ with $\pi(x')=\tilde{x}'$. By assumption, there is $z \in \val^{-1}(\tau) \cap \Scal_1^{\rm an}$ with 
$$\pi(z)=\tilde{\psi}_1'(\tilde{x}')=\tilde{\psi}_1'(\pi(x'))=\pi \circ \psi_1^{\rm an}(x').$$
By \eqref{correspondence for factor}, we  get $\pi \circ\psi_1^{\rm an}(x') \in T_1'$. An  application of \eqref{correspondence for factor} shows $\val(\psi_1^{\rm an}(x')) \in \tau$. We conclude that $\tilde{x}'=\pi(x')\in \pi ( (\psi_1^{\rm an})^{-1} (\val^{-1}(\tau) ) )$ proving \eqref{vertauschen von psi1 mit pi}.

Using \eqref{psi1-formel} and \eqref{vertauschen von psi1 mit pi}, we get finally the desired relation between  $R$ and $T'$:
\begin{equation} \label{Urbildrelation}
(\tilde{\psi}')^{-1}(T') = \pi \left( (\psi_1^{\rm an})^{-1} (\val^{-1}(\tau) ) \right)
=\pi\left(\Val^{-1}(\tau)\right) = R.
\end{equation}
By \cite{Ber4}, Lemma 2.2, the preimage of the stratum $T'$ with respect to the \'etale morphism $\tilde{\psi}'$ is the union of strata of the same dimension. This argument was already used in the proof of Proposition \ref{normalization of covering}. To prove that $R$ is a stratum, it is enough to show that $R$ is irreducible. Note that $\tilde{\iota}_1(T_1')=\{\tilde{\mathbf 0}\}$ in $\tilde{\Scal}_1$ and hence
\begin{equation} \label{preimage of torus orbit}
(\tilde{\psi}')^{-1}(T') = T' \times_{\tilde{\Scal}} \tilde{\Xcal}' = T' \times_{\{\tilde{\mathbf 0}\}\times \tilde{\Scal}_2} \tilde{\psi}^{-1}(\{\tilde{\mathbf 0}\}\times \tilde{\Scal}_2) \cong \mathbb G_m^{r-t} \times S,
\end{equation}
where $t:= \dim(\tau)$. Here, we have used that $T_1'$ is an $(r-t)$-dimensional torus orbit  and that $S = \tilde{\psi}^{-1}(\{\tilde{\mathbf 0}\}\times \tilde{\Scal}_2)$ (see Proposition \ref{normalization of covering} and \cite{Gu4}, Proposition 4.4). We conclude that $R=(\tilde{\psi}')^{-1}(T')$ is irreducible proving that $R \in \str(\tilde{\X}'')$.

Since the open faces of $\Dcal$ form a covering of the skeleton $\xskel$, we conclude that every $R \in \str(\tilde{\X}'')$ has the form $R=\pi(\Val^{-1}(\tau))$ for an open face $\tau$ of $\Dcal'$. 

It remains to show that $\tau$ may be reconstructed from $R$ by the second formula in \eqref{orbit correspondence}.  By the same argument as used in the paragraph after 
\eqref{correspondence for factor}, it is enough to prove this for $Y=\{\tilde{y}\}$ for any $\ktilde$-rational point $\tilde{y}$ of $R$. Since $\psi'$ is \'etale, the formal fibre $X_+'(\tilde{y}):=\pi^{-1}(\tilde{y})$ is isomorphic to the formal fibre over  $\tilde{z}:=\tilde{\psi}'(\tilde{y})$ in $\Scal'$ (see \cite{Gu4}, Proposition 2.9). For $\tilde{z}_1:=\tilde{p}_1(\tilde{z})$, we get the following isomorphism of formal fibres:
\begin{equation} \label{factorization of formal fibre}
X_+'(\tilde{y}) \cong (\Scal_1')^{\rm an}_+(\tilde{z}_1) \times (\Scal_2)^{\rm an}_+(\tilde{\mathbf 0}).
\end{equation}
The $(d-r)$-dimensional ball $\Scal_2^{\rm an}$ does not contribute to $\Val$. 
Using the analogue of the claim for the polytopal domain $\Scal_1^{\rm an}$ deduced after \eqref{correspondence for factor}, we get
$$\Val\left(X_+'(\tilde{y})\right)=\val\left((\Scal_1')^{\rm an}_+(\tilde{z}_1)\right) =\val\left(\pi^{-1}(T_1')\right) = \tau.$$
This proves the second formula in \eqref{orbit correspondence}. \qed

\begin{rem} \label{preimage of strata wrt psi} \rm
Let $\Ucal'$ be a formal affine open subset of  $\Xcal'$ as in Proposition \ref{normalization of covering} and let $\psi:\Ucal' \rightarrow  \Scal = \Scal_1 \times \Scal_2$ be the \'etale morphism from \ref{strictly semistable}. Let us consider the composition $\psi_1:\Ucal' \rightarrow \Scal_1$  of the first projection with $\psi$ and let $\psi_1':\Ucal'' \rightarrow \Scal_1'$ be the base change of $\psi_1$ induced by the polytopal decomposition $\Dcal$. We have seen in the above proof that the preimage of any stratum  of $\tilde{\Scal}_1'$ with respect to  $\tilde{\psi}_1'$ is a stratum  of $\tilde{\Ucal}''$. 
\end{rem}

Recall that $\Xcal''=(\X'')^{\rm f-sch}$ and that we have a canonical morphism $\iota':\Xcal'' \rightarrow \Xcal'$ extending the identity on the generic fibre.

\begin{cor} \label{strata corollary}
Let $R \in \str(\tilde{\X}'')$ with corresponding open face $\tau$ of $\Dcal$. 
\begin{itemize}
 \item[(a)] $\dim(\tau)= \codim(R,\tilde{\X}'')$.
 \item[(b)] $S:=\tilde{\iota}'(R) \in \xstr$.
 \item[(c)] $R \stackrel{\tilde{\iota}'}{\rightarrow}S$ is a fibre bundle with fibre $(\mathbb G_m)_\ktilde^{\dim(R)-\dim(S)}$.
 \item[(d)] Every stratum of $\tilde{\X}''$ is smooth.
 \item[(e)] The closure $\overline{R}$ is the union of the strata of $\tilde{\X}''$ corresponding to the open faces $\sigma$ of $\Dcal$ with $\tau \subset \overline{\sigma}$.
 \item[(f)] For open faces $\tau_1$, $\tau_2$ of $\Dcal$ with corresponding strata $R_1$, $R_2$ of $\tilde{\X}''$, we have $\overline{\tau_1} \subset \overline{\tau_2}$ if and only if $\overline{R_1} \supset \overline{R_2}$.
 \item[(g)] For an irreducible component $Y$ of $\tilde{\X}''$, let $\xi_Y$ be the unique point of $X'$ with reduction equal to the generic point of $Y$. Then $Y \mapsto \xi_Y$ is a bijection between the irreducible components of $\tilde{\X}''$ and the vertices of $\Dcal$.
\end{itemize}
\end{cor}

\proof We use the proof of Proposition \ref{orbits and refinement}. By \cite{Gu4}, Proposition 4.4, we have $\dim(\tau)=\codim(T_1',\tilde{\Scal}_1')$. Using that $R$ is locally equal to $(\tilde{\psi}_1')^{-1}(T_1')$ and the smoothness of $\tilde{\psi}_1'$, we get (a). Let $S \in \xstr$ with $\tau \subset \relint(\Delta_S)$. By (a) and \eqref{preimage of torus orbit}, we deduce (b) and (c). Since $S$ is smooth by Proposition \ref{normalization of covering}(c), we get (d) from (c).

Since strata are compatible with localization, it is enough to prove (e) in case of $\Xcal'=\Ucal'$ for a formal affine $\Ucal'$ as in Proposition \ref{normalization of covering} such that $S$ is the distinguished stratum of $\tilde{\Xcal}'$ associated to $\Ucal'$. Since $\tilde{\psi}_1'$ is flat, we get $\overline{R}= (\tilde{\psi}_1')^{-1}(\overline{T_1'})$. By \cite{Gu4}, Remark 4.8, $\overline{T_1'}$ is the union of all strata of $\tilde{\Scal}_1'$ corresponding to the open faces $\sigma$ of $\Dcal$ with $\overline{\sigma} \supset \overline{\tau}$. If we take preimages of this decomposition, we get (e). Note that (f) is a consequence of (e). 

By using the unique dense stratum of an irreducible component, it follows from (a) and Proposition \ref{orbits and refinement} that the map $Y \mapsto \Val(\xi_Y)$ is a bijection between the irreducible components of $\tilde{\X}''$ and the vertices of $\Dcal$. To prove (g), it remains to see that $\xi_Y \in \xskel$. There is a formal affine open subset $\Ucal'$ of $\Xcal'$ as in Proposition \ref{normalization of covering} with $Y \cap \Ucal' \neq \emptyset$ and an \'etale morphism $\psi:\Ucal' \rightarrow \Scal=\Scal_1 \times \Scal_2$. By Remark \ref{preimage of strata wrt psi}, there is an irreducible component $Z$ of $\tilde{\Scal}'$ such that $\psi(\xi_Y)=\xi_Z$. Since $\psi$ is \'etale, it is enough to prove that $\xi_Z \in S(\Scal)$ (see \cite{Ber5}, Corollary 4.3.2). 
Since $Z = Z_{\ub} \times \tilde{\Scal}_2$ for the irreducible component $Z_{\ub}$ of $\tilde{\Scal}_1$ corresponding to the vertex $\ub=\Val(\xi_Y)$ of $\Dcal$ (see \cite{Gu3}, Proposition 4.7), it is easy to see that the point $\xi_{\ub}$ from \ref{skeleton} reduces to the generic point of $Z$ and hence we get $\xi_{\ub}=\xi_Z$ proving the claim. \qed

\begin{art} \rm \label{skeleton and alteration}
For the remaining part of this section, we fix the following situation: Let $A$ be an abelian variety over $\kdop$ with uniformization $E$ such that $A^{\rm an}=E/M$ as in \ref{Raynaud extension}. We recall that $M$ is a discrete subgroup of $E(\kdop)$ such that $\val:E \rightarrow \rdop^n$ maps $M$ isomorphically onto a complete lattice $\Lambda$ of $\rdop^n$.

We assume that we have a morphism $\varphi_0:\Xcal' \rightarrow \Acal_0$, where $\Xcal'$ is still a strictly semistable scheme over $\kcirc$ and $\Acal_0$ is the Mumford model of $A$ associated to a $\Gamma$-rational polytopal decomposition $\overline{\Ccal_0}$ of $\rtor$. Let  $f:X' \rightarrow A$ be the generic fibre of $\varphi_0$. 
\end{art}

\begin{prop} \label{f_aff}
There is a unique map $\overline{f}_{\rm aff}:\xskel \rightarrow \rtor$ with $\overline{f}_{\rm aff} \circ \Val = \valbar \circ f$ on $X'$. The map $\overline{f}_{\rm aff}$ is continuous. For every $S \in \xstr$, the restriction of $\overline{f}_{\rm aff}$ to the canonical simplex $\Delta_S$ is an affine map and there is a unique $\Deltabar \in \overline{\Ccal_0}$ with $\overline{f}_{\rm aff}(\relint(\Delta_S))\subset \relint(\Deltabar)$.
\end{prop}

\proof We recall the construction of $\Acal_0$ given in \ref{Mumford's construction}. Let $\Vcal$ be a formal affine open subset of the formal abelian scheme $\Bcal$ which trivializes the
 Raynaud extension \eqref{Raynaud extension 1} of $A$. For the generic fibre $V$ of $\Vcal$ and $\Delta \in {\Ccal_0}$, we get a formal affinoid subdomain $U_{V,\Delta}$ of $E$ with associated affine formal schemes $\Ucal_{V,\Delta}$. With varying $V$ and $\Delta$, we get a formal affinoid atlas on $E$ with associated $\kcirc$-model $\Ecal_0$ of $E$ which is covered by the formal open subsets $\Ucal_{V,\Delta}$. By passing to the quotient by $M$, we get $\Acal_0=\Ecal_0/M$ and the quotient morphism maps $\Ucal_{V,\Delta}$ isomorphically onto the formal open chart $\Ucal_{[V,\Delta]}$ of $\Acal_0$.

There is a formal open covering $\{\Ucal'\}$ of $\Xcal'$ as in Proposition  \ref{normalization of covering} such that for any $\Ucal'$ of the covering, there are $V,\Delta$ as above with $\Ucal' \subset \varphi_0^{-1}(\Ucal_{[V,\Delta]})$.  We denote the generic fibre of $\Ucal'$ by $U'$. By construction, there is a unique lift $F:U' \rightarrow U_{V,\Delta}$ of $f$. Now we use the coordinates $x_1, \dots, x_n$ of the torus $T$ from the Raynaud extension \eqref{Raynaud extension 2} of $A$. They are defined on $U_{V,\Delta}$ by using the trivialization $U_{V,\Delta} \cong V \times U_\Delta$ from \ref{Raynaud trivializations} for the polytopal domain $U_\Delta$ of $T$. Note that $F^*(x_i)$ is a unit of $U'$. By \cite{Gu4}, Proposition 2.11, there are $u_i \in \Ocal(\Ucal')^\times$, $\mb_i \in \zdop^{r+1}$ and $\lambda_i \in \kdop^\times$ with
\begin{equation} \label{x in terms of x'}
F^*(x_i) = \lambda_i u_i \psi^*(\xb')^{\mb_i}
\end{equation}
for $i=1, \dots ,n$, where $\psi:\Ucal' \rightarrow \Scal=\Scal_1 \times \Scal_2$ is again the \'etale morphism and $\xb'=(x_0', \dots, x_r')$ are the semistable coordinates from \ref{strictly semistable}. 
Let $S \in \xstr$ be the distinguished stratum  of $\tilde{\Xcal}'$ associated to $\Ucal'$. Then the canonical simplex $\Delta_S$ may be identified with the simplex $\{u_0'+ \dots + u_r'=v(\pi)\}$ in  $\rdop_+^{r+1}$ and we define $f_{\rm aff}:\Delta_S \rightarrow \rdop^n$ by 
\begin{equation} \label{definition of f_aff}
f_{\rm aff}(u_0', \dots, u_r'):= \left(\mb_i \cdot \ub' + v(\lambda_i) \right)_{i=1, \dots, n}.
\end{equation}
We note that this definition depends only on $|\lambda_i|$, $|x_i'|$ and $\mb_i$, hence it is independent of the trivialization of $U_{V,\Delta}$. If we change $V$ and $\Delta$, then the new lift is obtained from $F$ by an $M$-translation. We deduce easily that the locally defined maps $f_{\rm aff}$ induce a well-defined map $\overline{f}_{\rm aff}:S(\Xcal') \rightarrow \rtor$. All the claimed properties follow from the construction and uniqueness is clear from surjectivity of $\Val$. \qed

\begin{rem} \label{strata are preserved} \rm
More generally, Berkovich (\cite{Ber5}, Corollary 6.1.2) has shown that a morphism between strongly non-degenerate pluristable formal schemes over $\kcirc$ induces a piecewise linear map between the skeletons. Now Proposition \ref{f_aff} describes precisely the domain of affineness and we will see in Remark \ref{strictly pluristable} that this holds also if $\Xcal'$ is a strongly non-degenerate strictly pluristable formal scheme over $\kcirc$. 

By Proposition \ref{strata on Mumford models} and Proposition \ref{orbits and refinement}, we conclude easily that every stratum of $\tilde{\Xcal}'$ is mapped into a stratum of $\tilde{\Acal}_0$. This will be proved in a more general context in Lemma \ref{lift of phi}. The preservance of strata is a key fact which  will allow us to describe canonical measures on $X=f(X')$ in terms of the skeleton of $\Xcal'$.
 \end{rem}

\begin{prop} \label{formal refinement for alteration}
Let $\Ccalbar$ be a $\Gamma$-rational polytopal subdivision of $\overline{\Ccal_0}$ with associated Mumford model $\Acal$ of $A$. Then $\left(\Xcal' \times_{\Acal_0} \Acal   \right)^{\rm f-an}$ is the formal analytic variety $\X''$ from Proposition \ref{refinement and formal analytic structure} associated to the $\Gamma$-rational subdivision $\Dcal$ of $\xskel$ given by 
$$\Dcal:=\{\Delta_S \cap \overline{f}_{\rm aff}^{-1}(\sigmabar) \mid S \in \xstr, \sigmabar \in \Ccalbar\}.$$
\end{prop}

\proof We will use the notation from the  proof of Proposition \ref{f_aff}. Let $\sigma \in \Ccal$ be contained in $\Delta \in \Ccal_0$, let $V$ be the generic fibre of a formal affine open subset $\Vcal$ of $\Bcal$ and let $\Ucal'$ be a formal affine open subset of $\Xcal'$ as in Proposition \ref{normalization of covering} with $\varphi_0(\Ucal') \subset \Ucal_{[V,\Delta]}$. Then the sets
$$\left( \Ucal' \times_{\Ucal_{[V,\Delta]}} \Ucal_{[V,\sigma]} \right)^{\rm an}
= U' \times_{A^{\rm an}} U_{[V,\sigma]}$$
form a formal affinoid atlas of $\left(\Xcal' \times_{\Acal_0} \Acal \right)^{\rm f-an}$. We have $U_{[V,\sigma]}=U_{[V,\Delta]} \cap  \valbar^{-1}(\sigmabar)$ and hence Proposition \ref{f_aff} yields 
$$U' \times_{A^{\rm an}} U_{[V,\sigma]}= U' \cap f^{-1}\left( \valbar^{-1}(\sigmabar) \right) = U' \cap \Val^{-1} \left(\overline{f}_{\rm aff}^{-1}(\sigmabar) \right).$$
Let $S$ be the distinguished stratum of $\tilde{\Xcal}'$ associated to $\Ucal'$. Therefore, we have  $\Val(U')=\Delta_S$ and we deduce 
$$\left( \Ucal' \times_{\Ucal_{[V,\Delta]}} \Ucal_{[V,\sigma]} \right)^{\rm an} = U' \cap \Val^{-1}(\sigma')$$
for the polytope $\sigma':= \overline{f}_{\rm aff}^{-1}(\sigmabar) \cap \Delta_S \in \Dcal$. These sets form the formal affinoid atlas \eqref{formal affinoid atlas of U'} for $\X''$ proving $\X''=\left(\Xcal' \times_{\Acal_0} \Acal   \right)^{\rm f-an}$. \qed

\begin{prop} \label{minimal formal analytic structure}
We keep the above assumptions and notation. Then $\X''$ is the coarsest formal analytic variety on $X'$ such that  $f:X' \rightarrow A^{\rm an}$ induces a formal analytic morphism $\phi:\X'' \rightarrow \Acal^{\rm f-an}$. If $R$ is the stratum of $\tilde{\X}''$ corresponding to the open face $\tau$ of $\Dcal$, then $\tilde{\phi}(R)$ is contained in the stratum of $\tilde{\Acal}$ corresponding to the unique open face $\sigmabar$ of $\Ccalbar$ with $\overline{f}_{\rm aff}(\tau) \subset \sigmabar$.
\end{prop}

\proof The first claim is clear by construction. By definition of $\Dcal$, there is an open face $\sigmabar$ of $\Ccalbar$ with $\overline{f}_{\rm aff}(\tau) \subset \sigmabar$. If $\pi$ denotes the reduction map, then we get
$$\tilde{\phi}(R)=\tilde{\phi}(\pi(\Val^{-1}(\tau)))=\pi(f(\Val^{-1}(\tau))) \stackrel{\ref{f_aff}}{\subset} \pi(\valbar^{-1}(\sigmabar)).$$
By Proposition \ref{strata on Mumford models}, we deduce that $\tilde{\phi}(R)$ is contained in the stratum of $\tilde{\Acal}$ corresponding to $\sigmabar$. \qed

\vspace{3mm}
Again, let $\Xcal''$ be the admissible formal $\kcirc$-scheme associated to the formal analytic variety $\X''$ from Propositions \ref{formal refinement for alteration} and \ref{refinement and formal analytic structure}. The following commutative diagram gives an overview of the occuring canonical morphisms of admissible formal schemes, where $\Ecal_0$ (resp. $\Ecal$) is the $\kcirc$-model of the uniformization $E$ associated to $\Ccal_0$ (resp. $\Ccal$) as in \ref{Mumford's construction} and where the vertical maps extend the identity on the generic fibre. 

\begin{equation}  \label{com diagram}
\begin{CD} 
\Xcal'' @>\varphi>> \Acal @<{p}<< \Ecal\\
@VV{\iota'}V    @VV{\iota_0}V  @VV{\iota}V\\
\Xcal' @>{\varphi_0}>> \Acal_0 @<{p_0}<<\Ecal_0
\end{CD}
\end{equation}


\begin{lem} \label{lift of phi}
Let $R \in \str(\tilde{\Xcal}'')$. Then $S:=\tilde{\iota}'(R)$ is 
a stratum  of $\tilde{\Xcal}'$. The restricted morphism $\tilde{\varphi}_0:\overline{S} \rightarrow \tilde{\Acal}_0=\tilde{\Ecal}_0/M$ has a lift $\tilde{\Phi}_0:\overline{S} \rightarrow \tilde{\Ecal}_0$, unique up to the $M$-action on $\tilde{\Ecal}_0$. Moreover, there is a unique lift  $\tilde{\Phi}:\overline{R}\rightarrow \tilde{\Ecal}$ of $\tilde{\varphi}:\overline{R} \rightarrow \tilde{\Acal}=\tilde{\Ecal}/M$ with $\tilde{\Phi}_0 \circ \tilde{\iota}' = \tilde{\iota} \circ \tilde{\Phi}$ on $\overline{R}$.
\end{lem}

\proof The first claim was proved in Corollary \ref{strata corollary}. The proof of the remaining claims follows standard arguments from the theory of coverings (applied to the quotient maps $\tilde{p}_0:\tilde{\Ecal}_0 \rightarrow \tilde{\Acal}_0=\tilde{\Ecal}_0/M$ and 
$\tilde{p}:\tilde{\Ecal} \rightarrow \tilde{\Acal}=\tilde{\Ecal}/M$):

Let $Y$ be an irreducible component of $\tilde{\Acal}_0$ with $\tilde{\varphi}_0(\overline{S}) \subset Y$. By Proposition \ref{strata on Mumford models}, $Y$ corresponds to a vertex $\ubb$ of $\overline{\Ccal_0}$ and $\tilde{p}_0^{-1}(Y)$ is the disjoint union of the irreducible components $Y_\ub$ of $\tilde{\Ecal_0}$ associated to the vertices 
$\ub$ of $\Ccal_0$ with residue class $\ubb \in \rtor$. Moreover, $Y_\ub$ is mapped isomorphically onto $Y$ by  $\tilde{p}_0$. Using composition with the inverse $Y \rightarrow Y_\ub$ of this isomorphism, we get the desired lift $\tilde{\Phi}_0$ of the restriction of $\tilde{\varphi}_0$ to ${\overline{S}}$. Uniqueness up to the $M$-action  on $\tilde{\Ecal_0}$ is obvious. Similarly, we get a lift $\tilde{\Phi}$ of the restriction of $\tilde{\varphi}$ to ${\overline{R}}$ by working with $\Ccalbar$ instead of $\overline{\Ccal_0}$. The lift $\tilde{\Phi}$ is also unique up to the $M$-action on $\tilde{\Ecal}$.

It remains to prove that $\tilde{\Phi}_0$ determines $\tilde{\Phi}$ uniquely by the condition $\tilde{\Phi}_0 \circ \tilde{\iota}' = \tilde{\iota} \circ \tilde{\Phi}$. Let us choose $\tilde{x}' \in R$ and let $\tilde{x}:=\tilde{\iota}'(\tilde{x}')$. Note that the lift $\tilde{\Phi}_0$ is uniquely determined by choosing an element $\tilde{y} \in \tilde{p}_0^{-1}(\tilde{\varphi}_0(\tilde{x}))$ if we require $\tilde{\Phi}_0(\tilde{x})=\tilde{y}$. Similarly, $\tilde{\Phi}$ is determined by $\tilde{\Phi}(\tilde{x}')=\tilde{y}'$ for some $\tilde{y}'  \in \tilde{p}^{-1}(\tilde{\varphi}(\tilde{x}))$. Since $M$ acts faithfully and transitively on $\tilde{p}_0^{-1}(\tilde{\varphi}_0(\tilde{x}))$ (resp. on $\tilde{p}^{-1}(\tilde{\varphi}(\tilde{x}))$), there is a unique $\tilde{y}'$ with $\tilde{\iota}(\tilde{y}')=\tilde{y}$. Note that $\tilde{p}_0$ and $\tilde{p}$ are local isomorphisms. Since $\tilde{\varphi}_0 \circ \tilde{\iota}'=\tilde{\iota}_0 \circ\tilde{\varphi}$, this lifts to the identity
  $\tilde{\Phi}_0 \circ \tilde{\iota}' = \tilde{\iota} \circ \tilde{\Phi}$ on $\overline{S}'$ for a unique $\tilde{\Phi}$. \qed

\begin{rem} \label{remarks to the lift of f} \rm 
We may use the same techniques to construct a lift of the morphism $f:X' \rightarrow  A^{\rm an}=E/M$ to the
uniformization $E$ of $A$. In general, such a lift does not exist globally on $X'$. By \cite{BL2}, Theorem 1.2, such a lift exists if $H^1(X',\zdop)=0$. Let us consider the formal open subset $U_S:= \Val^{-1}(\Delta_S)$ of $X'$ for the canonical simplex $\Delta_S$ associated to a stratum $S$ of $\tilde{\Xcal}'$. Then $U_S$ is the generic fibre of a formal open subset $\Ucal_S$ of $\Xcal'$. Obviously, $\Ucal_S$ is strictly semistable with skeleton $\Delta_S$. Since the skeleton is a proper deformation retraction of the generic fibre, we get $H^1(U_S,\zdop)=0$ and hence we may apply the above results to get the desired lift $F:U_S \rightarrow E$ of $f|_{U_S}$. 

Note that it is not necessary to appeal to such sophisticated results. We may just use Proposition \ref{f_aff} to conclude that $\overline{f}_{\rm aff}(\Delta_S)$ is contained in a polytope $\Deltabar \in \overline{\Ccal_0}$ and hence $f(U_S)$ is contained in the formal open subset $\valbar^{-1}(\Deltabar)$ of $A^{\rm an}$. The preimage of $\valbar^{-1}(\Deltabar)$ in $E$ is the disjoint union of the formal open subsets $\val^{-1}(\Delta)$, where $\Delta$ ranges over all polytopes of $\Ccal_0$ mapping (bijectively) onto $\Deltabar$ with respect to the residue map $\rdop^n \rightarrow \rdop^n/\Lambda$. Obviously, $M$ acts faithfully and transitively on the set of all such $\val^{-1}(\Delta)$. Since the quotient morphism $p:E \rightarrow A^{\rm an}=E/M$ maps $\val^{-1}(\Delta)$ isomorphically onto $\valbar^{-1}(\Deltabar)$, we get a lift $F:U_S \rightarrow E$ of the restriction of $f$ to $U_S$, unique up to the $M$-action. Note that this construction was partially used in the proof of Proposition \ref{f_aff}.

By Proposition \ref{f_aff}, we get a unique map  $f_{\rm aff}:\Delta_S \rightarrow \rdop^n$ such that $f_{\rm aff} \circ \Val = \val \circ F$ on $U_S$. Moreover, $f_{\rm aff}$ is affine on $\Delta_S$. Conversely, every lift of $\overline{f}_{\rm aff}:\Delta_S \rightarrow \rtor$ to $\rdop^n$ is an affine map $f_{\rm aff}:\Delta_S\rightarrow  \rdop^n$ and there is a unique lift $F:U_S \rightarrow E$ of the restriction of $f$ to $U_S$ such that $f_{\rm aff} \circ \Val = \val \circ F$ on $U_S$. This follows from the fact that the lift of $\overline{f}_{\rm aff}$ is unique up to $\Lambda=\val(M)$-translation.

Finally, we note that we may use such lifts $F$ to construct the lifts $\tilde{\Phi}_0:\overline{S} \rightarrow \tilde{\Ecal}_0$  and $\tilde{\Phi}:\overline{R}\rightarrow \tilde{\Ecal}$ from Lemma \ref{lift of phi}. We will give the construction for $\tilde{\Phi}_0$, but everything works similarly for $\tilde{\Phi}$. The map $f$ is the generic fibre of the formal morphism $\varphi_0:\Xcal' \rightarrow \Acal_0$ and  $\Ucal_S$ is a formal open subset of $\Xcal'$, hence the lift $F$ is the generic fibre of a formal lift $\Phi_S:\Ucal_S \rightarrow \Ecal_0$ of $\varphi_0$. We conclude that the reduction $\tilde{\Phi}_S$ agrees with a lift $\tilde{\Phi}_0$ from Lemma \ref{lift of phi} on the dense stratum $S$ of $\overline{S}$. Similarly, we could argue for every other stratum $T \subset \overline{S}$ to describe the  restriction of $\tilde{\Phi}_0$ to $T$ as the reduction of a formal lift of $\varphi_0$. However, it is not always possible to describe $\tilde{\Phi}_0$ by the reduction of a single formal lift defined on a formal open subset of $\Xcal'$. The problem arises if there are two strata $T_1,T_2$ in $\overline{S}$ such that $f_{\rm aff}(\Delta_{T_1}) \cup f_{\rm aff}(\Delta_{T_2})$ does not map bijectively onto $\overline{f}_{\rm aff}(\Delta_{T_1}) \cup \overline{f}_{\rm aff}(\Delta_{T_2})$. In this case, the lift $F$ will be multivalued on $U_{T_1} \cup U_{T_2}$ and the above covering argument breaks down. This problem can be omitted if we start with a sufficiently fine polytopal decomposition $\overline{\Ccal_0}$ of $\rtor$ and then $\tilde{\Phi}_0$ is indeed the reduction of a single formal lift.

\end{rem}

\begin{art} \label{setup for degree} \rm 
Our goal is to compute the degree of an irreducible component $Y$ of $\tilde{\X}''$ with respect to a line bundle $\Lcal$ on $\Acal$. This can be done in terms of convex geometry under the following  hypotheses fulfilled in our applications:

We still have our abelian variety $A$ over $\kdop$ with uniformization $E$  and the morphism $\varphi_0:\Xcal' \rightarrow \Acal_0$, where $\Acal_0$ is the Mumford model of $A$ associated to the $\Gamma$-rational polytopal decomposition $\overline{\Ccal_0}$ of $\rtor$ and where $\Xcal'$ is a strictly semistable formal scheme over $\kdop^\circ$ with connected generic fibre $X'$. We assume that the generic fibre $f:X' \rightarrow A^{\rm an}$ of $\varphi_0$ is proper and hence the special fibre $\tilde{\varphi_0}$ is also proper (see \cite{Gu2}, Remark 3.14). Let $\overline{\Ccal_1}$ be a $\Gamma$-rational polytopal decomposition of $\rtor$ with associated Mumford model $\Acal_1$ of $A$. We choose now 
$\Ccalbar := \overline{\Ccal_0} \cap \overline{\Ccal_1}:=\{\overline{\Delta_0} \cap  \overline{\Delta_1} \mid \overline{\Delta_0} \in \overline{\Ccal_0} \, , \overline{\Delta_1} \in \overline{\Ccal_1}\}$.
Let $\Acal$ be the Mumford model of $A$ associated to $\Ccalbar$. We apply Propositions \ref{formal refinement for alteration} and \ref{minimal formal analytic structure} to this setup. By \eqref{com diagram}, we get the following commutative diagram of canonical morphisms of admissible formal schemes over $\kcirc$:
\begin{equation}  \label{diagram of can morphisms}
\begin{CD} 
\Xcal'' @>\varphi>> \Acal @>{\iota_1}>> \Acal_1\\
@VV{\iota'}V    @VV{\iota_0}V\\
\Xcal' @>{\varphi_0}>> \Acal_0
\end{CD}
\end{equation}
Recall that all admissible formal schemes in \eqref{diagram of can morphisms} are associated to formal analytic varieties and that the morphism $\varphi$ is determined by the fact that the rectangle is cartesian on the level of formal analytic varieties.

By Corollary \ref{strata corollary}, the irreducible component $Y$ of $\tilde{\Xcal}''$ corresponds to the vertex $\ub'=\xi_Y$ of the $\Gamma$-rational subdivision  
$$\Dcal= \{ \Delta_S \cap \overline{f}_{\rm aff}^{-1}(\sigmabar) \mid S \in \xstr, \, \sigmabar \in \Ccalbar\}$$
of $\xskel$. There is a unique $S \in \xstr$ such that $\ub' \in \relint(\Delta_S)$. If $S'$ is the dense stratum in $Y$, then Corollary \ref{strata corollary} yields $S=\tilde{\iota}'(S')$. We choose a lift $f_{\rm aff}:\Delta_S \rightarrow \rdop^n$ of $\overline{f}_{\rm aff}$. By Lemma \ref{lift of phi}, there is a lift  $\tilde{\Phi}_0:\overline{S} \rightarrow \tilde{\Ecal}_0$ (resp. $\tilde{\Phi}:Y \rightarrow \tilde{\Ecal}$)  of $\tilde{\varphi}_0$ (resp. $\tilde{\varphi}$) to the special fibre of the formal $\kdop^\circ$-model $\Ecal_0$ (resp. $\Ecal$) of $E$ associated to $\Ccal_0$ (resp. $\Ccal$) with $\tilde{\Phi}_0 \circ \tilde{\iota}'=\tilde{\iota} \circ \tilde{\Phi}$. 

Let $L$ be a line bundle on $A$. The role of $\Acal_1$ becomes now clear as we assume that $L$ has a formal $\kdop^\circ$-model $\Lcal$ of $L$ on $\Acal_1$ corresponding to a continuous piecewise affine function $f_\Lcal$ as in Proposition \ref{Mumford line models} (applied to $\overline{\Ccal_1}$). We assume that $g:=f_\Lcal \circ f_{\rm aff}$ is convex in a neighbourhood of $\ub'$. In the light  of Proposition \ref{relatively ampel}, this is a natural positivity assumption for $\Lcal$. We have seen in \ref{skeleton} that we may identify $\Delta_S$ with the  simplex $\{w_0'+ \dots + w_r' = v(\pi)\}$ in $\rdop_+^{r+1}$. In the following, it is more convenient to identify $\Delta_S$ with the  simplex $\{w_1' + \dots + w_r' \leq v(\pi) \}$ in $\rdop_+^{r}$ obtained by omitting the coordinate $w_0'$. Then we define a polytope $\{\ub'\}^g$ in $\rdop^r$ by
$$\{\ub'\}^g:= \{\omega \in \rdop^r \mid \wb' \in \Delta \in {\rm star}_r(\ub') \Rightarrow \omega \cdot (\wb'- \ub') \leq g(\wb')-g(\ub') \},$$
where ${\rm star}_r(\ub')$ is the set of $r$-dimensional polytopes in $\Dcal$ with vertex $\ub'$. The volume of $\{\ub'\}^g$ with respect to the Lebesgue measure on $\rdop^r$ will be denoted by $\vol(\{\ub'\}^g)$. 
By \ref{line bundles}, there is  a line bundle $\Hcal$ on the formal abelian scheme $\Bcal$ from the Raynaud extension \eqref{Raynaud extension 2} such that $p^*(L^{\rm an})=q^*(H)$ for the generic fibre $H$ of $\Hcal$ and the canonical morphisms $p:E \rightarrow A^{\rm an}=E/M$, $q:E \rightarrow B = \Bcal^{\rm an}$.  We have now the following commutative diagram of varieties over $\ktilde$:
\begin{equation}  \label{reduction diagram}
\begin{CD}
Y @>{\tilde{\Phi}}>> \tilde{\Ecal} @>{\tilde{q}}>> \tilde{\Bcal}\\
@VV{\tilde{\iota}'}V    @VV{\tilde{i}_0}V   @VV{\id}V\\
\overline{S}  @>{\tilde{\Phi}_0}>> \tilde{\Ecal}_0  @>{\tilde{q}_0}>> \tilde{\Bcal}
\end{CD}
\end{equation}
For simplicity of notation, we will write $\deg_\Lcal(Y)$ for the degree of $Y$ with respect to the line bundle $(\tilde{\iota}_1 \circ \tilde{\varphi})^*(\tilde{\Lcal})$ and similarly for other degrees. It is always understood that we use the pull-backs of the line bundles with respect to the canonical morphisms from \eqref{diagram of can morphisms} or \eqref{reduction diagram}.

By Proposition \ref{f_aff}, it is easy to deduce that
\begin{equation} \label{decompositon identity}
\Dcal= \{ \Delta_S \cap \overline{f}_{\rm aff}^{-1}(\sigmabar) \mid S \in \xstr, \, \sigmabar \in \overline{\Ccal_1}\}.
\end{equation}
There is a unique $\overline{\Delta_1} \in \overline{\Ccal_1}$ with $\ubb := \overline{f}_{\rm aff}(\ub') \in \relint(\overline{\Delta_1})$. Since the vertex $\ub'$ of $\Dcal$ is contained in $\relint(\Delta_S)$, it follows from \eqref{decompositon identity} that
$$\{\ub'\}= \overline{f}_{\rm aff}^{-1}(\overline{\Delta_1})\cap \Delta_S, \quad
\{\ubb\}= \overline{\Delta_1} \cap \overline{f}_{\rm aff}(\Delta_S).$$
The first equality yields that the affine map $\overline{f}_{\rm aff}|_{\Delta_S}$ is injective and hence $\overline{f}_{\rm aff}(\Delta_S)$ is a $(d-e)$-dimensional simplex in $\rtor$, where $d:=\dim(X')$ and $e:=\dim(S)=d-\dim(\Delta_S)$. We make now the  {\it transversality assumption}
\begin{equation} \label{transversality assumption}
d-e = \codim(\Delta_1, \rdop^n). 
\end{equation}
\end{art}

\begin{prop} \label{degree and dual polytope} 
Using the assumptions from \ref{setup for degree}, we have
$$\deg_\Lcal(Y)= \frac{d!}{e!} \cdot \deg_\Hcal(\overline{S}) \cdot \vol(\{\ub'\}^g).$$
\end{prop}

\proof Let $\Ecal_1$ be the $\kcirc$-model of $E$ associated to $\overline{\Ccal_1}$. In the following, we will always use the canonical morphisms $p_1:\Ecal_1 \rightarrow \Acal_1$, $q_1:\Ecal_1 \rightarrow \Bcal$, $i_1:\Ecal \rightarrow \Ecal_1$ and $\tilde{\Phi}_1:=\tilde{i}_1 \circ \tilde{\Phi}$ to compute degrees. Using $p^*(L^{\rm an})=q^*(H)$, we have the decomposition
\begin{equation} \label{line decomposition}
p_1^*(\Lcal) = q_1^*(\Hcal) \otimes \Ocal_{\Ecal_1}(f_\Lcal)
\end{equation}
for a formal $\kcirc$-model $\Ocal_{\Ecal_1}(f_\Lcal)$ of $O_E$ on $\Ecal_1$. The reason behind the notation is that the formal metric on the trivial bundle $O_E$ associated to the formal model $\Ocal_{\Ecal_1}(f_\Lcal)$ (see Example \ref{formal models}) satisfies
\begin{equation} \label{divisor notation and metric}
-\log \|s\|_{\Ocal_{\Ecal_1}(f_\Lcal)}=f_\Lcal \circ \val,
\end{equation}
where $s$ is the unique meromorphic section of $\Ocal_{\Ecal_1}(f_\Lcal)$ extending the canonical section $1$ of $O_E$. This follows immediately from the definition of $f_\Lcal$ in \eqref{f-char}. In the decomposition \eqref{line decomposition}, $q_1^*(\Hcal)$ reflects the contribution of the abelian part $\Bcal$ to $\Lcal$ and $\Ocal_{\Ecal_1}(f_\Lcal)$ measures the combinatorial contribution from the polytopal decomposition $\Ccal_1$ and from the piecewise affine function $f_\Lcal$. We deduce 
\begin{equation} \label{degree decomposition}
\deg_\Lcal(Y)= \sum_{\ell=0}^d \binom{d}{\ell} d_\ell(Y), 
\end{equation}
from \eqref{line decomposition}, where 
$d_\ell(Y):= \deg_{\underbrace{\Hcal, \dots ,\Hcal}_{\ell}, \underbrace{\Ocal_{\Ecal_1}(f_\Lcal), \dots , \Ocal_{\Ecal_1}(f_\Lcal)}_{d-\ell}}(Y)$.

Our goal is now to prove that $d_\ell(Y)=0$ for $\ell \neq e$ and to compute $d_e(Y)$ using the projection formula with respect to $\tilde{\iota}':Y \rightarrow \overline{S}$ and 
\begin{equation} \label{Hcal commutative}
\tilde{\Phi}^*\left(\tilde{q}^*(\tilde{\Hcal})\right) = \left(\tilde{\iota}'\right)^*\left(\tilde{\Phi}_0^*\left(\tilde{q}_0^*\left(\tilde{\Hcal}\right)\right)\right)
\end{equation}
obtained from \eqref{reduction diagram}.

\vspace{2mm}
\noindent {\it Step 1: The cycle class $c_1(\Ocal_{\Ecal_1}(f_\Lcal))^{d-\ell}. Y$ in $CH(Y)$ is algebraically equivalent to a strata cycle of $Y$. }
\vspace{2mm}

It is always understood that $c_1(\Ocal_{\Ecal_1}(f_\Lcal))$ operates by pull-back with respect to $\tilde{\Phi}_1$ on $Y$. Again, we denote by $s$  the unique meromorphic section of $\Ocal_{\Ecal_1}(f_\Lcal)$ extending the canonical section $1$ of $O_E$. It is enough to show that $\Div(s).\overline{S''}$ is algebraically equivalent to a strata cycle of $Y$ for every $S'' \in \str(\tilde{\Xcal}'')$ with $S'' \subset Y$. We have seen in Proposition \ref{minimal formal analytic structure} that $\tilde{\varphi}=\tilde{\phi}$ maps strata into strata. By Proposition \ref{strata on Mumford models}, we easily deduce the same property for $\tilde{\iota}_1$. Passing to the lift $\tilde{\Phi}_1$, we see that $\tilde{\Phi}_1(S'')$ is contained in a stratum $Z$ of $\tilde{\Ecal}_1$. By Proposition \ref{strata on Mumford models}, $Z$ corresponds to $\relint(\sigma)$ for a unique $\sigma \in \Ccal_1$. Using
$$f_\Lcal(\ub)= \mb_\sigma \cdot \ub + v({a_\sigma})$$
on $\sigma$ with $\mb_\sigma \in \zdop^n$ and ${a_\sigma} \in \kdop^\times$, we deduce from \eqref{divisor notation and metric} that the Cartier divisor $\Div(s)$ is given on $\val^{-1}(\sigma)$ by ${a_\sigma} \cdot \xb^{\mb_\sigma}$. Here, we consider $\chi=\xb^{\mb_\sigma}$ as a meromorphic section of $q_1^*(\Ocal_\chi)$ which restricts to a nowhere vanishing global section on the generic fibre $q^*(O_\chi)$ (see \ref{Raynaud trivializations}). We consider the Cartier divisor $D:= \Div(s/({a_\sigma} \cdot \xb^{\mb_\sigma}))$ on $\Ecal_1$. It  has a well-defined reduction $\tilde{D}$ on a neighbourhood of $\overline{Z}$ which is trivial on $(\val^{-1}(\sigma))\sptilde$ and hence $\tilde{\Phi}_1^*(\tilde{D})$ is a Cartier divisor on $Y$ which is trivial on $S''$. Since $\overline{S''}$ is a strata subset, $\tilde{\Phi}_1^*(\tilde{D}).\overline{S''}$ is a strata cycle  in $Y$. 
We have seen in \ref{Raynaud trivializations} that $\tilde{\Ocal}_\chi$ is algebraically equivalent to $0$ and hence $\Ocal(\tilde{D})|_{\overline Z}$ is algebraically equivalent to $\Ocal_{\Ecal_1}(f_\Lcal)\sptilde|_{\overline Z}$.
By construction, 
$\Div(s).\overline{S''}$ is algebraically equivalent to  $\tilde{\Phi}_1^*(\tilde{D}).\overline{S''}$
 proving the first step.

We need an explicit description of the Cartier divisor $\tilde{\Phi}_1^*(\tilde{D})$ on $Y$ from step 1. 
Let $\Ucal'$ be a formal affine open subset of $\tilde{\Xcal}'$ as in Proposition \ref{normalization of covering} with   \'etale morphism $\psi:\Ucal' \rightarrow \Scal=\Scal_1 \times \Scal_2$ such that $S$ is the distinguished stratum of $\tilde{\Xcal}'$ associated to $\Ucal'$. Passing to a formal open refinement, we may assume that $\varphi_0(\Ucal')$ is contained in a formal trivialization of the Raynaud extension \eqref{Raynaud extension 1} and hence the torus coordinates $x_1, \dots, x_n$ make sense on $\varphi_0(\Ucal')$.

We denote by $\Scal_1'$ the $\kcirc$-model of  the polytopal domain $\Scal_1^{\rm an}$ in $({\mathbb G}_m^r)^{\rm an}$ associated to the refinement $\Dcal \cap \Delta_S$ and let ${\Ucal''}:=(\iota')^{-1}(\Ucal')$. We have seen in Remark \ref{local description by base change} that $\psi':{\Ucal''} \rightarrow \Scal'=\Scal_1' \times \Scal_2$ is the base change of $\psi$ to $\Scal'$ and hence $\psi'$ is \'etale. Let $\psi_1:\Ucal' \rightarrow \Scal_1'$ be the composition of the first projection with $\psi$ and let $\psi_1'$ be  
the base change of $\psi_1$ to $\Scal_1'$. Then $\tilde{\psi}_1'$ is a smooth morphism such that the preimage of a stratum of $\tilde{\Scal}_1'$ is a stratum of $\tilde{\Ucal}''$ (see Remark \ref{preimage of strata wrt psi}). We conclude that $\tilde{\psi}_1'(Y \cap \tilde{\Ucal}'')$ is dense in an irreducible component of $\tilde{\Scal}_1'$ which we denote by $Y_{\ub'}$. This notation is justified by the fact that the irreducible components of $\Scal_1'$ are in bijective correspondence with the vertices of $\Dcal \cap \Delta_S$ (see \cite{Gu1}, Proposition 4.7).

\vspace{2mm}
\noindent {\it Step 2: There is a Cartier divisor $\tilde{D}_1$ on $Y_{\ub'}$ with $(\tilde{\psi}_1')^*(\tilde{D}_1)= \tilde{\Phi}_1^*(\tilde{D})|_{Y \cap \tilde{\Ucal}''}$.}
\vspace{2mm}

By Remark \ref{remarks to the lift of f}, the lift $\tilde{\Phi}_1:Y \cap \tilde{\Ucal}'' \rightarrow \tilde{\Ecal}_1$ is equal to the reduction of a suitable lift $F:U' \rightarrow E$ of $f$. Moreover, $F$ induces a lift $\Delta_S \rightarrow \rdop^n$ of $\overline{f}_{\rm aff}$. We may assume that the lift is equal to the $f_{\rm aff}$ from \ref{setup for degree}. Indeed, $f_{\rm aff}$ is determined up to $\Lambda$-translation and hence the polytope $\{\ub'\}^g$ is also determined up to translation which does not affect the volume in Proposition \ref{degree and dual polytope}. 

We consider the polytopes $\nu$ of $\Ccal_1$ with closed face $\Delta_1$ from \ref{setup for degree}. We have seen in the first step that the Cartier divisor $D$ is given on $\val^{-1}(\nu)$ by $({a_\nu}/{a_\sigma})\cdot \xb^{\mb_\nu-\mb_\sigma}$. It follows easily from the definitions that the polytopes $\mu:=(f_{\rm aff})^{-1}(\nu) \cap \Delta_S$ are just the polytopes of $\Dcal \cap \Delta_S$ with vertex $\ub'$. By Proposition \ref{orbits and refinement} and Corollary \ref{strata corollary}, the reductions of the formal open subsets $\Val^{-1}(\mu)$ cover $Y$. By construction, $F$ induces a formal morphism $\Phi_1:\Ucal'' \rightarrow \Acal_1$ with reduction $\tilde{\Phi}_1$ and hence $\tilde{\Phi}_1^*(\tilde{D})= (\Phi_1^*(D))\sptilde$ on $Y \cap \tilde{\Ucal}''$. By \cite{Gu4}, Proposition 2.11, there is an $n \times r$ matrix $M$  with entries in $\zdop$, $\gamma_i \in \Ocal(\Ucal')^\times$ and $\lambda \in (\kdop^\times)^n$ with
$$F^*(x_i)= \lambda_i \cdot \gamma_i \cdot (\psi^{\rm an})^*(x_1')^{M_{i1}} \dots (\psi^{\rm an})^*(x_r')^{M_{ir}}$$
on the generic fibre $U'$ of $\Ucal'$ for $i=1, \dots, n$. We have seen in \eqref{definition of f_aff} that
$$f_{\rm aff}({\mathbf w'})=M {\mathbf w'} + \lambda$$
for ${\mathbf w'}=(w_1', \dots, w_r') \in \Delta_S$. Note that we use here the identification of $\Delta_S$ with the  simplex $\Sigma_S:= \{w_1' + \dots + w_r' \leq v(\pi)\}$ in $\rdop^r_+$ which is different from the one used in \eqref{definition of f_aff}. Let $\mathbf y: = (x_1', \dots, x_r')$.  We conclude that $\Phi_1^*(D)$ is given on $\Val^{-1}(\mu)\cap U'$ by 
$$\frac{{a_\nu}}{{a_\sigma}} \cdot F^*(\xb^{\mb_\nu -\mb_\sigma})
= \frac{{a_\nu}}{{a_\sigma}} \cdot \lambda^{\mb_\nu -\mb_\sigma} \cdot \gamma \cdot (\psi^{\rm an})^*(\yb)^{(\mb_\nu-\mb_\sigma)^t \cdot M}$$
for some $\gamma \in \Ocal(\Ucal')^\times$. For a Cartier divisor, such a unit $\gamma$  can be omitted. Let $\alpha_\mu \in \kdop^\times$ with $v(\alpha_\mu) + ({(\mb_\nu-\mb_\sigma)^t \cdot M})\cdot \ub'=0$. Then $\Phi_1^*(D)$ is  given by $\alpha_\mu \cdot (\psi^{\rm an})^*(\yb)^{(\mb_\nu-\mb_\sigma)^t \cdot M}$ on $\Val^{-1}(\mu)\cap U'$. These functions are also defined on the formal open subsets $\val^{-1}(\mu)$ of $(\Scal_1')^{\rm f-an}$. Let $\Ucal_\mu$ be the formal affine open subset of $\Scal_1'$ associated to $\val^{-1}(\mu)$ and let $D_1:=\{\Ucal_\mu, \alpha_\mu \cdot \yb^{(\mb_\nu-\mb_\sigma)^t \cdot M}\}$ with $\mu$ ranging over the polytopes of $\Dcal \cap \Delta_S$ with vertex $\ub'$. It is easy to see that $D_1$ is a Cartier divisor on the open subset $\bigcup_\mu \Ucal_\mu$ of $\Scal_1'$ containing $Y_\ub$. We conclude that $(\tilde{\psi}_1')^*(\tilde{D}_1)= \tilde{\Phi}_1^*(\tilde{D})|_{Y \cap \tilde{\Ucal}''}$ proving the second step.

We note that the Cartier divisor ${D}_1$ depends on the choice of the stratum $S''$, but the linear  equivalence class of the Cartier divisor $D_1$ is independent of $S''$ and hence the same is true for the linear equivalence class of 
$\tilde{D}_1$ on $Y_{\ub'}$.

\vspace{2mm}
\noindent {\it Step 3: $d_\ell(Y)=0$ for $\ell \neq e$.}
\vspace{2mm}

If $\ell > e$, then $c_1(\Ocal_{\Ecal_1}(f_\Lcal))^{d-\ell}.Y$ has dimension $\ell > e = \dim(S)$ and hence the projection formula with respect to $\tilde{\iota}':Y \rightarrow \overline{S}$ and \eqref{Hcal commutative} prove 
\begin{equation*}
\begin{split}
d_\ell(Y) &= \deg \left( c_1(\Hcal)^{\ell} . c_1\left( \Ocal_{\Ecal_1}(f_\Lcal) \right)^{d-\ell}.Y \right) \\
&= \deg \left( c_1(\Hcal)^\ell . \tilde{\iota}_*' \left( c_1\left( \Ocal_{\Ecal_1}(f_\Lcal) \right)^{d-\ell}.Y \right) \right) =0.
\end{split}
\end{equation*}
It remains to consider $\ell<e$. We will use the first step for the dense stratum $S'$ in $Y$ (instead of $S''$). We conclude that $\tilde{\Phi}_1(S')$ is contained in the stratum $Z$ of $\tilde{\Ecal}_1$ corresponding to $\relint(\Delta_1)$. By Proposition \ref{strata on Mumford models}, we have $\dim(Z)=\dim(A)-\dim(\Delta_1)$. By the construction in the first step, 
 the cycle class
$$\alpha := (\tilde{\Phi}_1)_*\left(c_1\left(\Ocal_{\Ecal_1}(f_\Lcal)\right)^{d-\ell}.Y\right)=
 c_1\left(\Ocal_{\Ecal_1}(f_\Lcal)\right)^{d-\ell}. (\tilde{\Phi}_1)_*(Y)$$
is algebraically equivalent to a cycle supported in  $\overline{Z_1}$ for a strata subset $Z_1$ of codimension $\geq d-\ell$ in $\overline{Z}$.
We have
$$e=d-\codim(\Delta_1,\rdop^n)=d+\dim(B)-\dim(Z)$$
by our transversality assumption. Using $\ell < e$, we get
$$\dim(Z_1)\leq \dim(Z)-(d-\ell) = \dim(B) + \ell -e <\dim(B).$$
By Proposition \ref{strata on Mumford models}, all strata of $\tilde{\Ecal}_1$ have dimension $\geq \dim(B)$ and hence $\alpha$ is algebraically equivalent to $0$. The projection formula  shows now that \begin{equation*}
d_\ell(Y)= \deg \left( c_1(\Hcal)^\ell . (\tilde{\Phi}_1)_* \left( c_1\left( \Ocal_{\Ecal_1}(f_\Lcal) \right)^{d-\ell}.Y \right) \right) =\deg \left( c_1(\Hcal)^\ell .\alpha \right)= 0.
\end{equation*}

\vspace{2mm}
\noindent {\it Step 4: $d_e(Y)= (d-e)! \cdot \deg_\Hcal(\overline{S}) \cdot \vol(\{\ub'\}^g).$}
\vspace{2mm}

Recall that $r=d-e$. By the first step, we know that $c_1\left(\Ocal_{\Ecal_1}\left(f_\Lcal\right)\right)^{r}.Y$ is algebraically equivalent to an $e$-dimensional strata cycle $W$ of $\tilde{\Xcal}''$. Since $W$ has support in $Y$, its components  have the form $\overline{S_i}$, where $S_i \in \str(\tilde{\Xcal}'')$ corresponding to an open face of $\Dcal_S \cap \Delta_S$ with vertex $\ub'$. This follows from Corollary \ref{strata corollary} as well as the fact that $\tilde{\iota}'$ maps $S_i$ isomorphically onto $S$.

Now we will use the formal open subset ${\Ucal''}:= (\iota')^{-1}(\Ucal')$ of $\Xcal''$ from step 2. Since $\tilde{\iota}'(S_i)=S$, we have $S_i \cap {\Ucal''} \neq \emptyset$. The same holds for every stratum relevant in the intersection process for $W$ described in step 1. We conclude that we may compute $W$ on ${\Ucal''}$. We note that $\tilde{\psi}_1':\tilde{\Ucal}'' \rightarrow \tilde{\Scal}_1'$ is a smooth morphism and that the stratification of $\tilde{\Ucal
}''$ is obtained by the preimages of the stratas of $\tilde{\Scal}_1'$ (see Remark \ref{preimage of strata wrt psi}). By the second step and the compatibility of flat pull-back with the intersection operations (see \cite{Fu1}, Proposition 2.3), the intersection process on $Y \cap {\Ucal''}$ leading to $W \cap {\Ucal''}$ may be first performed on $Y_{\ub'}$ giving a cycle $W'$ and then $W \cap {\Ucal''} = (\tilde{\psi}_1')^*(W')$. To obtain $W'$, we have just to replace $Y$ by $Y_{\ub'}$ and the Cartier divisors $\tilde{\Phi}_1^*(\tilde{D})$ by $\tilde{D}_1$. It is clear that $W' = \sum n_i{S_i'}$, where the $0$-dimensional strata $S_i'$ of $\tilde{\Scal}_1'$ correspond to the same open faces as the $S_i$. We deduce $W=\sum n_i \overline{S_i}$ and $\sum n_i = \deg_{D_1}(Y_{\ub'})$. Using projection formula with respect to $\tilde{\iota}':Y \rightarrow \overline{S}$ and  \eqref{Hcal commutative}, we get
$$d_e(Y)=\deg \left(c_1(\Hcal)^e . W \right) = \deg \left(c_1(\Hcal)^e . \tilde{\iota}_*'(W) \right).$$
We have noticed that  $S_i \cong S$ and hence we get
\begin{equation} \label{degree formula}
d_e(Y) = \deg_\Hcal(\overline{S}) \sum_i n_i [\overline{S_i}:\overline{S}] = \deg_\Hcal(\overline{S}) \deg_{D_1}(Y_{\ub'}).
\end{equation}
To compute the degree of $Y_{\ub'}$, we will use the theory of toric varieties. The projection ${\mathbb G}_m^{r+1} \rightarrow {\mathbb G}_m^{r}$, given by $(x_0', \dots, x_r') \mapsto (x_1', \dots, x_r')$, leads to an isomorphism of $\Scal_1^{\rm f-an}$ with the polytopal domain $\val^{-1}(\Sigma_S)$ for the simplex $\Sigma_S := \{w_1' + \dots + w_r' \leq v(\pi)\}$ in $\rdop_+^r$. 
We recall from \ref{setup for degree} that we identify $\Sigma_S$ with $\Delta_S$ and hence  $Y_{\ub'}$ is equal to the $({\mathbb G}_m^{r})_{\ktilde}$-toric variety associated to the vertex $\ub'$ of $\Dcal \cap \Sigma_S$ (see \cite{Gu4}, Proposition 4.7). As in the second step, let $\nu \in \Ccal_1$ with closed face $\Delta_1$ and let $\mu := f_{\rm aff}^{-1}(\nu) \cap \Delta_S$. Then the polytopes $\mu$ are just the polytopes of $\Dcal \cap \Sigma_S$ with vertex $\ub'$. We have seen in the second step that the Cartier divisor $D_1$ is given on the formal open subset 
$\val^{-1}(\mu)$ of $(\Scal_1')^{\rm an}$ by $\alpha_\mu \cdot\yb^{(\mb_\nu-\mb_\sigma)^t \cdot M}$. In the theory of toric varieties, the Cartier divisor $\tilde{D}_1|_{Y_{\ub'}}$ induces a polyhedron $P$ as the set of all $\omega \in \rdop^r$ with 
$$ \forall \wb' \in \mu \in {\rm star}_r(\ub') \Rightarrow \omega \cdot (\wb'-\ub') \leq (\mb_\nu - \mb_\sigma)^t  \cdot M \cdot (\wb' - \ub') .$$
It is easy to see that $P$ is a translate of our polytope $\{\ub'\}^g$. By \cite{Fu2}, \S 5.3, Corollary on p. 111, we get
\begin{equation*} \label{final conclusion for degree}
\deg_{D_1}(Y_{\ub'}) = r! \cdot \vol(P) = r! \cdot \vol(\{\ub'\}^g).
\end{equation*}
Together with \eqref{degree formula}, this proves the fourth step. Finally, the proposition is a consequence of \eqref{degree decomposition}, step 3 and step 4.\qed 

\begin{rem} \label{strictly pluristable} \rm
Berkovich has defined the skeleton $S(\Xcal')$ more generally for a non-degenerate pluristable formal scheme $\Xcal'$ over $\kcirc$ and he has shown that $S(\Xcal')$ has a canonical piecewise linear structure (see \cite{Ber5}). If $\Xcal'$ is strongly non-degenerate, then there is a well defined proper strong deformation retraction from the generic fibre $X'$ to $S(\Xcal')$ which generalizes the map $\Val$.

All the results of Section 5 can be generalized to a strongly non-degenerate {\it strictly} pluristable $\Xcal'$. This is based on the following facts proved in  the appendix: The linear pieces of $S(\Xcal')$ are given by canonical plurisimplices $\Delta_S$ corresponding to the strata $S$ of $\tilde{\Xcal}'$. Moreover, $\Delta_S$ is a polytope with associated polytopal domain $U_{\Delta_S}$ (see \ref{polytopal domains}). In analogy to Proposition \ref{normalization of covering}, $\Xcal'$ consists locally of open building blocks $\Ucal'$ such that $S(\Ucal')$ is a canonical plurisimplex $\Delta_S$ of $S(\Xcal')$ and there is an \'etale morphism $\psi:\Ucal'\rightarrow U_{\Delta_S}^{\rm f-sch}$. 

Similarly as in the strictly semistable case, this allows us to prove the results of this section by using well-known results for polytopal domains. Moreover, we could replace strictly semistable formal scheme in Section 6 by strongly non-degenerate strictly pluristable formal schemes. This is straight forward and we leave the details to the reader.
\end{rem}

\section{Canonical measures}

In this section, $K$ is a field with a discrete valuation $v$. We denote by $\kdop$ the completion of the algebraic closure of $K$. Note that $\kdop$ is algebraically closed (\cite{BGR}, Proposition 3.4.1/3) and the value group $\Gamma$ is equal to $\qdop$. 

We consider a geometrically integral $d$-dimensional closed subvariety $X$ of $A$ over $K$. In \S 3, we have defined canonical measures on $X$. Now we will compute them explicitly in terms of convex geometry. The main idea is to choose a Mumford model of $A$ and a semistable alteration of $X$ to apply the results from \S 4 and \S 5. Note that the restriction to geometrically integral varieties is no serious restriction. In general, we may perform a finite base change and then we can proceed by linearity in the components. 

\begin{art} \rm \label{infinite extension}

For our computations, Proposition \ref{degree and dual polytope} will be crucial. To fulfill its trans\-versality assumption \eqref{transversality assumption}, we have to choose the polytopal decomposition of the Mumford model ``completely irrational''. We choose 
an infinite dimensional $\qdop$-subspace $\Gamma'$ of $\rdop$ containing $\qdop$. By \cite{Bou}, Ch. VI, $n^\circ$ 10, Prop. 1, there is an algebraically closed field $\kdop'$, complete with respect to a valuation $v'$  extending $v$ such that the value group $v'((\kdop')^\times)$ is $\Gamma'$.  
\end{art}

\begin{art} \rm \label{Semistable alteration}
We denote the analytic space over $\kdop$ associated to $X$ by $\Xan$. Let $\overline{\Ccal_0}$ be a rational polytopal decomposition of $\rtor$ with associated Mumford model $\Acal_0$ of $A^{\rm an}$ over $\kcirc$. We denote the closure of $\Xan$ in $\Acal_0$ by $\Xcal_0$ which is a formal $\kcirc$-model of $\Xan$ (see \cite{Gu2}, Proposition 3.3). By de Jong's alteration theorem (\cite{dJ}, Theorem 6.5) applied to a projective $\kcirc$-model of $\Xan$ dominating $\Xcal_0$ (see \cite{Gu3}, Proposition 10.5), there is always a {\it semistable alteration} $\varphi_0:\Xcal' \rightarrow \Xcal_0$ which means that the generic fibre $f:X' \rightarrow \Xan$ of $\varphi$ is a proper surjective morphism and $X'$ is an irreducible $d$-dimensional analytic space which is the generic fibre of a strictly semistable admissible formal scheme $\Xcal'$ over $\kcirc$. It follows from \cite{Gu2}, Remark 3.14, that $\tilde{\varphi}$ is a proper surjective morphism between the special fibres. 
\end{art}

\begin{art} \rm \label{non-degenerate canonical simplex}
We will use the notations from the previous sections. Let $E$ be the uniformization of $A$, i.e. $A^{\rm an}=E/M$ for a discrete subgroup $M$ in $E$ with complete lattice $\Lambda = \val(M)$ in $\rdop^n$. Let $\Ecal_0$ be the $\kcirc$-model of $E$ associated to the polytopal decomposition $\Ccal_0$ of $\rdop^n$ (see \ref{Mumford's construction}). 

Let $S \in \str(\tilde{\Xcal}')$ with canonical simplex $\Delta_S$ in the skeleton $S(\Xcal')$. By Lemma \ref{lift of phi}, there is a lift $\tilde{\Phi}_0:\overline{S} \rightarrow \tilde{\Ecal}_0$ of $\tilde{\varphi}_0: \overline{S} \rightarrow \tilde{\Acal}_0$, unique up to the $M$-action on $\tilde{\Ecal}_0$. If $q_0:\Ecal_0 \rightarrow \Bcal$ denotes the unique morphism extending $q:E \rightarrow B=\Bcal^{\rm an}$ from the Raynaud extension \eqref{Raynaud extension 2}, then $\tilde{q}_0 \circ \tilde{\Phi}_0$ is unique up to $q(M)$-translation on the abelian variety $\tilde{\Bcal}$ over $\ktilde$.

A canonical simplex $\Delta_S$ is called {\it non-degenerate} with respect to the morphism $f$ if $\dim(\overline{f}_{\rm aff}(\Delta_S))=\dim(\Delta_S)$ and $\dim(\tilde{q}_0 \circ \tilde{\Phi}_0(S))= \dim(S)$. This definition does not depend on the choice of the lift $\tilde{\Phi}_0$. Moreover, it  depends only on $\Xcal'$ and $f$, but not on the choice of $\Ccal_0$. This means that if we have a second rational polytopal decomposition $\overline{\Ccal_0'}$ of $\rtor$ with associated Mumford model $\Acal_0'$ and with a semistable alteration $\varphi_0':\Xcal' \rightarrow \Acal_0'$ such that the generic fibre is again $f$, then the definitions of non-degenerate canonical simplices agree. Indeed, the independence of the first condition is obvious and the invariance of the second condition follows from an easy diagram chase involving Lemma \ref{lift of phi} by passing to the common refinement $\Ccal_0 \cap \Ccal_0'$. 
\end{art}

\begin{art} \rm \label{sigma-generic and transversal}
Let $\Sigma$ be a $\Lambda$-periodic set of polytopes such that $\Sigmabar:=\{\sigmabar \subset \rtor \mid \sigma \in \Sigma \}$ is a finite set. If $\sigma$ is a polytope in $\Sigma$, then we assume that all closed faces of $\sigma$ are also in $\Sigma$. Let $\adop_\sigma$ be the affine space in $\rdop^n$ generated by the polytope $\sigma$.

The polytopal decomposition $\Ccalbar$ of $\rtor$ is said to be {\it $\overline{\Sigma}$-generic}  if the following conditions hold for every $\sigma \in \Sigma$,  $\Delta \in \Ccal$:
\begin{itemize}
\item[(a)] $\dim(\adop_\sigma \cap \adop_\Delta) = D$ if $D:=\dim(\sigma)+\dim(\Delta)-n \geq 0$,
\item[(b)] $\adop_\sigma \cap \adop_\Delta = \emptyset$ if $D <0$.
\end{itemize}
By \cite{Gu4}, Proposition 8.2, every $\overline{\Sigma}$-generic $\Ccalbar$ is {\it $\Sigmabar$-transversal} which means that $\Delta \cap \sigma$ is either empty or of dimension $\dim(\Delta) + \dim(\sigma) - n$ for all $\Delta \in \Ccal$, $\sigma \in \Sigma$. 
\end{art}

\begin{lem} \label{totally irrational construction}
Let $L$ be an ample line bundle on $A$. Then there is a $\Gamma'$-rational polytopal decomposition $\overline{\Ccal_1}$ of $\rtor$ with the following properties:
\begin{itemize}
\item[(a)] $\overline{\frac{1}{m} \Ccal_1}$ is $\overline{\Sigma}$-generic and hence $\Sigmabar$-transversal for all $m \in \ndop \setminus \{0\}$.
\item[(b)] If $\Acal_1$ denotes the formal $(\kdop')^\circ$-model of $A_{\kdop'}$ associated to $\overline{\Ccal_1}$, then there are $N \in \ndop \setminus \{0\}$ and a formal $(\kdop')^\circ$-model $\Lcal$ of $L^{\otimes N}$ on $\Acal_1$ corresponding to a function $f_\Lcal$ as in Proposition \ref{Mumford line models} which is a strongly polyhedral convex function with respect to $\Ccal_1$. 
\end{itemize}
\end{lem}

\proof In \cite{Gu4}, Lemma 8.4, this was proved for a totally degenerate abelian variety $A$. Using Pro\-po\-si\-tion \ref{Mumford line models}, the same proof applies here. \qed

\begin{art} \rm \label{setup for explicit}
We keep the assumptions from \ref{Semistable alteration} and we consider an ample line bundle $\overline{L}$ on $A$ endowed with a canonical metric.
 
Let $\Delta_S$ be a canonical simplex of the skeleton $S(\Xcal')$ which is non-degenerate with respect to $f$. By \ref{skeleton}, we may identify $\Delta_S$ with the simplex $\{u_0'+ \dots + u_r' = v(\pi)\}$ in $\rdop_+^{r+1}$. In the following, it is more convenient to identify $\Delta_S$ with the  simplex
$\Sigma_S := \{\ub' \in \rdop_+^r \mid u_1' + \dots + u_r' \leq v(\pi)\}$
by omitting the  coordinate $u_0'$. Let us choose an affine lift $f_{\rm aff}: \Delta_S \rightarrow \rdop^n$ of the map $\overline{f}_{\rm aff}$ from Proposition \ref{f_aff}. Using the identification $\Delta_S=\Sigma_S$, there is a unique injective linear map $\ell_S^{(0)}:\rdop^r \rightarrow \rdop^n$ extending $f_{\rm aff} - f_{\rm aff}(\mathbf 0)$. By \eqref{definition of f_aff}, $\ell_S^{(0)}$ is defined over $\zdop$ and hence $\Lambda_S:= (\ell_S^{(0)})^{-1}(\Lambda)$ is a complete rational lattice in $\rdop^r$. The positive definite bilinear form $b$ associated to $L$ (see \ref{line bundles}) induces a complete lattice
$$\Lambda_S^L:=\{b(\ell_S^{(0)}(\cdot),\lambda) \mid \lambda \in \Lambda\}$$
on $(\rdop^r)^*=\rdop^r$. We denote by $\vol$ the  volume with respect to the  Lebesgue measure on $\rdop^r$.

There is an ample line bundle $\Hcal$ on the abelian scheme $\Bcal$ from the Raynaud extension \eqref{Raynaud extension 2} of $A$ with generic fibre $H$ such that $p^*(L)=q^*(H)$ on $E$ (see \ref{line bundles}). As in \ref{setup for degree}, we define the degree $\deg_\Hcal(\overline{S})$ of $S \in \str(\tilde{\Xcal}')$ by using the  lift  $\tilde{\Phi}_0:\overline{S} \rightarrow \tilde{\Ecal}_0$ and $\tilde{q}_0:\tilde{\Ecal}_0 \rightarrow \tilde{\Bcal}$.  
\end{art}

\begin{thm} \label{explicit formula for canonical measure}
Under the hypothesis in \ref{setup for explicit}, the support of the positive measure $\mu:= c_1(f^*(\overline{L}))^{\wedge d}$ is equal to the union of the canonical simplices of $S(\Xcal')$ which are non-degenerate with respect to $f$. For a measurable subset $\Omega$ contained in the relative interior of such a simplex $\Delta_S$ and $r:=\dim(\Delta_S)$, we have
\begin{equation} \label{explicit formula}
\mu(\Omega)=\frac{d!}{(d-r)!} \cdot \deg_\Hcal(\overline{S}) \cdot \frac{\vol({\Lambda_S^L})}{\vol(\Lambda_S)} \cdot \vol(\Omega).
\end{equation}
\end{thm}

\begin{rem} \rm \label{generalization to several line bundles}
The theorem generalizes easily to several canonically metrized ample line bundles $\overline{L_1}, \dots , \overline{L_d}$ on $A$. Let $\mu := c_1(f^*(\overline{L_1})) \wedge \dots \wedge c_1(f^*(\overline{L_d}))$ and let $\Hcal_j$ be an ample line bundle on $\Bcal$ with $p^*(L_j)=q^*(H_j)$. Then the support of $\mu$ will be again equal to the union of all canonical simplices $\Delta_S$ which are non-degenerate with respect to $f$. 

We are going to  describe the canonical measure $\mu(\Omega)$ for any measurable subset $\Omega$ of a canonical simplex $\Delta_S$. For $r:=\dim(\Delta_S)$, let $\vol(\Lambda_S^{L_1}, \dots , \Lambda_S^{L_r})$ be the mixed volume in $\rdop^r$ of the corresponding fundamental lattices. This is a positive number which agrees with $\vol(\Lambda')$ if all lattices $\Lambda_S^{L_i}$ are equal to a single lattice $\Lambda'$. Moreover, the mixed volume is symmetric and multilinear with respect to Minkowski sum of fundamental lattices. We conclude that $\vol(\Lambda_S^{L_1}, \dots , \Lambda_S^{L_r})$ is multilinear and symmetric with respect to the line bundles $L_1, \dots, L_r$. For more details, we refer to \cite{Gu4}, A6. 

The generalization of Theorem \ref{explicit formula for canonical measure}  can be stated now as
\begin{equation} \label{explicit formula for several line bundles}
\mu(\Omega) = r! \sum_{\mathbf i} \deg_{\Hcal_{j_1}, \dots, \Hcal_{j_s}}(\overline{S}) \cdot 
\frac{\vol(\Lambda_S^{L_{i_1}}, \dots , \Lambda_S^{L_{i_r}})}{\vol(\Lambda_S)} \cdot \vol(\Omega),
\end{equation}
where $\mathbf i$ ranges over $\{1, \dots, d\}^r$ with $i_1<i_2< \dots <i_r$ and where $j_1 < \dots <j_s$ is the complement of $\mathbf i$ in $\{1, \dots , d\}$. Both sides of \eqref{explicit formula for several line bundles} are multilinear and symmetric with respect to $\overline{L_1}, \dots , \overline{L_d}$. Since symmetric real-valued multilinear forms are determined by the restriction to the diagonal, \eqref{explicit formula for several line bundles} follows from \eqref{explicit formula}. 
\end{rem}

\begin{cor} \label{explicit for non-ample}
If $\overline{L_1}, \dots , \overline{L_d}$ are arbitrary line bundles on $A$ endowed with ca\-non\-ical metrics, then $\mu:=c_1(f^*(\overline{L_1})) \wedge \dots \wedge c_1(f^*(\overline{L_d}))$ is  supported in the union of canonical simplices of $S(\Xcal')$ which are non-degenerate with respect to $f$ and the restriction of $\mu$ to such a simplex is a multiple of the Lebesgue measure.
\end{cor}

\proof This follows from \eqref{explicit formula for several line bundles} and multilinearity. \qed

\begin{rem} \label{extreme cases} \rm
Theorem \ref{explicit formula for canonical measure} is well known in the two extreme cases of abelian varieties. If $A$ is an abelian variety of potentally good reduction, then \eqref{explicit formula} shows that $\mu = \sum_Y \deg_\Hcal(Y) \delta_{\xi_Y}$ with $Y$ ranging over the irreducible components of $\tilde{\Xcal}'$. This is a special case of Proposition \ref{Chambert-Loir's measures}(d) and was first proved by Chambert-Loir in \cite{Ch}.

If $A$ is an abelian variety which is totally degenerate at $v$, then Theorem \ref{explicit formula for canonical measure} shows that the support of $\mu$ is equal to the union of all $d$-dimensional canonical simplices $\Delta_S$ of $S(\Xcal')$ with $\dim(\overline{f}_{\rm aff}(\Delta_S))=d$ and we have $\mu(\Omega)=d! \cdot \vol(\Lambda_S^L)\cdot \vol(\Omega)/\vol(\Lambda_S)$ for any measurable subset $\Omega$ of $\Delta_S$. This was proved in \cite{Gu4}, Theorem 9.6.
\end{rem}

\noindent{\bf Proof of Theorem \ref{explicit formula for canonical measure}: \/}  \rm
By \cite{Gu4}, Remark 3.14, the measure $\mu$ is independent of the odd part $L_-$ of $L$. Moreover, $L_-$ does not influence the bilinear form $b$ of $L$ and hence we may assume that $L$ is a symmetric ample line bundle. It will be crucial for the proof to choose a Mumford model $\Acal$ of $A$ as ``generic'' as possible. Let $\Sigmabar$ be the set of simplices $\{\overline{f}_{\rm aff}(\Delta_S) \mid S \in \str(\tilde{\Xcal}') \}$ together with all their closed faces. Then we will use the $\Gamma'$-rational polytopal decomposition $\overline{\Ccal_1}$ of $\rtor$ from Lemma \ref{totally irrational construction} with associated Mumford model $\Acal_1$. By multilinearity, we may assume that the strongly polyhedral convex function $f_\Lcal$ from Lemma \ref{totally irrational construction} induces a model $\Lcal$ of $L$ on $\Acal_1$. For $m \geq 1$,  let $\Acal_m'$ be the Mumford model of $A$  associated to the $\Gamma'$-rational polytopal decomposition 
$\overline{\Ccal_m'}:=\overline{\Ccal_0} \cap \overline{ \frac{1}{m} \Ccal_1}$ (see \ref{setup for degree}). 
Note that $\Acal_m'$, $\Lcal$ and $\Acal$ are only defined over the valuation ring of the ``large'' field extension $\kdop'$. Since $\mu$ is invariant under base change (\cite{Gu4}, Remark 3.10), we may perform analytic calculations for $\mu$ over $\kdop'$.

We fix a rigidification on $L$ such that the given canonical metric $\canmetr$ on $L$ is given by \eqref{Tate's limit} in Example \ref{canonical metrics on abelian var}. 
Let $\Xcal_m$ be the closure of $X$ in $\Acal_m'$. If we apply Propositions \ref{formal refinement for alteration} and \ref{minimal formal analytic structure} to the polytopal decomposition $\Ccal_m'$ instead of $\Ccal$, then we get a minimal formal analytic structure $\X_m''$ on $X'$ which refines $(\Xcal')^{\rm f-an}$ such that our given morphism $f:X' \rightarrow A^{\rm an}$ extends to a morphism $\phi_m:\X_m'' \rightarrow (\Acal_m')^{\rm f-an}$. Let $\varphi_m:\Xcal_m'' \rightarrow \Acal_m'$ be the associated morphism of admissible formal schemes over $\kcirc$.

\vspace{2mm}
\noindent {\it Step 1:  $\mu$ is the weak limit of discrete measures $\mu_m$ on $X'$ which are supported in the preimages of the generic points of the irreducible components of $\tilde{\Xcal}_m''$ with respect to the reduction map.}
\vspace{2mm} 

\noindent This will be a consequence of Tate's limit argument (see \eqref{Tate's limit} in Example \ref{canonical metrics on abelian var}). We may assume that the metric $\metr$ in \eqref{Tate's limit} is equal to the formal metric $\metr_\Lcal$. We note that multiplication by $m$ extends uniquely to a morphism $\psi_m: \Acal_m' \rightarrow \Acal_1$ (see \cite{Gu4}, Proposition 6.4). Then $\Lcal_m :=  \psi_m^*(\Lcal)$ is a $(\kdop')^\circ$-model of $[m]^*(L)$ on $\Acal_m'$ with associated formal metric $[m]^*\metr$ and hence we have
$f^*([m]^*\metr)= \metr_{\Lcal_m''}$
for the formal $(\kdop')^\circ$-model $\Lcal_m'':= \varphi_m^*(\Lcal_m)$ of $f^*([m]^*(L))=f^*(L)^{\otimes m^2}$. By \eqref{Tate's limit}, we get
$$f^*\canmetr= \lim_{m \to \infty} \metr_{\Lcal_m''}^{1/m^2}.$$
If we use this uniform limit together with Proposition \ref{Chambert-Loir's measures}, then we get
\begin{equation} \label{limit of discrete measures}
\mu = \lim_{m \to \infty} m^{-2d} \sum_Z \deg_{\tilde{\Lcal}_m''}(Z) \, \delta_{\xi_Z},
\end{equation}
where $Z$ ranges over all irreducible components of $\tilde{\Xcal}_m''$.

\vspace{2mm}
\noindent {\it Step 2:  A first determination of $\supp(\mu)$.}
\vspace{2mm}

\noindent By Corollary \ref{strata corollary}(g) and Proposition \ref{formal refinement for alteration}, the points $\xi_Z$  are the vertices of the subdivision
$\Dcal_m:= \{\Delta_S \cap \overline{f}_{\rm aff}^{-1}(\sigmabar) \mid S \in \str(\tilde{\Xcal}'),\, \sigmabar \in \overline{\Ccal_m'} \}$
of $S(\Xcal')$. As we have seen in \ref{setup for degree}, a vertex may only occur in the interior of a canonical simplex $\Delta_S$ with $\dim(\overline{f}_{\rm aff}(\Delta_S))=\dim(\Delta_S)$. By \eqref{limit of discrete measures}, we conclude that the support of $\mu$ is contained in the union of such $\Delta_S$.

\vspace{2mm}
\noindent {\it Step 3:  Transformation of the limit in \eqref{limit of discrete measures} into a multiple of $\vol(\Omega)$.}
\vspace{2mm}

\noindent To prove \eqref{explicit formula}, we may assume that $\Omega$ is a polytope contained in the interior of a canonical simplex $\Delta_S$ with $\dim(\overline{f}_{\rm aff}(\Delta_S))=\dim(\Delta_S)$. Using the identification $\Delta_S=\Sigma_S$, the lift $f_{\rm aff}:\Delta_S \rightarrow \rdop^n$ extends to an affine map $f_0:\rdop^r \rightarrow \rdop^n$ which is also one-to-one and the  polytopal decomposition $\Dcal:=\{f_0^{-1}(\Delta) \mid \Delta \in \Ccal_1\}$ is periodic with respect to the lattice $\Lambda_S$ from \ref{setup for explicit}. Similarly as in \eqref{decompositon identity}, we have 
$\textstyle{  \Dcal_m = \{\Delta_S \cap \overline{f}_{\rm aff}^{-1}(\sigmabar) \mid S \in \str(\tilde{\Xcal}'),\, \sigmabar \in \overline{\frac{1}{m}\Ccal_1} \}}$. 
We conclude that there is a bijective correspondence between the irreducible components $Z$ of $\tilde{\Xcal}_m''$  with $\xi_Z \in \Omega$ and the vertices $\ub'$ of $\frac{1}{m}\Dcal$ contained in $\Omega$. We note that our situation matches with \ref{setup for degree}. By our above choice of $\Sigmabar$, the transversality assumption \eqref{transversality assumption} in the vertex $f_0(\ub')$ follows easily from $\Sigmabar$-transversality in Lemma \ref{totally irrational construction}. From Proposition \ref{degree and dual polytope}, we get
\begin{equation*} 
\deg_{\tilde{\Lcal}_m''}(Z)=\frac{d!}{e!} \cdot \deg_{\Hcal^{\otimes m^2}}(\overline{S}) \cdot \vol(\{\ub'\}^{g_m}),
 \end{equation*}
where $e:=\dim(S)=d-r$ and $g_m:=f_{\Lcal_m} \circ f_0$. We deduce that
\begin{equation} \label{degree formula in 3rd step}
\deg_{\tilde{\Lcal}_m''}(Z)= \frac{d!}{e!} \cdot \deg_\Hcal(\overline{S}) \cdot \vol(\{\ub'\}^{g_m}) \cdot m^{2e}.
 \end{equation}
We define the dual polytope of the vertex $\ub := m \ub'$ of $\Dcal$ with respect to the convex function $g:=f_\Lcal \circ f_0:\rdop^r \rightarrow \rdop$ by 
$$\{\ub \}^g :=  \{\omega \in \rdop^r \mid \omega \cdot (\wb- \ub) \leq g(\wb)-g(\ub) \; \forall \wb \in U\},$$
where $U$ is a sufficiently small neighbourhood of $\ub$ in $\rdop^r$. Since  $\{\ub'\}^{g_m}=m\{\ub\}^g$, formula \eqref{degree formula in 3rd step} yields
\begin{equation} \label{second degree formula in 3rd step}
\deg_{\tilde{\Lcal}_m''}(Z)=\frac{d!}{e!} \cdot \deg_{\Hcal}(\overline{S}) \cdot \vol(\{\ub\}^{g})\cdot m^{2e+r}.
\end{equation}
Let $F$ be the fundamental domain of the lattice $\Lambda_S$ in $\rdop^r$. For $m \gg 0$, the number of $\frac{1}{m} \Lambda_S$-translates of $\frac{1}{m}F$ contained in $\Omega$ (resp. intersecting $\partial \Omega)$ is $m^r \vol(\Omega)/\vol(F)+O(m^{r-1})$ (resp. $O(m^{r-1})$). By \eqref{limit of discrete measures} and \eqref{second degree formula in 3rd step}, we deduce
\begin{equation} \label{first formula for the measure in terms of the volume}
\mu(\Omega)=\frac{d!}{e!} \cdot \deg_\Hcal(\overline{S}) \cdot \sum_\ub \vol(\{\ub\}^g) \cdot \frac{\vol(\Omega)}{\vol(\Lambda_S)},
\end{equation}
where $\ub$ ranges over all vertices of $\Dcal$ modulo $\Lambda_S$. 
The set $\{\{\ub\}^g \mid \text{$\ub$ vertex of $\Dcal$}\}$ is invariant under $\Lambda_S$-translation. By \cite{McM}, Theorem 3.1, this set is a $\Lambda_S^L$-periodic tiling of $\rdop^r$ which means that $\rdop^r$ is covered by these $r$-dimensional polytopes and they meet face-to-face. Together with \eqref{first formula for the measure in terms of the volume}, this proves \eqref{explicit formula}. Since $\tilde{\Hcal}$ is ample on $\tilde{\Bcal}$, we have $\deg_\Hcal(\overline S) \neq 0$ if and only if $\Delta_S$ is non-degenerate with respect to $f$. By step 2, we get also the claim about the support. \qed 

\begin{rem} \rm \label{canonical subset}
By the projection formula (b) in Proposition \ref{Chambert-Loir's measures}, Theorem \ref{explicit formula for canonical measure} gives also an explicit description for the canonical measure $$c_1(\overline{L}|_X)^{\wedge d}=f_*(c_1(f^*\overline{L})^{\wedge d})$$ on $X$. We conclude that  the support of such a canonical measure is equal to the union of all $f(\Delta_S)$, where $\Delta_S$ ranges over all canonical simplices of $S(\Xcal')$ which are non-degenerate with respect to $f$. Note that this set is independent of the choice of $L$. We call it the {\it canonical subset} of $\Xan$. 
\end{rem} 

The referee has suggested that the canonical subset is a piecewise linear space. We refer to \cite{Ber5}, chapter 1, for  the definition of a piecewise $R_{\zdop_+}$-piecewise linear spaces for $R:=\qdop \cap (0,1]$. We will always skip  $R_{\zdop_+}$ for briefity.

\begin{thm} \label{piecewise linear}
The canonical subset  of $\Xan$ has a unique structure as a piecewise linear space $T$ such that for any semistable alteration $\varphi_0:\Xcal' \rightarrow \Acal_0$ as in \ref{Semistable alteration} with generic fibre $f:X'\rightarrow A^{\rm an}$, the restriction of $f$ to the union of all canonical simplices which are non-degenerate with respect to $f$ induces a piecewise linear map to $T$ with finite fibres.
\end{thm}

\proof 
Let $\Xcal_0$ be the closure of $X$ in a Mumford model $\Acal_0$ of $A$ over $\kcirc$ associated to the rational polytopal decomposition $\overline{\Ccal_0}$ of $\rtor$.
By a result of de Jong, there is a finite group $G$ acting on a strongly non-degenerate pluristable formal scheme $\Ycal$ over $\kcirc$ with the following properties (see 
\cite{Ber4}, Lemma 9.2): 
\begin{itemize}
\item[(a)] We endow $\Xcal_0$ with the trivial $G$-action. Then there is a dominant $G$-equivariant morphism $\gamma:\Ycal \rightarrow \Xcal_0$.
\item[(b)] The generic fibre $Y$ of $\Ycal$ is the analytic space associated to an irreducible smooth projective variety over $\kdop$.
\item[(c)] The generic fibre $g:Y \rightarrow \Xan$ of $\gamma$ is a generically finite proper morphism.
\item[(d)] The fixed field $\kdop(Y)^G$ is a purely inseparable extension of the field of rational functions $\kdop(X)$.
\end{itemize} 

Now we choose a semistable alteration $\eta:\Xcal' \rightarrow \Ycal$ with generic fibre $h:X' \rightarrow Y$. Then $\varphi_0:=\gamma \circ \eta$ plays the role of the semistable alteration in \ref{Semistable alteration} and $f:=g \circ h$ is its generic fibre.

Let $\Delta_S$ be a canonical simplex of $S(\Xcal')$ which is non-degenerate with respect to $f$. Since $\overline{f}_{\rm aff} \circ \Val=\valbar \circ g \circ h$, it is clear that $h$ is one-to-one on $\Delta_S$. We claim that $h(\Delta_S)$ is contained in the skeleton $S(\Ycal)$ of $\Ycal$. 


By continuity, it is enough to prove that $h(\ub') \in S(\Ycal)$ for every $\ub' \in \relint(\Delta_S)$ with rational coordinates. We choose a rational polytopal decomposition $\overline{\Ccal_1}$ of $\rtor$ with associated Mumford model $\Acal_1$ such that $\ub'$ is a vertex of the subdivision $\Dcal:=\{\Delta_R \cap \overline{f}_{\rm aff}^{-1}(\sigmabar) \mid R \in \xstr, \sigmabar \in \overline{\Ccal_1}\}$ satisfying the transversality condition \eqref{transversality assumption} and such that $g:=f_\Lcal \circ f_{\rm aff}$ is a strictly convex polyhedral function  in  $\ub'$ for a symmetric ample line bundle $L$ on $A$ with a formal $\kcirc$-model $\Lcal$ on $\Acal_1$ associated to a piecewise affine function $f_\Lcal$ as in Proposition \ref{Mumford line bundles}. This is much easier to construct than the simultaneous transversality conditions in Lemma \ref{totally irrational construction} and does not require a base change. 

Let $\Acal_1'$ be the Mumford model associated to $\overline{\Ccal_0} \cap \overline{\Ccal_1}$. 
We get a commutative  diagram of admissible formal schemes over $\kcirc$ with reduced special fibres 
\begin{equation}  \label{diagram of can morphisms 2}
\begin{CD} 
\Xcal'' @>\eta_1>> \Ycal_1 @>{\gamma_1}>> \Acal_1' @>{\psi_1}>> \Acal_1\\
@VV{\iota_1'}V    @VV{j_1}V  @VV{\iota_1}V\\
\Xcal' @>{\eta}>> \Ycal @>{\gamma}>> \Acal_0
\end{CD}
\end{equation}
by assuming that the rectangles are cartesian on the level of formal analytic varieties. The vertical maps and $\psi_1$ are the identity on the generic fibre. 

By Corollary \ref{strata corollary}, there is a unique irreducible component $Z$ of $\tilde{\Xcal}''$ with $\ub'=\xi_Z$. Since the assumptions of \ref{setup for degree} are satisfied, Proposition \ref{degree and dual polytope} yields
$$\deg_\Lcal(Z)= \frac{d!}{e!} \cdot \deg_\Hcal(\overline{S}) \cdot \vol(\{\ub'\}^g).$$
Since $\Hcal$ is ample  (see \ref{line bundles}) and $\Delta_S$ is non-degenerate with respect to $f$, we have $\deg_\Hcal(\overline{S})>0$. By strict convexity of $g$ in $\ub'$, we get also $\vol(\{\ub'\}^g)>0$ and hence $\deg_\Lcal(Z)>0$. By projection formula, we have $\deg_\Lcal(Z)=\deg_\Lcal(\tilde{\beta}_*(Z))$ for $\beta:=\psi_1\circ \gamma_1\circ \eta_1$. Note that $\dim(\tilde{\beta}(Z))=\dim(Z)=d$ is necessary for the positivity of the degree.  Now \eqref{diagram of can morphisms 2} yields that $\tilde{\eta}_1(Z)$ is also $d$-dimensional. Since $Y$ is $d$-dimensional, we conclude that $\tilde{\eta}_1(Z)$ is  an irreducible component $W$ of $\tilde{\Ycal}_1$ and hence $h(\xi_Z)=\xi_W$. By the generalization of Corollary \ref{strata corollary}(g) to $\Ycal$ (see Lemma \ref{pluristable} below), we know that $\xi_W$ is a vertex of a subdivision $\Dcal_1$ of $S(\Ycal)$ and hence $h(\ub')=\xi_W \in S(\Ycal)$ proving $h(\Delta_S) \subset S(\Ycal)$. 


By Remark \ref{canonical subset} and the above,  $h$ maps the support of $\mu:=c_1(f^*\overline{L}^{\wedge d})$ into $S(\Ycal)$. By \cite{Ber5}, Corollary 6.1.3, $h$ restricts to a piecewise linear map from the piecewise linear subspace $\supp(\mu)$ to $S(\Ycal)$. Moreover, the skeleton $S(\Ycal)$ is invariant under $G$ and the $G$-transformations induce piecewise linear automorphisms of the skeleton (\cite{Ber5}, Corollary 6.1.2). 

There is a Zariski dense  open subset $U$ of $\Xan$ such that $g:V \rightarrow U$ is finite for $V:=g^{-1}(U)$ (see (c)). By \cite{Ber4}, Corollary 8.6, the quotient $V/G$ exists. By \cite{Ber4}, Corollary 8.4, we have $S(\Ycal) \subset V$. We note that the compact subset  $S(\Ycal)/G$  of  $V/G$ has a canonical structure as a piecewise linear space. Indeed, the skeleton $S(\Ycal)$ is a piecewise linear space because it is the geometric realization of a polysimplicial set $D$ (see \cite{Ber5}, Theorem 5.1.1). As $S(\Ycal)/G$ is the geometric realization of the polysimplicial set $D/G$, we deduce that $S(\Ycal)/G$ is also a piecewise linear space (see \cite{Ber5}, Proposition 3.5.3). 
We conclude that $h(\supp(\mu))$ is a piecewise linear subspace of $S(\Ycal)$ which maps onto a piecewise linear subspace of $S(\Ycal)/G$. By shrinking $U$ and using (d), we may assume that the canonical morphism $V/G \rightarrow U$ is radicial. In particular, it is a homeomorphism of the underlying topological spaces (see \cite{Ber6}, Remark 2.2.2). As a consequence, we get a piecewise linear structure on $f(\supp(\mu))=g(h(\supp(\mu)))$. By Remark \ref{canonical subset}, this  is the canonical subset $T$ of $\Xan$. 

The domains of linearity for the piecewise linear map $f:\supp(\mu) \rightarrow T$ are subsets of the canonical simplices of $S(\Xcal')$ which are non-degenerate with respect to $f$. By Proposition \ref{refinement and formal analytic structure}, they induce a finer formal analytic structure on $X'$ and we may apply de Jong's alteration theorem also to the associated formal scheme over $\kcirc$. Replacing the alteration $\eta$ by the composition of the two alterations, we may assume that the domains of linearity are really equal to the canonical simplices of $S(\Xcal')$ which are non-degenerate with respect to $f$. Then a linear atlas of $T$ is given by the charts $f(\Delta_S)$, where $\Delta_S$ ranges over all such canonical simplices, and $f$ is a linear isomorphism from $\Delta_S$ onto $f(\Delta_S)$. 

Let us consider now any semistable alteration $\varphi_0':\Zcal' \rightarrow \Acal_0'$ as in \ref{Semistable alteration} with generic fibre $f':Z'\rightarrow A^{\rm an}$. Then there is a semistable alteration which factors through $\varphi_0$ and $\varphi_0'$. Using the above atlas, it follows easily that $f'$ induces a piecewise linear map $S(\Zcal')\rightarrow T$. Uniqueness is obvious. \qed

\vspace{3mm} We sketch now how the results of the first part of Section 5 generalize to the strongly non-degenerate pluristable formal scheme $\Ycal$ from the above proof. By definition, there is a strongly non-degenerate strictly pluristable formal scheme $\Ycal'$ and a surjective \'etale morphism $\rho:\Ycal' \rightarrow \Ycal$. By \cite{Ber5}, Theorem 5.1.1, the skeleton $S(\Ycal)$ has a piecewise linear structure which is a cokern of the piecewise linear structure on $S(\Ycal')$ described in Remark \ref{strictly pluristable} and in the appendix. Moreover, there is ``polytopal'' subdivision of $S(\Ycal)$ given by canonical ``plurisimplices'' $\Delta_S$ which are in bijective correspondence to the strata $S$ of $\tilde{\Ycal}$. We use quotation marks because $\Delta_S$ is only a quotient of a canonical plurisimplex $\Delta_{S'}$ of $S(\Ycal')$ for any stratum $S'$ with $\tilde{\rho}(S') \subset S$. To construct $\Delta_S$, we have to identify closed faces $\Delta_P$ and $\Delta_Q$ in the boundary of $\Delta_{S'}$ if and only if the strata $P$ and $Q$ map into the same stratum of $\tilde{\Ycal}$. By \cite{Ber5}, there are well defined proper strong deformation retractions $\Val:Y \rightarrow S(\Ycal)$ and $\Val':Y' \rightarrow S(\Ycal')$, where $Y$ and $Y'$ are the generic fibres of $\Ycal$ and $\Ycal'$. 

\begin{lem} \label{pluristable}
There is a unique map $\overline{g}_{\rm aff}:S(\Ycal)\rightarrow \rtor$ with $\overline{g}_{\rm aff} \circ \Val = \valbar \circ g$ on $Y$. We get a ``polytopal'' subdivision $\Dcal_1:=\{\Delta_S \cap \overline{g}_{\rm aff}^{-1}(\sigmabar) \mid S \in \str(\tilde{\Ycal}), \sigmabar \in \overline{\Ccal_1}\}$ of $S(\Ycal)$ defining a formal analytic structure $\Y_1$ on $Y$ as in Proposition \ref{refinement and formal analytic structure} with $\Y_1 = \Ycal_1^{\rm f-an}$. Moreover, Proposition \ref{orbits and refinement} and Corollary \ref{strata corollary}(a),(b),(e),(f),(g) hold for $\Y_1$ (instead of $\X''$).
\end{lem}

\proof By Remark \ref{strictly pluristable}, the lemma holds for the strictly pluristable $\Ycal'$. The idea is now to use the \'etale covering $\rho:\Ycal' \rightarrow \Ycal$ to deduce  the claim for $\Ycal$. It is necessary to define $\overline{g}_{\rm aff}:=\valbar \circ g$. Using $\rho^{\rm an}$  surjective (\cite{Ber3}, Lemma 2.2) and ${\rho^{\rm an}} \circ \Val'=\Val \circ {\rho^{\rm an}}$, it is easy to prove $\overline{g}_{\rm aff} \circ \Val =\valbar \circ g$ from the corresponding property for $\overline{g}_{\rm aff}':S(\Ycal')\rightarrow \rtor$. 

In Appendix A, we have studied building blocks for strongly non-degenerate strictly pluristable formal schemes over $\kcirc$. We define a {\it building block} for $\Ycal$ as a formal affine open subscheme $\Ucal$ of $\Ycal$ such that a building block $\Ucal'$ of $\Ycal'$ exists with $\rho(\Ucal')=\Ucal$. Since an \'etale map is open, the building blocks cover $\Ycal$. By definition, $\tilde{\Ucal}'$ has a smallest stratum $S'$ and hence $\tilde{\Ucal}$ has the smallest stratum $S:=\tilde{\rho}(S')$. We set $U:=\Ucal^{\rm an}$. If $\Ucal$ varies over all building blocks and $\Delta$ over $\Dcal_1$, then $U \cap \Val^{-1}(\Delta)$ is a formal affinoid atlas for $Y$. Here, we use that the ``polytopal subdivision'' $\Dcal_1$ is induced by the valuation of units, i.e. $\Delta$ is given by inequalities $v(b_\mb) + \mb \cdot \ub \geq 0$ induced by the units $g^*(b_\mb \xb^\mb)$ on $U$ coming from the inequalities of a corresponding polytope of $\overline{\Ccal_1}$. It follows now analogous to the proofs of Proposition \ref{refinement and formal analytic structure} and \ref{formal refinement for alteration} that the formal affinoid atlas induces a formal analytic variety $\Y_1$ isomorphic to $\Ycal_1^{\rm f-an}$.

A similar construction applies to $\Ycal'$ leading to an admissible formal scheme $\Ycal_1'$ over $\kcirc$ with reduced special fibre. Since $\rho$ is \'etale and surjective, the natural base change $\rho_1$ of $\rho$ to $\Ycal_1$ is \'etale and surjective. We deduce from \ref{admissible formal schemes} that $\rho_1$ is the  canonical map $\Ycal_1' \rightarrow \Ycal_1$. 

Let $R'$ be a stratum of $\tilde{\Ycal}_1'$. By \cite{Ber4}, Lemma 2.2, $\tilde{\rho}_1(R')$ is an open dense set of a stratum of $\tilde{\Ycal}_1$. Using the claim for $\Ycal'$, we know that $R'=\pi((\Val')^{-1}(\tau'))$ for a unique open face $\tau'$ of $\Dcal_1'$. Note that ${\rho^{\rm an}}$ maps $\tau'$ isomorphically onto an open face $\tau$ of $\Dcal_1$ and $\tilde{\rho}_1(R') \subset \pi(\Val^{-1}(\tau))$.

If $\tau$ varies over all open faces of $\Dcal_1$, then the definition of $\Y_1=\Ycal_1^{\rm f-an}$ yields that $\pi(\Val^{-1}(\tau))$ is a partition of $\tilde{\Y}_1=\tilde{\Ycal}_1$. The surjectivity of $\tilde{\rho}_1$ and a partition argument show that $R:=\pi(\Val^{-1}(\tau))$ is a strata subset equal to  $\bigcup\tilde{\rho}_1(\pi((\Val')^{-1}(\tau')))$, where $\tau'$ ranges over all open faces of $\Dcal_1'$ with ${\rho^{\rm an}}(\tau')=\tau$. 

We claim that $R$ is a stratum of $\tilde{\Ycal}_1$. By localizing, we may assume that $\Ycal$ and $\Ycal'$ are building blocks and that $\tau \subset \relint(\Delta_S)$ for the canonical ``plurisimplex'' $\Delta_S=S(\Ycal)$. Since ${\rho^{\rm an}}$ maps the interior of the plurisimplex $\Delta_{S'}=S(\Ycal')$ isomorphically onto $\relint(\Delta_S)$, there is a unique open face $\tau'$ of $\Dcal_1'$ with ${\rho^{\rm an}}(\tau')=\tau$. We conclude that $R=\tilde{\rho}_1(\pi((\Val')^{-1}(\tau'))) \in \str(\tilde{\Ycal}_1)$.

The remaining claims are easily deduced from the corresponding claims for $\Ycal'$ using that $\rho_1$ is \'etale and surjective. \qed

\section{Proof of the main theorem and examples}

First, we will give the proof of Theorem \ref{main theorem}.  Then we will describe the canonical measure  in two relevant examples.

\vspace{3mm}
\noindent {\bf Proof of Theorem \ref{main theorem}:\:}
By a finite base change and using linearity in the components, we may assume that $X$ is a $d$-dimensional geometrically integral closed subvariety of the abelian variety $A$. The argument will be based on the description of canonical measures in the previous section, hence we will use the notation from there.
By multilinearity, it is enough to consider ample line bundles in (c)  of Theorem \ref{main theorem}. We choose a semistable alteration $\varphi_0:\Xcal' \rightarrow \Xcal_0$ as in \ref{Semistable alteration} with generic fibre $f:X' \rightarrow \Xan$. Then we have the explicit description \eqref{explicit formula for several line bundles} of $\mu:=c_1(f^*(\overline{L_1})) \wedge \dots \wedge c_1(f^*(\overline{L_d}))$. Since $\mu$ is supported in $S(\Xcal')$, we get
\begin{equation} \label{tropical measure}
\valbar_* \left( c_1(\overline{L_1}|_{\Xan} )\wedge \dots \wedge c_1(\overline{L_d}|_{\Xan} ) \right) = \deg(f) \cdot \left(\overline{f}_{\rm aff} \right)_*(\mu)
\end{equation}
by Propositions \ref{Chambert-Loir's measures} and \ref{f_aff}. More precisely, Remark \ref{generalization to several line bundles} shows that $\mu$ is supported in the union of the canonical simplices $\Delta_S$ which are non-degenerate with respect to $f$. Since 
$$d-\dim(\Delta_S)=\dim(S)=\dim\left(\tilde{q}_0 \circ \tilde{\Phi}_0(S)\right) \leq \dim(\tilde{\Bcal})=b,$$
we get $\dim(\overline{f}_{\rm aff}(\Delta_S)=\dim(\Delta_S) \geq d-b$. By Theorem \ref{dimension of tropical variety}, the tropical variety $\valbar(\Xan)$ is a finite union of rational polytopes of dimension at most $d$ and at least $d- \min\{b,d\}\}$. Hence we may list the simplices $\overline{f}_{\rm aff}(\Delta_S)$ as in (a), where $\Delta_S$ ranges over the canonical simplices which are non-degenerate with respect to $f$. Then (c) and (d) follow from  \eqref{explicit formula for several line bundles} and \eqref{tropical measure}. Finally, (b) follows from the next lemma. \qed

\begin{lem} \label{tropical variety and non-degenerate simplices}
In the above notation, let us consider $\ubb \in \valbar(\Xan)$ and let $d-e$ be the dimension of the tropical variety $\valbar(\Xan)$ in a neighbourhood of $\ubb$. Then there is a $(d-e)$-dimensional canonical simplex $\Delta_S$ of $S(\Xcal')$ which is non-degenerate with respect to $f$ such that $\ubb \in  \overline{f}_{\rm aff}(\Delta_S)$.
\end{lem}

\proof By Proposition \ref{f_aff} and the surjectivity of $f$ and $\valbar$, it is clear that the simplices $\overline{f}_{\rm aff}(\Delta_T)$, $T \in \str(\tilde{\Xcal}')$, cover $\valbar(\Xan)$. 
Since $\overline{f}_{\rm aff}$ is locally an affine map defined over $\qdop$, we may assume that $\ubb$ is an element of  $\valbar(\Xan)$ with coordinates in $\qdop$. Moreover, by density of the $\qdop$-rational points in $\valbar(X^{\rm an})$, we may assume that $\ubb$  is not contained in any $\overline{f}_{\rm aff}(\Delta_T)$ of dimension $<d-e$.  We have to prove that there is a $(d-e)$-dimensional canonical simplex $\Delta_S$   with $\ubb \in \overline{f}_{\rm aff}(\Delta_S)$ and 
\begin{equation} \label{dimension and e}
\dim\left( \tilde{q}_0 \circ \tilde{\Phi}_0(\overline{S}) \right) =e,
\end{equation}
which yields that $\Delta_S$ is non-degenerate with respect to $f$.
We choose a rational simplex $\Deltabar$ in $\rtor$ of codimension $d-e$ such that $\ubb \in \relint(\Deltabar)$ and 
\begin{equation} \label{transvers}
\Deltabar \cap \valbar(\Xan) = \{\ubb\}.
\end{equation}
We extend $\Deltabar$ to a rational polytopal subdivision $\Ccalbar_1$ of $\rtor$ (which is not assumed to have the properties of Lemma \ref{totally irrational construction}). We denote the associated Mumford model of $A$ by $\Acal_1$. 

Let $\Ccalbar:= \overline{\Ccal_0} \cap \overline{\Ccal_1}$ be the minimal polytopal  decomposition of $\rtor$ containing $\overline{\Ccal_0}$ and $\overline{\Ccal_1}$. 
Let $\Acal$ be the Mumford model of $A$ associated to $\Ccalbar$ and let $\phi:\X'' \rightarrow \Acal^{\rm f-an}$ be the morphism obtained from $\varphi_0$ by base change as in Propositions \ref{formal refinement for alteration} and \ref{minimal formal analytic structure}. Since $\Ccalbar$ is a polytopal subdivison of $\Ccalbar_1$, we get a canonical formal analytic morphism $\phi_1:\X'' \rightarrow \Acal_1^{\rm f-an}$. Passing to the associated admissible formal schemes over $\kcirc$, this induces a morphism $\varphi_1:\Xcal'' \rightarrow \Acal_1$ (see \ref{admissible formal schemes}). 

Note that $U:=\valbar^{-1}(\Deltabar)$ is a formal open subset of $\Acal_1^{\rm f-an}$. By \eqref{transvers}, we get $\Xan \cap U \neq \emptyset$. Let $\Xcal_1$ be the closure of $\Xan$ in $\Acal_1$. Then the special fibre $\Xcal_1$ has an irreducible component $Y$ with $Y \cap \tilde{U} \neq \emptyset$. We use here that the reduction of $\Acal_1^{\rm f-an}$ is equal to the special fibre of $\Acal_1$ (see \ref{admissible formal schemes}). Since $\tilde{\varphi}_1$ maps $\tilde{\Xcal}''$ onto  $\tilde{\Xcal}_1$, there is an irreducible component $Y_{\ub'}$ of $\tilde{\Xcal}''$ mapping onto $Y$. As the notation already indicates, $Y_{\ub'}$ is the irreducible component associated to a vertex $\ub'$ of the rational subdivision 
$\Dcal= \{ \Delta_S \cap \overline{f}_{\rm aff}^{-1}(\sigmabar) \mid S \in \xstr, \, \sigmabar \in \Ccalbar\}$
of $\xskel$ (see Corollary \ref{strata corollary}). More precisely, $\ub'=\Val(\xi)$ for the unique point $\xi$ of $X'$ which reduces to the generic point of this irreducible component $Y_{\ub'}$. From $\tilde{\varphi}_1(Y_{\ub'})=Y$, we deduce that the reduction of $f(\xi)$ to the special fibre $\tilde{\Acal}_1$ is equal to the generic point of $Y$. We conclude that $f(\xi)\in U$ and hence $\overline{f}_{\rm aff}(\ub')=\valbar(f(\xi)) \in \Deltabar$ by Proposition \ref{f_aff} and by the definition of $U$. Since $f(\xi) \in X^{\rm an}$, we get  $\overline{f}_{\rm aff}(\ub')=\ubb$ from \eqref{transvers}.
 
 Let $S$ be the unique stratum of the chosen strictly semistable $\kcirc$-model $\Xcal'$ of $X'$ with $\ub' \in \relint(\Delta_S)$. We note first that $\ubb=\overline{f}_{\rm aff}(\ub') \in \overline{f}_{\rm aff}(\Delta_S)$. As we have remarked at the end of \ref{setup for degree}, the fact that $\relint(\Delta_S)$ contains a vertex of $\Dcal$ implies that $\dim(\overline{f}_{\rm aff}(\Delta_S))=\dim(\Delta_S)$. From the non-degeneracy assumption on $\ubb$, we deduce that $\dim(\Delta_S)=d-e$ and hence $S$ is an $e$-dimensional stratum. 

It remains to prove \eqref{dimension and e}. By Corollary \ref{strata corollary}, the canonical morphism $\tilde{\Xcal}'' \rightarrow \tilde{\Xcal}'$ maps $Y_{\ub'}$ onto $\overline{S}$. By 
 Lemma  \ref{lift of phi}, we have  lifts $\tilde{\Phi}_0:\overline{S} \rightarrow \tilde{\Ecal}_0$ and $\tilde{\Phi}:Y_{\ub'}\rightarrow \tilde{\Ecal}$ of $\tilde{\varphi}_0$ and $\tilde{\varphi}$, where $\Ecal_0$ and $\Ecal$ are the $\kcirc$-models of the uniformization $E$ of $A$ associated to the polytopal decompositions $\Ccal_0$ and $\Ccal$. 
Using  that $\Ccal$ is a polytopal subdivision of $\Ccal_1$, the map $\tilde{\Phi}$ induces a  canonical morphism $\tilde{\Phi}_1:Y_{\ub'}\rightarrow \tilde{\Ecal}_1$ which is a lift of the restriction of $\tilde{\varphi}_1$ to ${Y_{\ubb'}}$. This lift may be also constructed by the fact that $\Ecal_1$ and $\Acal_1=\Ecal_1/M$ are locally isomorphic. 
 Since $\tilde{\varphi}_1(Y_{\ub'})=Y$, Proposition \ref{strata on Mumford models}(d) yields that $Y':=\tilde{\Phi}_1(Y_{\ub'})$ is isomorphic to $Y$. Lemma  \ref{lift of phi} and an easy diagram chase involving the canonical morphisms $q_i:\Ecal_i\rightarrow \Bcal$ to the formal abelian scheme $\Bcal$ of the Raynaud extension  show 
$$\tilde{q}_0 \circ \tilde{\Phi}_0(\overline{S})= \tilde{q} \circ \tilde{\Phi}(Y_{\ub'})=\tilde{q}_1 \circ \tilde{\Phi}_1(Y_{\ub'}).$$
We conclude that the dimension of
\begin{equation} \label{strata and component}
\tilde{q}_0 \circ \tilde{\Phi}_0(\overline{S})= \tilde{q}_1(Y')
\end{equation}
is at most $e$-dimensional. To show equality, we consider a basic formal affinoid subdomain $U_{V,\Delta}\cong V \times U_\Delta$ from the construction of $\Ecal_1^{\rm f-an}$ (see \ref{Mumford's construction}). Here, $V$ is the generic fibre of a formal affine open subset  of $\Bcal$ which trivializes the Raynaud extension \eqref{Raynaud extension 1} of $A$ and $\Delta$ is a simplex in $\rdop^n$ lifting the simplex $\Deltabar$ considered at the beginning of the proof. We may choose $V$ such that $\tilde{U}_{V,\Delta} \cap Y' \neq \emptyset$. Recall that $U_\Delta$ is  the polytopal subdomain $\val^{-1}(\Delta)$ of $({\mathbb G}_m^n)^{\rm an}$. By the choice of $V$, we have $\tilde{U}_{V,\Delta} \cong \tilde{V} \times \tilde{U}_\Delta$. 

We claim that the second projection $\tilde{p}_2$ maps the generic point of $Y'$ into the torus orbit $Z$ of $\tilde{U}_\Delta$ corresponding to $\relint(\Delta)$. Recall from \cite{Gu4}, Proposition 4.4, that $Z=\pi(\val^{-1}(\Delta))$, where $\pi$ is the reduction map. We have seen above that $f(\xi)$ reduces to the generic point of $Y$ and that $\valbar(f(\xi))=\ubb$. Let $\xi'$ be the unique lift of $f(\xi)$ to $E$ whose reduction $\pi(\xi')$ is the generic point of $Y'$. We conclude that $\val(p_2(\xi'))=\val(\xi')$ is the unique point $\ub \in \relint(\Delta)$ which lifts $\ubb$. Therefore we have $\pi(p_2(\xi')) \in Z$.  
Since $\pi(p_2(\xi'))=\tilde{p}_2(\pi(\xi'))$, we get the above claim. Since $Z$ is the closed orbit of $\tilde{U}_\Delta$, we conclude that $\tilde{p}_2(Y') \subset Z$.

By \cite{Gu4}, Proposition 4.4, the dimension of the torus orbit $Z$ is $\codim(\Delta)=d-e$. The above claim shows that $Y'\cap \tilde{U}_{V,\Delta}$ is contained in the closed subset $(\tilde{q}_1(Y')\cap \tilde{V})\times Z$ of $\tilde{U}_{V,\Delta} \cong \tilde{V} \times \tilde{U}_\Delta$. Since $Y'$ is $d$-dimensional and the product is at most $d$-dimensional, we get 
$Y' \cong (\tilde{q}_1(Y')\cap \tilde{V}) \times Z$. 
Moreover, we deduce $\dim(\tilde{q}_1(Y'))=e$. By \eqref{strata and component}, we get \eqref{dimension and e} and the lemma. \qed

\vspace{3mm}
\noindent {\bf Remark to the proof of Lemma \ref{tropical variety and non-degenerate simplices}:\:}
The argument shows that the irreducible component $Y'$ of $\tilde{\Ecal}_1$ is a fibre bundle over $\tilde{q}_1(Y')$ with fibre isomorphic to the toric variety $Y_\Delta$ associated to ${\rm star}(\Delta)=\{\sigma\in \Ccal_1 \mid \Delta \subset \sigma\}$ (see \cite{Gu4}, Remark 4.8). Indeed, the choice of $V$ as a trivialization yields that $\tilde{q}_1^{-1}(\tilde{V}) \cong \tilde{V} \times \tilde{\Zcal}$, where $\Zcal$ is the formal $\kcirc$-model of $({\mathbb G}_m^n)^{\rm an}$ associated to the polytopal decomposition $\Ccal_1$. Now the proof of the lemma shows that $Y' \cap \tilde{q}_1^{-1}(\tilde{V}) \cong ( \tilde{q}_1(Y') \cap \tilde{V}) \times Y_\Delta$ proving the claim.

\begin{ex} \rm \label{canonical measure on abelian variety}
Let us consider the special case $X= A$ in Theorem \ref{explicit formula for canonical measure}. For every $\ubb \in \rtor$, there is a canonical point $\xi_{\ubb} \in A^{\rm an}$ which we describe as follows: Let $V$ be the generic fibre of a non-empty formal affine open subset of the abelian scheme $\Bcal$ which trivializes the Raynaud extension \eqref{Raynaud extension 1} of $A$. Then $U_{V,\{\ub\}}=\val^{-1}(\ub) \cap q^{-1}(V) \cong V \times U_{\{\ub \}}$ is an affinoid subdomain of the uniformization $E$. Using $A^{\rm an} \cong E/M$, it is obvious that $U_{V,\{\ub\}}$ is isomorphic to an affinoid subdomain $U_{[V,\ub]}$ of $A^{\rm an}$. By Lemma \ref{relative polytopal domain}, we may write every analytic function $h$ on $U_{[V,\ub]}$ as a strictly convergent Laurent series
\begin{equation*} 
h= \sum_{\mb \in \zdop^n} a_\mb x_1^{m_1} \cdots x_n^{m_n} 
\end{equation*}
in the torus coordinates $x_1, \dots, x_n$ on the polytopal domain $U_{\{\ub\}}$ in $(\Tor)^{\rm an}$, where the $a_\mb \in \Ocal(V)$ are uniquely determined by $h$. Then we define $\xi_{\ubb} \in U_{[V,\ub]}$ by
$$|h(\xi_{\ubb})| = \sup_{\mb \in \zdop^n} |a_\mb|_{\rm sup}\cdot e^{-\mb \cdot \ub}.$$
It is easy to see that $\xi_\ubb$ does not depend on the choice of $V$ and the representative $\ub$. The subset $S(A):=\{\xi_\ubb \mid \ubb \in \rtor\}$ of $A^{\rm an}$ is called the {\it skeleton} of $A$ (see \cite{Ber}, \S 6.5). 

By a combinatorial result of Knudson and Mumford (\cite{KKMS}, Chapter III), there is a rational triangulation $\Ccalbar$ of $\rtor$ (even refining any given rational polytopal decomposition) and $m_\Ccal \in \ndop \setminus \{0\}$ such that for every  maximal $\Delta \in \Ccal$, the simplex $m_\Ccal \Delta$ is ${\rm GL}(n,\zdop)$-isomorphic to a $\zdop^n$-translate of the standard simplex $\{\ub \in \rdop_+^n \mid u_1 + \dots + u_n \leq 1 \}$. Then the Mumford model $\Acal$ of $A$ associated to $\Ccalbar$ is strictly semistable.  K\"unnemann used this to construct projective strictly semistable $\kdop^\circ$-models for abelian varieties (see \cite{Ku1} and also the erratum in \cite{Ku2}, 5.8). 

Let $\Delta \in \Ccal$ with $\ub \in \Delta$. A similar application of Lemma \ref{relative polytopal domain} as above shows that $\xi_\ubb$ is contained in the affinoid chart $U_{[V,\Delta]} \cong V \times U_\Delta$ of $A^{\rm an}$ and that $|h(\xi_\ubb)| \geq |h(x)|$ for all $h \in \Ocal(U_{[V,\Delta]})$ and all $x \in U_{[V,\Delta]}$ with $\valbar(x)=\ubb$. By \cite{Ber4}, Theorem 5.2, this maximality implies that $\xi_\ubb$ is contained in the skeleton $S(\Acal)$ of $\Acal$. We conclude that $S(A)=S(\Acal)$ and $\valbar=\Val$ maps the skeleton homeomorphically onto $\rtor$.

We apply Theorem \ref{explicit formula for canonical measure} with $X'=X=A$ and $\Xcal'=\Acal$. The canonical simplices of $S(\Acal)$ are just the elements of $\Ccalbar$. By Proposition \ref{strata on Mumford models} and its proof,  a stratum $S$ of $\tilde{\Acal}$ has locally the form $S \cap \tilde{U}_{[V,\Delta_S]} \cong \tilde{V} \times Z_\tau$ for $V$ as above and the open face $\tau= \relint(\Delta_S)$  with corresponding stratum $Z_\tau:=\pi(\val^{-1}(\tau))$ in $\tilde{U}_{\Delta_S}$. Hence $\Delta_S$ is a non-degenerate simplex of $S(\Acal)$ (with respect to $f=\id$) if and only if $\dim(Z_\tau)=0$. We conclude that the non-degenerate canonical simplices of $S(\Acal)$ are just the $n$-dimensional simplices of $\Ccalbar$. The lattice $\Lambda_S^L$ in \ref{setup for explicit} does not depend on the choice of such a simplex $\Delta=\Delta_S$. By Proposition \ref{Chambert-Loir's measures} and Theorem \ref{explicit formula for canonical measure}, we conclude that $c_1( 
 \overline{L})^{\wedge d}$ is supported in $S(A)$ and corresponds to the unique Haar measure $\nu$ on $\rtor$ with $\nu(A)=\deg_L(A)$. Using multilinearity for non-ample line bundles and Remark \ref{generalization to several line bundles}, we deduce easily:
\end{ex}

\begin{cor} \label{cor for abelian varieties}
Let $\overline{L_1}, \dots , \overline{L_d}$ be canonically metrized line bundles on the a\-be\-lian variety $A$ over $K$ of dimension $d$. Then $c_1(\overline{L_1}) \wedge \cdots \wedge c_1(\overline{L_d})$ is supported in the skeleton $S(A)$ and corresponds to the Haar measure on $\rtor$ with total measure $\deg_{L_1, \dots, L_d}(A)$.
\end{cor}

\begin{ex} \rm \label{spectrum example} 
We will show that the whole spectrum of values $\{d-b, \dots, d-e\}$ in Theorem \ref{main theorem} may occur for a single canonical measure, where $d-e$ denotes the dimension of the tropical variety. We assume that $K$ is the function field $k(C)$ for an irreducible regular projective curve $C$ over an algebraically closed field $k$ of characteristic $0$. Let $v$ be the discrete valuation on $K$ defined by the order in a given closed point $P \in C$. It is easy to use the construction below to give similar examples for other fields.
 
We consider a product $A=B_1 \times B_2$ of abelian varieties over $K$, where $B_1$ has good reduction at $v$ and where $B_2$ is totally degenerate at $v$. As usual, let $\kdop$ be a minimal algebraically closed field containing $K$ which is complete with respect to a valuation extending $v$. The analytic considerations will be performed over $\kdop$. Totally degenerate at $v$ means that the Raynaud extension of $B_2$ is an analytic torus and hence $B_2^{\rm an} \cong (\Tor)^{\rm an}/M$ for a discrete subgroup $M$ with $\Lambda=\val(M)$ a complete lattice in $\rdop^n$. Then $E \cong B_1^{\rm an} \times (\Tor)^{\rm an}$ is the Raynaud extension of $A$ and we have $A^{\rm an} \cong E/M$.

By assumption, $B_1$ is the generic fibre of an abelian scheme $\Bmcal_1$ over the discrete valuation ring $K^\circ$. The associated admissible formal scheme $\Bcal_1:=\hat{\Bmcal}_1$ over $\kcirc$ (see \ref{schemes})) is just the formal abelian scheme $\Bcal$ over $\kcirc$ in the Raynaud extension \eqref{Raynaud extension 1} for $A$. 
To get a Mumford model $\Bcal_2$ for $B_2^{\rm an}$, we will use a similar polytopal decomposition $\Ccalbar$ of $\rtor$ as in Example \ref{canonical measure on abelian variety}. There is a rational triangulation 
$\Ccalbar$ of $\rtor$ such that the strictly semistable $\kdop^\circ$-model $\Bcal_2$ is projective (\cite{Ku1}, \S3, \S4). K\"unnemann's proofs show that $\Ccalbar$ can be chosen as a refinement of any given rational polytopal decomposition of $\rtor$ (see also \cite{Ku2}, 5.5). We get a strictly semistable formal $\kdop^\circ$-model $\Acal:=\Bcal_1 \times \Bcal_2$ of $A^{\rm an}$. 

By K\"unnemann's construction, $\Bcal_2$ is defined algebraically over the valuation ring $F^\circ$ of a finite extension $F$ of the completion $K_v$, i.e. we have a strictly semistable algebraic $F^\circ$-model $\Bmcal_2$ of $B_2$ with associated admissible formal scheme $\Bcal_2=\hat{\Bmcal}_2$. We choose  ample line bundles $\Lmcal_1$ on $\Bmcal_1$ and $\Lmcal_2$ on $\Bmcal_2$. Then  $\Lmcal:=p_1^*(\Lmcal_1) \otimes p_2^*(\Lmcal_2)$ is an ample line bundle on $\Amcal:=\Bmcal_1 \times_{K^\circ} \Bmcal_2$. By passing to a suitable tensor power of $\Lmcal$, we may assume that $\Lmcal$ is very ample and that $H^0(\Amcal, \Lmcal) \rightarrow H^0(\tilde{\Amcal},\tilde{\Lmcal})$ is surjective for the reduction $\tilde{\Lmcal}$ of $\Lmcal$ to the special fibre $\tilde{\Amcal}$.  

Let $b:=\dim(B_1)$ and let us fix $m \in \{0, \dots, \min(b,n)\}$. Let us choose generic global sections $\tilde{s}_1, \dots , \tilde{s}_m \in H^0(\tilde{\Amcal},\tilde{\Lmcal})$. By assumption, they are the reductions of global sections $s_1, \dots, s_m$ of $\Lmcal$. The generic choice of the sections leads to a closed subscheme 
$$\Xmcal := \Div(s_1) \cap \dots\cap \Div(s_m)$$
of  codimension $m$ in $\Amcal$ which is flat over $F^\circ$. By Bertini's theorem, the generic fibre $X$ of $\Xmcal$ is an irreducible smooth variety over $F$ of dimension $d:=b+n-m$. The same argument shows that the irreducible components $Y_i$ of the special fibre $\tilde{\Xmcal}$ are Cartier divisiors and $\cap_{i \in I} Y_i$ is a smooth variety over $\tilde{F}$ of pure dimension $\dim(X)-|I|+1$ for any non-empty subset $I$. By a criterion of Hartl and L\"utkebohmert (\cite{HL}, Proposition 1.3), $\Xmcal$ is strictly semistable. Since $m \leq b$, the fibre of $X$ over any point of $B_2$ is non-empty and hence
$$\valbar(\Xan)=\rtor.$$
We conclude that the excess $e$ of the tropical variety $\valbar(\Xan)$ is given by
$$e:=\dim(X)-\dim(\valbar(\Xan))=b-m.$$
Now we switch from the algebraic point of view to the analytic and formal category. Then $\Xmcal$ has an associated admissible formal scheme $\Xcal:=\hat{\Xmcal}$  over $\kcirc$ which is a strictly semistable formal $\kcirc$-model  of $\Xan$ and  a closed formal subscheme of the Mumford model $\Acal=\Bcal_1 \times \Bcal_2$ of $A^{\rm an}$. Let $\Lcal,\Lcal_1,\Lcal_2$ be the formal line bundles on $\Acal,\Acal_1,\Acal_2$ induced by $\Lmcal,\Lmcal_1,\Lmcal_2$.

If $S \in \str(\tilde{\Acal})$, then $S= \tilde{\Bcal}_1 \times S_2$ for $S_2 \in \str(\tilde{\Bcal}_2)$ corresponding to an open face $\overline{\tau}:=\relint(\Deltabar)$ for a unique $\Deltabar \in \Ccalbar$ (see Proposition \ref{strata on Mumford models}). We note that 
$$\dim(S) = \codim(\Delta, \rdop^n)+b \geq b \geq m.$$ We consider first the case $\dim(S)>m$. Using the generic choice of $s_1,\dots,s_m$ again,  Bertini's theorem yields that 
\begin{equation} \label{S'}
 S':= \Div(s_1) \cap \dots \cap \Div(s_m) \cap S
\end{equation}
is  a stratum of $\tilde{\Xcal}$ with $\dim(S')=\dim(S)-m$. If $\dim(S)=m=b$, then $S'$ is a strata subset of $\tilde{\Xcal}$ consisting of $\deg_{\tilde{\Lcal}_1}(\tilde{\Bcal}_1)=\deg_{L_1}(B_1)$ points. Therefore the skeleton $S(\Xcal)$ may be identified with the triangulation $\Ccalbar$ by using the map $\Val$ except in the case $m=b$ where we have to count the $n$-dimensional simplices of $\Ccalbar$ with multiplicity $\deg_{{L_1}}(B_1)$. By construction, $S'$ is non-degenerate (with respect to $f=\id$ in the sense of \ref{non-degenerate canonical simplex}) if and only if $\dim(S_2) \leq m$. If we endow the generic fibre $L$ of $\Lmcal$ with a canonical metric, then Theorem \ref{explicit formula for canonical measure} shows that
\begin{equation} \label{explicit formula in example}
\nu:= \valbar_*\left(c_1(\overline{L}|_X )^{\wedge d} \right) = \sum_\Deltabar \lambda_\Deltabar \cdot \delta_\Deltabar
\end{equation}
where $\Deltabar$ ranges over all simplices of $\Ccalbar$ with $\dim(\Deltabar) \geq n-m=d-b$ and where $\lambda_\Deltabar > 0$. Similarly as in Example \ref{canonical measure on abelian variety}, we deduce that the contribution of the $n$-dimensional simplices of $\Ccalbar$ to $\nu$ is equal to a strictly positive Haar measure $\nu_n$ on $\rtor$. Note however that for $m>0$, we have 
$$\nu_n(\rtor) < \nu(\rtor)=\deg_L(X).$$

Finally, we show that the multiplicities $\lambda_\Deltabar$ are given completely in terms of convex geometry. The simplex $\Deltabar \in \Ccalbar$ of dimension $r \geq n-m$ corresponds to a stratum $S'$ of $\tilde{\Xcal}$ as above. Let $\mb_\Delta \in \zdop^n$ and $c_\Delta \in \qdop$ such that $f_\Lcal(\ub)=\mb_\Delta \cdot \ub + c_\Delta$ for all $\ub \in \Delta$. Then the dual polytope $\Delta^g$ of $\Delta$ with respect to $g:=f_\Lcal$ is given by the face
$$\Delta^g:= \{\ub\}^g \cap \left(\mb_\Delta + \Delta^\bot \right)$$
of the dual polytope $\{\ub\}^g$ of the vertex $\ub$ of $\Delta$ (see \ref{setup for degree}), where $\Delta^\bot$ is the orthogonal complement of $\Delta$ in $\rdop^n$ (see \cite{Gu4}, Appendix A). Let ${\mathbb L}_\Delta$ be the linear space such that $\ub +  {\mathbb L}_\Delta$ is the affine space spanned by  $\Delta$ and let $\Lambda_\Delta := \Lambda \cap {\mathbb L}_\Delta$ be the complete lattice in  ${\mathbb L}_\Delta$ induced by $\Lambda$. 

Recall that $m_\Ccal$ is the natural number such that a suitable translate of $m_\Ccal \sigma$ is  ${\rm GL}(n,\zdop)$-isomorphic to the standard simplex $\{\ub \in \rdop_+^n \mid u_1 + \dots + u_n \leq 1 \}$ for every $\sigma \in \Ccal$. Let $\vol_\Delta$ be the Haar measure on $\ldop_\Delta$ such that $\vol_\Delta(\Delta-\ub)=1/(r!m_\Ccal)$. On $\ldop_\Delta^*$, we will use the dual measure also denoted by $\vol_\Delta$. These are the volumes from Theorem \ref{explicit formula for canonical measure}. On the other hand, we have the relative Lebesgue measure on $\ldop_\Delta \subset \rdop^n$ which is used for the Dirac measure $\delta_\Deltabar$ and which we denote now by $\vol_{\rdop^n}$.
Formula  \eqref{explicit formula} in Theorem \ref{explicit formula for canonical measure} yields
$$ \lambda_{\Deltabar}= \frac{d!}{(d-r)!} \cdot \deg_{\Lcal_1}(\overline{S'}) \cdot
\frac{\vol_\Delta\left( (\Lambda_\Delta)^L \right)}{\vol_{\rdop^n}(\Lambda_\Delta)}.$$
Here, we have used the complete lattice $(\Lambda_\Delta)^L:=\{b(\cdot, \lambda) \mid \lambda \in \Lambda_\Delta\}$ in $({\mathbb L}_\Delta)^*$ defined by the bilinear form $b$ associated to $L$ (see \ref{line bundles}).  By using \eqref{S'}, we get
$$\dim(S')=d-r, \quad \deg_{\Lcal_1}(\overline{S'}) = \binom{m}{n-r} \cdot \deg_{\Lcal_1}(\tilde{\Bcal}_1) \cdot \deg_{\Lcal_2}(\overline{S_2}).$$
By the theory of toric varieties, the degree of the toric variety $\overline{S_2}$ with respect to $\tilde{\Lcal}_2$ is given in terms of $\vol(\Delta^g)$. As in formula (36) in \cite{Gu4}, we get
$$\deg_{\Lcal_2}(\overline{S_2}) = (n-r)! \cdot \vol_{\rdop^n}(\Delta^g) \cdot \vol_{\rdop^n}(\zdop^n \cap \Delta^\bot)^{-1}$$
and hence $\deg_{\Lcal_1}(\tilde{\Bcal}_1)=\deg_{L_1}(B_1)$ yields
$$\lambda_{\overline{\Delta}} =\frac{d! \cdot m! \cdot \vol_{\rdop^n}(\Delta^g) \cdot \vol_\Delta((\Lambda_\Delta)^L) \cdot \deg_{L_1}(B_1)}{(d-r)! \cdot  (m+r-n)! \cdot \vol_{\rdop^n}(\zdop^n \cap \Delta^\bot) \cdot \vol_{\rdop^n}(\Lambda_\Delta)}.$$
\end{ex}

\appendix
\section{Building blocks}

Let $\kdop$ be  an algebraically closed field  with a non-trivial non-archimedean complete absolute value $|\phantom{a}|$ and valuation $v:=-\log |\phantom{a}|$. In the appendix, we will study  building blocks of  strongly non-degenerate strictly pluristable formal schemes  of length $l \in \ndop$ over the valuation ring $\kcirc$.

\begin{art} \label{strongly non-degenerate} \rm
Such a {\it building block} $\Ucal_l$ is recursively defined by $\Ucal_0:= \Spf(\kcirc)$ and the following conditions:

\begin{itemize}
\item[(a)] $\Ucal_l$ is an affine formal scheme over $\kcirc$ with generic fibre $U_l$.
\item[(b)] There is an \'etale morphism $\psi_l: \Ucal_l \rightarrow \Ucal_{l-1}(\nb^{(l)}, \ab^{(l)})$ over $\kcirc$ for a building block $\Ucal_{l-1}$ of length $l-1$ and $\nb^{(l)} \in (\ndop\setminus \{0\})^{p_l}, \ab^{(l)}  \in \Ocal(\Ucal_{l-1})^{p_l}$.
\item[(c)] The entries of $\ab^{(l)} =(a_1^{(l)}, \dots, a_{p_l}^{(l)})$ are units in $\Ocal(U_{l-1})$.
\item[(d)] The special fibre $\tilde{\Ucal}_l$ has a smallest stratum which maps into the smallest stratum of $\tilde{\Ucal}_{l-1}(\nb^{(l)}, \ab^{(l)})$.
\end{itemize}

Let $\Dcal= \Spf( \kcirc \langle x \rangle)$ be the formal unit disk. For $k=1, \dots,l$, we recall that $\Ucal_{k-1}(\nb^{(k)}, \ab^{(k)})$ is the closed formal subscheme of 
$$\Ucal_{k-1} \times \Dcal^{n_1^{(k)}+1} \times \dots \times \Dcal^{n_{p_k}^{(k)}+1} $$
given by the following equations:
\begin{equation} \label{pluristable equations}
x_{i0}^{(k)} \dots x_{in_i^{(k)}}^{(k)}=a_i^{(k)}  \quad (i=1, \dots, p_k)
\end{equation}
Recursively, we know that $\tilde{\Ucal}_{l-1}$ has a smallest stratum. It follows from \cite{Ber4}, Lemma 2.3, that $\tilde{\Ucal}_{l-1}(\nb^{(l)}, \ab^{(l)})$ has a smallest stratum and hence (d) makes sense. By \cite{Ber4}, Lemma 2.2, the smallest stratum of $\tilde{\Ucal}_{l}$  maps onto an open dense subset of the smallest stratum of $\tilde{\Ucal}_{l-1}(\nb^{(l)}, \ab^{(l)})$. By \cite{Ber4}, Lemma 2.10, we have an isomorphism from $\str(\tilde{\Ucal}_{l-1}(\nb^{(l)}, \ab^{(l)}))$ onto $\str(\tilde{\Ucal}_{l})$ given by taking preimages with respect to $\tilde{\psi}_l$. It is easy to see that every strongly non-degenerate strictly  pluristable formal scheme is covered by open building blocks.

For $i=1,\dots,p_l$ and $j=0,\dots ,n_i^{(l)}$, let $z_{ij}^{(l)}=\psi_l^*(x_{ij}^{(l)})$ and let $\zb^{(l)}$ be the resulting vector. Recursively, we define $\zb=(\zb^{(1)}, \dots, \zb^{(l)})$, where we use the natural pull-backs of the coordinates $\xb=(\xb^{(1)}, \dots, \xb^{(l)})$ from the definition of the building blocks $\Ucal_1, \dots, \Ucal_l$ to $\Ucal_l$. 
It will be convenient to skip all entries with index $j=0$, i.e. let
$$\hat{\xb}:=\left(x_{ij}^{(k)}\right)_{k=1,\dots,l;i=1,\dots,p_k;j=1,\dots ,n_i^{(k)}}.$$
Similarly, we define $\nb=(\nb^{(1)}, \dots, \nb^{(l)})$ and $\hat{\nb}=(\hat{\nb}^{(1)}, \dots, \hat{\nb}^{(l)})$. We set
$$\Val:U_l \longrightarrow \rdop^{|\hat{\nb}|}, \quad p \mapsto p(\hat{\xb})$$
and $\Delta_l:=\Val(U_l)$. By \cite{Ber5}, Sections 4 and 5, the map $\Val$ restricts to a homeomorphism from the skeleton $S(\Ucal_l)$ onto $\Delta_l$. It gives $S(\Ucal_l)$ a canonical piecewise linear structure and induces a canonical proper strong deformation retraction $U_l \rightarrow S(\Ucal_l)$. 

These constructions can be globalized for any strongly non-degenerate strictly pluri\-stable formal scheme $\Xcal$ over $\kcirc$. We will show in the next proposition that the building blocks induce the linear pieces of the skeleton $S(\Xcal)$. 
\end{art}

\begin{prop} \label{unit lemma}
Let $\Ucal_l$ be a strongly non-degenerate strictly pluristable building block of length $l$ with generic fibre $U_l$ and let $\Delta_l:=\val(U_l) \subset \rdop^{|\hat{\nb}|}$ as above. Then the following properties hold:
\begin{itemize}
\item[(a)] $\Delta_l$ is a polytope in $\rdop^{|\hat{\nb}|}$ defining a polytopal domain $U_{\Delta_l} :=\val^{-1}(\Delta_l)$ in $({\mathbb G}_m^{|\hat{\nb}|})^{\rm an}$ and $\Ucal_{\Delta_l}:=U_{\Delta_l}^{\rm f-sch}$(see \ref{polytopal domains}).
\item[(b)] The pull-backs $\hat{\zb}$ of the coordinates $\hat{\xb}$ define an \'etale morphism $\Ucal_l \stackrel{\psi}{\rightarrow} \Ucal_{\Delta_l}$.
\item[(c)] There is a bijective order reversing correspondence between strata $S$ of $\tilde{\Ucal}_l$ and open faces $\tau$ of $\Delta_l$. It is given by 
$$\tau= \Val(\pi^{-1}(T)), \quad S = \pi(\Val^{-1}(\tau)),$$
where $\pi:\Ucal_l \rightarrow \tilde{\Ucal}_l$ is the reduction map and $T$ is any non-empty subset of $S$. We have $\dim(\tau)=\codim(S,\tilde{\Ucal}_l)$.
\item[(d)] Let $Y$ be an irreducible component of $\tilde{\Ucal}_l$. Then there is a unique $\xi_Y \in U_l$ with $\pi(\xi_Y)$ dense in $Y$. Moreover, we have $\xi_Y \in S(\Ucal_l)$ and $\Val(\xi_Y)$ is the vertex of $\Delta_l$ corresponding to the dense open stratum of $Y$ by (c).
\item[(e)] If $f \in \Ocal(U_l)^\times$, then there is $\lambda \in \kdop^\times$, $\mu \in \Ocal(\Ucal_l)^\times$ and $\mb \in \zdop^{|\hat{\nb}|}$ with $f=\lambda\mu\hat{\zb}^\mb$. There is a unique affine function $F$ on $\Delta_l$ with $v \circ f=F\circ \Val$.
\end{itemize}
\end{prop}

\proof The proof is by induction on $l$. 
 $\ub=(\ub^{(1)}, \dots , \ub^{(l)})$ be the coordinates on $\rdop^{|\hat{\nb}|}=\rdop^{|\hat{\nb}^{(1)}|}\times \dots \times \rdop^{|\hat{\nb}^{(l)}|}$. By induction hypothesis, $a_i^{(k)}$ induces an affine function $A_i^{(k)}=A_i^{(k)}(\ub^{(1)}, \dots , \ub^{(k-1)})$ on the polytope $\Delta_{k-1}=\Val(U_{k-1} )$ in $\rdop_+^{|\hat{\nb}^{(1)}|}\times \dots \times \rdop_+^{|\hat{\nb}^{(k-1)}|}$ for $k=2,\dots, l$. Since $a_i^{(k)} \in \Ocal(\Ucal_{k-1})$, the values of $A_i^{(k)}$ are in $[0,1]$. It follows from \eqref{pluristable equations} that $\Delta_l$ is given as a subset of $\rdop_+^{|\hat{\nb}|}$ by
\begin{equation} \label{plurisimplex equations}
u_{i1}^{(k)} + \dots + u_{in_i^{(k)}}^{(k)} \leq  A_i^{(k)} \quad (k=1,\dots,l;i=1,\dots, p_k)
\end{equation}
proving (a). By induction again, we have $a_i^{(l)}=\lambda_i \mu_i \hat{\yb}^{\mb_i}$, where $\lambda_i \in \kdop^\times$, $\mu_i \in \Ocal(\Ucal_{l-1})^\times$ and $\hat{\yb}$ is the pull-back of the coordinates $(\hat{\xb}^{(1)}, \dots, \hat{\xb}^{(l-1)})$ to $\Ucal_{l-1}$. For $b_i^{(l)}:=\lambda_i \hat{\yb}^{\mb_i}$ and $c_i^{(l)}:=\lambda_i \hat{\xb}^{\mb_i}$, we get $\Ucal_{l-1}(\nb^{(l)},\ab^{(l)}) \cong \Ucal_{l-1}(\nb^{(l)},\bb^{(l)})$ and $\Ucal_{\Delta_l} \cong \Ucal_{\Delta_{l-1}}(\nb^{(l)},\cb^{(l)})$.
For the latter, we use \cite{Ber4}, Proposition 1.4. Therefore the canonical diagram
\begin{equation*} 
\begin{CD} 
\Ucal_{l-1}(\nb^{(l)},\ab^{(l)}) @> >> \Ucal_{\Delta_l} \\
@VV V    @VV V\\
\Ucal_{l-1} @> >> \Ucal_{\Delta_{l-1}}
\end{CD}
\end{equation*}
is cartesian. The bottom line is given by $\hat{\yb}$ and the induction hypothesis yields that this map is \'etale.  We conclude that the upper line is \'etale proving (b).

Note that (c) holds for any polytopal domain (see \cite{Gu4}, Proposition 4.4). It follows from \cite{Ber4}, Section 2, that
$$\str(\tilde{\Ucal}_{l-1}(\nb^{(l)}, \ab^{(l)}))\longrightarrow \str(\tilde{\Ucal}_l), \quad S' \mapsto \tilde{\psi}_l^{-1}(S')$$
is a bijective order preserving map. This proves easily (c).

By \cite{Ber4}, Proposition 1.4, we have $\Ocal(U_l)^\circ = \Ocal(\Ucal_l)$ and hence we may apply the theory of formal affinoid varieties to deduce the existence and uniqueness of $\xi_Y$ (see Section 2). Since $\Val$ maps $S(\Ucal_l)$ bijectively onto $\Delta_l$, there is  $\xi \in S(\Ucal_l)$ with $\Val(\xi)$ equal to the vertex $\Val(\xi_Y)$ of $\Delta_l$ given by the correspondence in (c). By the first paragraph on p. 332 of \cite{Ber5}, $\pi(\xi)$ is dense in $Y$ and hence $\xi=\xi_Y$ proving (d).

Let $\tilde{P}$ be a $\ktilde$-rational point in the smallest stratum of $\tilde{\Ucal}_l$. By \cite{Gu4}, $\tilde{\psi}$ induces an isomorphism $\pi^{-1}(\tilde{P}) \rightarrow \pi^{-1}(\tilde{\psi}(\tilde{P}))$ between formal fibres. This allows us to apply results for polytopal domains to the formal fibre $\pi^{-1}(\tilde{P})$. By Lemma \ref{relative polytopal domain}, we have a convergent Laurent expansion
$$f=\sum_{\mb \in \zdop^{|\hat{\nb}|}} a_\mb \hat{\zb}^\mb$$
on $\pi^{-1}(\tilde{P})$ and there is a dominant term $t:=a_\nu\hat{\zb}^\nu$ in the expansion, i.e.
$$|t(x)|>| a_\mb \hat{\zb}^\mb(x)|$$
for all $x \in \pi^{-1}(\tilde{P})$ and $ \mb \in \zdop^{|\hat{\nb}|}\setminus \{\nu\}$. Let $Y$ be an irreducible component of $\tilde{\Ucal}_l$. 
Applying (c) with $T=\{\tilde{P}\}$, we deduce that there is a sequence $x_n \in \pi^{-1}(\tilde{P})$ with $\Val(x_n) \in \relint(\Delta_l)$ converging to the vertex $\Val(\xi_Y)$ of $\Delta_l$. By compactness of $U_l$, we may assume that $x_n$ converges to some $x \in U_l$. By continuity, we have $\Val(x)=\Val(\xi_Y)$ and hence (c) again shows that $\pi(x)$ is not contained in any other irreducible component than $Y$. It is a basic fact for units in an affinoid algebra that this implies $|f(x)|=|f(\xi_Y)|$ (see \cite{Gu3}, Proposition 7.6). We conclude that
$$|f(\xi_Y)|=|f(x)|=\lim_{n \to \infty}|f(x_n)|=\lim_{n \to \infty}|t(x_n)|=|t(x)|=|t(\xi_Y)|.$$
If $Y$ ranges over all irreducible components of $\tilde{\Ucal}_l$, then the points $\xi_Y$ form the Shilov boundary of $U_l$ (see \cite{Ber}, Proposition 2.4.4). We conclude that $\mu:=ft^{-1}$ is a unit in $\Ocal(U_l)^\circ = \Ocal(\Ucal_l)$ proving (e). \qed

\begin{rem} \label{uniqueness of factors} \rm
By \eqref{plurisimplex equations}, the coordinate $u_{ij}^{(k)}$ is identically zero on $\Delta_l$ if and only if $a_{i}^{(k)} \in \Ocal(\Ucal_{k-1})^\times$. The corresponding $z_{ij}^{(k)}$ is a formal unit on $\Ucal_l$. 
We deduce easily from Proposition \ref{unit lemma}(e) that 
$$\Ocal(U_{l})^\times = \Ocal(\Ucal_{l})^\times \times \prod_{k,i,j} \left(z_{ij}^{(k)}\right)^\zdop,$$
where the indices of the basis range over $1 \leq k \leq l$, $1 \leq i \leq p_k$ with $a_{i}^{(k)} \not \in \Ocal(\Ucal_{k-1})^\times$ and $1 \leq j \leq n_i^{(k)}$.
\end{rem}

{\small Walter Gubler, Institut f\"ur Mathematik, Universit\"at T\"ubingen,
 Auf der Morgenstelle 10, D-72076 T\"ubingen, walter.gubler@uni-tuebingen.de}
\end{document}